\def\VersionDateTime{13/January/2026 (Version $1.0$)}

\documentclass[12pt, leqno
]{amsart}
\usepackage{amscd, amsmath, amssymb, amsfonts}
\usepackage{latexsym}
\usepackage{mathrsfs}

\usepackage{fullpage}
\usepackage{colonequals}

\usepackage{graphicx}
\usepackage{tikz}
\usepackage{tikz-cd}
\usetikzlibrary{calc, positioning}

\newtheorem{Theorem}{Theorem}[section]
\newtheorem{Proposition}[Theorem]{Proposition}
\newtheorem{Lemma}[Theorem]{Lemma}
\newtheorem{Corollary}[Theorem]{Corollary}
\newtheorem{Claim}{Claim}[Theorem]

\theoremstyle{definition}
\newtheorem{Definition}[Theorem]{Definition}
\newtheorem{Remark}[Theorem]{Remark}
\newtheorem{Example}[Theorem]{Example}

\numberwithin{equation}{section}

\newcommand{\TT}{{\mathbb{T}}}
\newcommand{\GG}{{\mathbb{G}}}

\newcommand{\ZZ}{{\mathbb{Z}}}

\newcommand{\RR}{{\mathbb{R}}}
\newcommand{\CC}{{\mathbb{C}}}
\newcommand{\PP}{{\mathbb{P}}}
\newcommand{\OO}{{\mathcal{O}}}

\newcommand{\Hom}{\operatorname{Hom}}
\newcommand{\Jac}{\operatorname{Jac}}

\newcommand{\GL}{\operatorname{GL}}

\newcommand{\Ker}{{\operatorname{ker}}}

\newcommand{\rest}[2]{\left.{#1}\right\vert_{{#2}}}  

\newcommand{\Spec}{{\operatorname{Spec}}}

\newcommand{\trop}{\operatorname{trop}}

\newcommand{\an}{{\operatorname{an}}}

\newcommand{\id}{\operatorname{id}}

\newcommand{\Sk}{\operatorname{Sk}}

\newcommand{\val}{\operatorname{val}}

\usepackage{mathrsfs}

\newcommand{\Proof}{{\sl Proof.}\quad}
\newcommand{\QED}{{\unskip\nobreak\hfil\penalty50\quad\null\nobreak\hfil
{$\Box$}\parfillskip0pt\finalhyphendemerits0\par\medskip}}

\newcommand{\ndot}{\raisebox{.4ex}{.}}

\begin{document}


\title{Tropical Kummer quartic surfaces}
\author{Shu Kawaguchi}
\address{Department of Mathematics, 
Kyoto University, Kyoto, 606-8502, Japan}
\email{kawaguch@math.kyoto-u.ac.jp}
\author{Kazuhiko Yamaki}
\address{Mathematical Institute, 
Tohoku University, Sendai, 980-8578, Japan}
\email{kazuhiko.yamaki.d6@tohoku.ac.jp}
\date{\VersionDateTime}


\newcommand{\Proj}{\operatorname{\operatorname{Proj}}}
\newcommand{\Prin}{\operatorname{Prin}}
\newcommand{\Rat}{\operatorname{Rat}}
\newcommand{\zero}{\operatorname{div}}
\newcommand{\Func}{\operatorname{Func}}
\newcommand{\red}{\operatorname{red}}
\newcommand{\pr}{\operatorname{pr}}
\newcommand{\can}{\operatorname{can}}
\newcommand{\relin}{\operatorname{relin}}

\begin{abstract}
We introduce tropical Kummer quartic surfaces in tropical projective $3$-space as the images of certain principally polarized tropical abelian surfaces under tropical theta functions of second order. We study some of their properties, showing that they are included in the tropicalizations of Kummer quartic surfaces defined over nonarchimdean valued fields. In the course of this work, 
we introduce the notion of a rational polyhedral orbifold and we provide faithful embeddings of tropical Kummer surfaces as such. Further, we show faithful tropicalizations of the canonical skeletons of certain Kummer surfaces over nonarchimdean valued fields. Under a suitable assumption on the base field, the canonical skeletons coincide with
the Kontsevich--Soibelman skeletons.
\end{abstract}

\maketitle

\section{Introduction}
Over $\CC$, a Kummer quartic surface is the quotient of a certain principally polarized abelian surface by the involution group $\langle -1\rangle$, embedded in projective $3$-space by the theta functions of second order. There is a long history of 
Kummer (quartic) surfaces, and they enjoy beautiful geometries. See \cite{Hudson} and \cite{Do}, for example. 

In this paper, we introduce tropical Kummer quartic surfaces in tropical projective $3$-space as the images of certain principally polarized tropical abelian surfaces under tropical theta functions of second order. We show that a tropical Kummer quartic surface is homeomorphic to the $2$-sphere $S^2$ and has $8$ distinguished points. We study its relation to a Kummer quartic surface defined over a nonarchimdean valued field (with residue characteristic $\neq 2$), and show that, under the tropicalization map, a tropical Kummer quartic surface is contained in the image of a nonarchimdean Kummer quartic surface. 

In the course of this study, we provide faithful tropicalizations of the canonical skeletons of certain Kummer surfaces over nonarchimdean valued fields via nonarchimedean theta functions of second order, which,  to our knowledge, is the first case to treat faithful tropicalizations of $K3$ surfaces except for a general theory such as \cite{GRW}. 
Over a suitable base field, these canonical skeletons coincide with the Kontsevich--Soibelman skeletons. 
We also introduce the notion of rational polyhedral orbifolds in the boundaryless case, show that 
tropical Kummer varieties are rational polyhedral orbifolds, and give faithful embeddings of tropical Kummer surfaces as such. 

In the following, we fix the notation that we use throughout this paper,  and then we explain our main results.

\subsection{Notation}
Let $\TT = \RR \cup\{+ \infty\}$ be the tropical semifield with addition $\min$ and multiplication $+$, endowed with the Euclidean topology. Let $\TT\PP^n = (\TT^n\setminus 
\{ (+ \infty, \ldots , + \infty) \}
)/\!\!\sim$ be the tropical projective $n$-space, 
where $x \sim y + (c, \ldots, c)$ for $x, y \in \TT^n\setminus 
\{ (+ \infty, \ldots , + \infty) \}$ and $c \in \RR$. 
The tropical projective $n$-space $\TT\RR^n$ has a dense open subset $(\RR^{n+1}/\!\!\sim )\cong \RR^n$, which we call the \emph{principal tropical algebraic torus} of $\TT\RR^n$.
For a $\ZZ$-module $M$, let $M_{\RR}$ denote $M \otimes \RR$. For a $\ZZ$-linear map $f\colon M\to M^\prime$ between $\ZZ$-modules, 
let $f_\RR\colon M_\RR \to M^\prime_\RR$ denote the $\RR$-linear extension.

\subsection{Tropical Kummer surfaces}
Let $M$ be a free $\ZZ$-module of rank $2$. We set $N_{\RR} := \Hom (M,\RR)$, which is equipped with a natural lattice $N := \Hom (M,\ZZ)$. Let $M^\prime$ be a lattice in $N_\RR$. Then the quotient $X \colonequals 
N_\RR/M^\prime$ is a real torus of dimension $2$ with an integral structure. Here, the integral structure means the lattice $N$ of $N_{\RR}$.
A symmetric $\RR$-bilinear form $Q\colon N_\RR \times N_\RR \to \RR$ is called a polarization on $X$ if it is positive-definite and satisfies the tropical Riemann's relation $Q(N \times M^\prime) \subseteq \ZZ$. We call $X$ a tropical abelian surface 
if it admits a polarization. A polarization $Q$ induces a $\ZZ$-linear map $\lambda\colon M^\prime \to M$ by $\lambda(u^\prime) = Q(\cdot, u^\prime)$, and it is called a principal polarization if $\lambda$ is a $\ZZ$-linear isomorphism. 

For a tropical abelian surface $X$, 
we define the tropical Kummer surface associated to $X$ to be
\[
Y = X/\langle -1\rangle. 
\]
We will introduce the notion of a rational polyhedral orbifold and show that $Y$ is a rational polyhedral orbifold. (More generally, we will show that tropical Kummer varieties in arbitrary dimensions are rational polyhedral orbifolds. See Section~\ref{subsec:trop:Kummer:orbifold:rps}.)  

Let $(X, Q)$ be a principally polarized abelian surface. Then 
one considers tropical theta functions of second order, i.e., the tropical theta functions with respect to $(2Q,0)$.  We have specified $4$ such tropical theta functions $\vartheta_{b_0}, \ldots, \vartheta_{b_3}$ that generate the $\TT$-semimodule of tropical theta functions of second order.

By the quasi-periodicity of the tropical theta functions, the map $(\vartheta_{b_0} : \cdots : \vartheta_{b_3})\colon N_{\RR} \to \TT\PP^3$ descends to a unique map $\varphi\colon  X \to \TT\PP^3$. Furthermore,
since $\vartheta_{b_0}, \ldots, \vartheta_{b_3}$ are even functions (see Equation~\eqref{eqn:theta:even}),  $\varphi\colon  X \to \TT\PP^3$ descend to a unique map 
$\psi \colon Y \to \TT\PP^3$ via the canonical projection $X \to Y$. Thus we have the following commutative diagram.
\[
\begin{tikzpicture}[auto,->]
\node (X) at (0,2) {$X$}; 
\node (P) at (4,2) {$\TT\PP^3$};
\node (Y) at (2,0) {$Y$}; 
\draw (X) --  node{$\scriptstyle \varphi 
$} (P);
\draw (X) -- node[swap]{
} (Y);
\draw (Y) -- node[swap]{$\scriptstyle \psi 
$} (P);
\end{tikzpicture}
\]
Since $\vartheta_{b_0}, \ldots \vartheta_{b_3}$ take values in $\RR$, the image of $\varphi$ and hence that of $\psi$ is contained in the principal tropical algebraic torus $(\RR^4 / \! \sim) =  \RR^3$.

One of the main results of this paper is the following; see Definition~\ref{def:faithful:embedding} for the notion of a faithful embedding. 

\begin{Theorem}[see Theorem~\ref{thm:faith:embeddings}]
\label{thm:intro:main}
Assume that $(X, Q)$ is not isomorphic to the product of 
principally polarized tropical elliptic curves. Then 
the map $\psi$ faithfully embeds $Y$ into $\TT\PP^3$, i.e., $\psi$ is a homeomorphism onto its image preserving the integral structures. 
\end{Theorem}

It will also turn out that if $(X, Q)$ is isomorphic to the product of principally polarized tropical elliptic curves, then $\psi$ can never be injective; see Remark~\ref{remark:product-type:case}. Thus the assumption in
Theorem~\ref{thm:intro:main} is essential.

Under the setting of Theorem~\ref{thm:intro:main}, we call $\psi(Y)$ the tropical Kummer quartic surface associated to $(X,Q)$. This is a terminology that is analogous to that in the classical algebraic geometry, as the image of a polarized abelian surface by the morphism given by a basis of the theta functions of second order is called a Kummer quartic surface.

The tropical Kummer surfaces can be described more explicitly.
Let $V_4$ be the Klein $4$-group which acts on the principal tropical algebraic torus $\RR^4/\!\!\sim$ by change of coordinates and let $\iota$ be the tropical Cremona involution on $\RR^4/\!\!\sim$ given by $x \mapsto -x$. Then the group~$G := V_4 \times \langle \iota \rangle$ 
acts on $\RR^4/ \!\!\sim$, and it is isomorphic to $(\ZZ/2\ZZ)\times(\ZZ/2\ZZ)\times (\ZZ/2\ZZ)$. We set 
$\tau = (\vartheta_{b_0}(0)\colon\cdots\colon \vartheta_{b_3}(0)) \in \RR^4/\!\sim$,  determined by the tropical theta constants. 
Then we will prove the following theorem.

\begin{Theorem}[see Theorem~\ref{thm:faith:embeddings2}]
\label{thm:intro:main2}
Assume that $(X, Q)$ is not isomorphic to the product of 
principally polarized tropical elliptic curves. Then the image $\psi(Y)$ is the parallelepiped whose
eight vertices are given by $\{g \cdot \tau \mid g \in G\}$. 
\end{Theorem}

\subsection{Nonarchimedean Kummer surfaces}
Let $k$ be an algebraically closed field that is complete with respect to a nontrivial 
nonarchimedean absolute value $|\ndot|_k$. 
Let $A$ be a totally degenerate abelian surface over $k$. 
Let $A^{\an}$ denote the analytification of $A$ in the sense of Berkovich. 
By the definition of totally degenerate abelian surfaces,
there exist an algebraic torus $T = \Spec(k[\chi^M])$ over $k$ with character lattice $M \cong \ZZ^2$ and a universal cover
 $p\colon T^{\an} \to A^{\an}$. 
We set $N_{\RR} := \Hom (M,\RR)$.
For $x \in T^{\an}$, we have a homomorphism $M \ni u \mapsto -\log|\chi^u(x)| \in \RR$, which defines the valuation map $\val\colon T^{\an}\to N_\RR$. This map restricts to a group homomorphism $T(k) \to N_{\RR}$.
We have $\ker(p) \subseteq T(k) \subseteq T^{\an}$ and $\val$ restricts to an isomorphism from $\ker (p)$ to a (full) lattice $M'$ in $N_{\RR}$.
Thus a commutative diagram 
\[
\begin{CD}
0 @>>>  \Ker(p) @>>> T^{\an} @>{p}>> A^{\an} @>>> 0\\
@. @VV{\cong}V @VV{\val}V @VV{\val_{A^{\an}}}V @. \\
0 @>>>  M^\prime @>>> N_\RR @>>> N_\RR/M^\prime @>>> 0\\
\end{CD}
\]
arises, in which the horizontal lines are exact.

We call the real torus $X \colonequals N_\RR/M^\prime$ with integral structure $N$ the tropicalization of $A$ or $A^{\an}$.
By the nonarchimedean Appell--Humbert theory \cite[Proposition~6.7]{BoschLutke-DAV}, a line bundle $L$ on $A$ determines a symmetric $\RR$-bilinear form $Q$ on $N_{\RR}$ such that $Q (M' \times N) \subseteq \ZZ$. Further, if $L$ is ample, then $Q$ becomes a polarization. Thus $X$ is a tropical abelian variety.
It is known that 
there is a canonical section 
$\sigma_{A^\an}\colon X \to A^{\an}$ 
of $\val_{A^{\an}}$, 
and the image $\Sigma_{A^\an}$ of 
$X$ under $\sigma_{A^\an}$ is called the canonical skeleton of $A^\an$. 

Let $\tilde{Y} \colonequals A/\langle -1 \rangle$ be the Kummer surface associated to $A$, and let $q\colon A \to \tilde{Y}$ denote the quotient morphism. A symmetric rigidified line  bundle $\tilde{L}$ on $A$ descends to a line bundle on $\tilde{L}_{\tilde{Y}}$, which we call the descent of $\tilde{L}$. Let $X = N_{\RR}/M'$ be the tropicalization of $A$ and let $Q_1$ be the polarization corresponding to $\tilde{L}_1$. We set $\Sigma_{\tilde{Y}^\an} \colonequals q^{\an}(\Sigma_{A^\an})$. We show in Proposition~\ref{prop:canonicalskeleton:kummer} that 
$\Sigma_{\tilde{Y}^\an}$ is the quotient of the canonical skeleton $\Sigma_{A^\an}$ by the action 
 of $\langle -1 \rangle$, and thus $\Sigma_{\tilde{Y}^\an}$ is homeomorphic to
 the tropical Kummer surface $Y \colonequals X/\langle -1\rangle$. 
 We call $\Sigma_{\tilde{Y}^\an}$ the {\em canonical skeleton} 
 of $\tilde{Y}^{\an}$. 
 
We will show a result on a faithful tropicalization of $\tilde{Y}^{\an}$. See Definition~\ref{def:faithful:tropicalization} for the notion of a faithful tropicalization.

\begin{Theorem}[see Theorem~\ref{thm:main1:FT}] 
\label{thm:intro:main:2}
Let $(A, \tilde{L}_1)$ be a pair of an abelian surface over $k$ and a symmetric ample line bundle $\tilde{L}_1$ on $A$. 
Assume that the polarization $Q_1$ corresponding to $\tilde{L}_1$ is a principal polarization. Further, assume that $(X,Q_1)$ is not 
a product of principally polarized tropical elliptic surfaces. Let $\tilde{L}_{\tilde{Y}}$ be the descent of $\tilde{L} \colonequals \tilde{L}_1^{\otimes 2}$. Then there exist a basis $s_0, \ldots, s_3$ of $H^0(\tilde{Y}, \tilde{L}_{\tilde{Y}})$ such that the morphism ${\varphi}_{s_0, \ldots, s_3}\colon \tilde{Y} \to \PP^3, y \mapsto (s_0(y):\cdots:s_3(y))$ faithfully tropicalizes the canonical skeleton $\Sigma_{\tilde{Y}^{\an}}$. 
\end{Theorem}

We remark that because the assumption in Theorem~\ref{thm:intro:main} is essential, the condition for $(A, \tilde{L}_1)$ in Theorem~\ref{thm:intro:main:2} is also essential. 

We also consider a lifting of tropical Kummer surface. Let 
\[
\trop\colon (\PP^3)^{\an} \ni x \mapsto (-\log|X_0(x)|:\cdots:-\log|X_3(x)|) \in \TT\PP^3
\]
be the tropicalization map, where $X_0, \ldots, X_3$ are the homogeneous coordinate functions on~$\PP^3$. 

\begin{Theorem}[see Theorem~\ref{thm:FTmain2}]
\label{thm:intro:lifting}
Let $\psi (Y)$ be a tropical Kummer quartic surface, i.e., that as in Theorem~\ref{thm:intro:main}. Then there exists a Kummer quartic surface $\tilde{Z}$ in $\PP^3_k$ over some nonarchimedean valued field $k$ such that 
$\psi(Y)$ is contained in the image $\trop(\tilde{Z}^{\an})$. 
\end{Theorem}

Since $\trop(\tilde{Z}^{\an})$ is a tropical quartic surface in $\TT\PP^3$,
the above theorem indicates that the tropical Kummer quartic surface is a part of a tropical quartic surface.

Let $A_0$ be an abelian surface over $\CC (\!( t )\!)$ and consider the Kummer surface $Y_0 := A_0 / \langle -1 \rangle$. Let $\tilde{Y}_0^\prime$ be the minimal resolution of $(\tilde{Y}_0)^{\an}$. In \cite{KS}, Konstevich and Soibelman define a certain subset of $\tilde{Y}_0^\prime$. It is called the Kontsevich--Soibelman skeleton in \cite{MN} and is of importance in the context of mirror symmetry. 
We will prove in Proposition~\ref{thm:canonical-essential} that the Kontsevich--Soibelman skeleton coincides with the canonical skeleton 
of $\tilde{Y}_0$. Then our main theorem on faithful tropicalization will lead us to a faithful tropicalization result for 
the  Kontsevich--Soibelman skeleton.
In the proof of Proposition~\ref{thm:canonical-essential}, we will use a strict Kulikov model of $\tilde{Y}_0^\prime$ together with the results in K\"unnemann \cite{Ku}, Musta\c{t}\u{a}--Nicaise \cite{MN}, and Overkamp \cite{Overkamp}

\subsection{Organization of this paper}
In Section~\ref{sec:trop:ab:theta}, we recall some basic facts on tropical abelian varieties and tropical theta functions. In Section~\ref{sec:orbifold:rps:ft}, 
after recalling the notion of a rational polyhedral space, we introduce 
the notion of a rational polyhedral orbifold and a faithful embedding. 
In Section~\ref{sec:fe:trop:Kum:surf}, after seeing that tropical Kummer varieties are 
rational polyhedral orbifolds, we state our result on 
faithful embeddings of tropical Kummer surfaces. In Section~\ref{section:proof:(1)}, we give a proof of Theorem~\ref{thm:intro:main}. In Section~\ref{sec:tropical:Kummer:quartic}, we show that 
tropical Kummer quartic surfaces are parallelepipeds (cf. Theorem~\ref{thm:intro:main2}). 

From Section~\ref{sec:uniformization}, 
we study interplay between tropical Kummer surfaces and their nonarchimedean counterparts. In Section~\ref{sec:uniformization}, we recall the uniformization theory for totally degenerate abelian varieties with a focus on the descent theory of line bundles and nonarchimedean theta functions. In Section~\ref{sec:faithful:nonarch:Kummer}, we define the canonical skeleton of a Kummer variety in arbitrary dimension, and then we show faithful tropicalizations of 
the canonical skeletons of Kummer surfaces. We also show that 
a tropical Kummer quartic surface in $\TT\PP^3$ is included in the image 
of the analytification of a nonarchimedean Kummer quartic surface in $\PP^3$ under 
the tropicalization map. In Section~\ref{can:skeleton:essential:skeleton}, we remark on the relation between the canonical skeleton and the Kontsevich--Soibelman skeleton of a Kummer surface and prove a faithful tropicalization for the Kontsevich--Soibelman skeleton. 

\bigskip
{\it Acknowledgment.}\quad 
The authors would like to thank Professors Klaus K\"unnemann and Walter Gubler for helpful discussions and kind hospitality while the authors visited the University of Regensburg in September, 2024. SK is partially supported by KAKENHI 23K03041 and 25H00587, and YK is 
partially supported by KAKENHI 23K03046 and 24K00519.

\section{Tropical abelian varieties and tropical theta functions}
\label{sec:trop:ab:theta}
We recall some basic facts on tropical abelian varieties and tropical theta functions. 
For details, we refer to \cite{MZ, FRSS, Sumi, GS, KY}. As is mentioned in \cite{FRSS}, while the name ``tropical abelian variety'' seems to have been used first in \cite{MZ}, tropical abelian varieties have been implicit in the literature for years as in Mumford \cite[\S6]{Mu}, Faltings--Chai \cite[Chap.~VI, \S1]{FC}, and Alexeev--Nakamura \cite{AN}, for example. 

\subsection{Tropical abelian varieties}
\label{subsec:trop:abelian:var}

Let $M$ be a free $\ZZ$-module of finite rank $g$. We set $N \colonequals \mathrm{Hom} (M , \ZZ)$ and $N_{\RR} \colonequals \mathrm{Hom} (M , \RR) = N \otimes \RR$. We naturally regard $N \subseteq N_\RR$ and call $N$ the \emph{integral structure} of $N_{\RR}$. Let $\langle \ndot, \ndot \rangle\colon  M \times N_{\RR} \to \RR$ be the canonical pairing. 
Let $M'$ be a lattice in $N_{\RR}$. Here, we mean by a lattice a full lattice in the vector space. Then 
the quotient $X = N_\RR/M^\prime$ is a {\em real torus with an integral structure} of dimension $g$, where the integral structure means the lattice $N$ of $N_\RR$. A real torus with an integral structure is naturally a rational polyhedral space; see Section~\ref{subsec:rps} for the definition of rational polyhedral spaces.

Let $X_1 = N_{1\RR}/M_1^\prime$ and $X_2 = N_{2\RR}/M_2^\prime$ be real tori with integral structures. 
We call a map $f\colon X_1 \to X_2$ a {\em homomorphism of real tori with integral structures} if there is a $\ZZ$-linear homomorphism $h\colon N_1 \to N_2$ such that $h_\RR(M_1^\prime) = M_2^\prime$ such that $f$ is induced by $h_\RR$, where $h_{\RR}$ denotes the $\RR$-linear extension of $h$. We note that such an $h$ is unique. We define the notion of isomorphisms of real tori with integral structures in an obvious way.
If there exists an isomorphism $X_1 \to X_2$ of real tori with integral structures, we say that $X_1$ and $X_2$ are {\em isomorphic as real tori with integral structures}. 
By definition, the isomorphisms of real tori with integral structures are isomorphism of groups. 

A symmetric $\RR$-bilinear form $Q\colon N_\RR \times N_\RR \to \RR$ is said to be \emph{integral} over $M' \times N$ if $Q(M^\prime \times N) \subseteq \ZZ$. The condition $Q(M' \times N) \subseteq \ZZ$ is called a \emph{tropical Riemann's relation}. 
For such a $Q$ that satisfies tropical Riemann's relation, 
we define a $\ZZ$-linear map $\lambda\colon M^\prime \to M$ by
\begin{equation}
\label{eqn:lambda}
\lambda(u^\prime) = Q(u^\prime, \ndot) \in \Hom(N, \ZZ) = M.
\end{equation} 
We say that such a $\ZZ$-linear map $\lambda$ \emph{corresponds} to $Q$.

We call a symmetric $\RR$-bilinear form $Q\colon N_\RR \times N_\RR \to \RR$ that satisfies the tropical Riemann's relation a \emph{polarization} on $X= N_{\RR}/M'$ if it is positive-definite. 

\begin{Definition}[Tropical abelian variety, torpical elliptic curve, torpical abelian surface]
\label{def:trop:ab}
A {\em tropical abelian variety} is a real torus with an integral structure that admits a polarization.
We say that two tropical abelian varieties are {\em isomorphic} if 
they are isomorphic as real tori with integral structures. 
A tropical abelian variety of dimension $1$ (resp. $2$) is called a {\em tropical elliptic curve} (resp. {\em a tropical abelian surface}). 
\end{Definition}

Let $Q\colon N_\RR \times N_\RR \to \RR$ be a polarization on a tropical abelian variety $X \colonequals N_\RR/M^\prime$ of dimension $g$. 
Since $Q$ is nondegenerate, the corresponding $\ZZ$-linear map $\lambda\colon  M' \to M$ given by (\ref{eqn:lambda}) is injective. Let $(d_1, d_2, \ldots, d_g) \in (\ZZ_{\geq 1})^g$ with 
$d_1 \mid d_2 \mid \cdots \mid d_g$ be the elementary divisors of $\lambda$. This $(d_1, d_2, \ldots, d_g)$ is called 
the {\em type} of $Q$. A polarization of type $(1, 1, \ldots, 1)$ is called a {\em principal polarization}. 

\begin{Definition}[(Principally) polarized tropical abelian variety]
A {\em polarized tropical abelian variety} is a pair $(X, Q)$ of a tropical abelian variety $X$ and a polarization $Q$ on $X$. 
If $Q$ is a principal polarization, then the pair $(X, Q)$ is called a {\em principally polarized tropical abelian variety}. 
\end{Definition}

Let $(X_1, Q_1)$ and $(X_2, Q_2)$ be polarized tropical abelian varieties. We write $X_1 = N_{1\RR}/M_1^\prime$ and $X_2 = N_{2\RR}/M_2^\prime$. We say that $(X_1, Q_1)$ and $(X_2, Q_2)$ are {\em isomorphic as polarized tropical abelian varieties} if there is an isomorphism $f\colon X_1 \to X_2$ of 
tropical abelian varieties such that if $h\colon N_1 \to N_2$ is 
the $\ZZ$-linear isomorphism that induces $f$, then 
$Q_1(x, y) = Q_2(h_\RR(x), h_\RR(y))$ for any $x, y \in N_{1\RR}$. 

\begin{Definition}[products]
\label{def:products}
Let $X_1 = N_{1\RR}/M_1^\prime$ and $X_2 = N_{2\RR}/M_2^\prime$ be 
real tori with integral structures. We note that $M'_1 \times M'_2$ is naturally a lattice in $(N_1 \times N_2)_{\RR}$.
\begin{enumerate}
\item
We call the real torus
$
X_1 \times X_2 \colonequals (N_1 \times N_2)_\RR/(M_1^\prime \times M_2^\prime) 
$
with an integral structure $N_1 \times N_2$ the
the {\em product} of $X_1$ and~$X_2$. 
\item
Let~$Q_1$ and $Q_2$ be polarizations on $X_1$ and $X_2$, respectively. We set 
\begin{align}
\label{eqn:Q1:Q2}
Q_1 \times Q_2\colon 
& (N_1 \times N_2)_\RR \times (N_1 \times N_2)_\RR \to \RR, 
\\
\notag
& \left((v_1, v_2), (v_1^\prime, v_2^\prime) \right) \mapsto 
Q_1(v_1, v_1^\prime) + Q_2(v_2, v_2^\prime). 
\end{align}
Then $Q_1 \times Q_2$ is a positive-definite, symmetric $\RR$-bilinear form and satisfies $(Q_1 \times Q_2)(M_1^\prime \times M_2^\prime, N_1 \times N_2) \subseteq \ZZ$, and thus $Q_1 \times Q_2$ is a polarization on $X_1 \times X_2$. We call this polarization the \emph{product} of $Q_1$ and $Q_2$.
\item
Let $(X_1, Q_1)$ and $(X_2, Q_2)$ be polarized tropical abelian varieties. 
We call $(X_1 \times X_2, Q_1 \times Q_2)$ the {\em product} of polarized tropical abelian varieties $(X_1, Q_1)$ and $(X_2, Q_2)$. 
\end{enumerate} 
\end{Definition}

\begin{Remark} \label{remark:product-PP}
Under the setting of Definition~\ref{def:products}~(3), 
for $i = 1, 2$, let $\lambda_i\colon M_i^\prime \to M$ be the $\ZZ$-linear map  corresponding to $Q_i$ by \eqref{eqn:lambda}. Then the $\ZZ$-linear map
\[
\lambda_1 \times \lambda_2\colon M_1^\prime \times M_2^\prime \to M_1 \times M_2
, \quad (u_1^\prime, u_2^\prime) \mapsto \left(\lambda_1(u_1^\prime), 
\lambda_2(u_2^\prime)\right). 
\]
corresponds to $Q_1 \times Q_2$.
It follows, in particular, that if $(X_1 \times X_2, Q_1 \times Q_2)$ is a principally polarized tropical abelian variety if and only if so are $(X_1, Q_1)$ and $(X_2, Q_2)$.
\end{Remark}

\subsection{Tropical theta functions}

In this subsection, we recall the notion of tropical theta functions and their basic properties. We refer to \cite{MZ, Sumi} for details.

A function $h\colon U \to \RR$ on an open subset $U$ of $N_{\RR}$ is called a {\em tropical function} if there exist an $m \in M$ and a $b \in \RR$ such that $h(u) = \langle m, u \rangle$ for any $u \in U$. 
Let $f\colon N_{\RR} \to \TT$ be a $\TT$-valued function on $N_{\RR}$. We say that $f$ is \emph{regular} at $v \in N_{\RR}$ if there exists an open neighborhood $U$ of $v$ in $N_{\RR}$ such that $\rest{f}{U}$ is constantly equal to $+ \infty$ or there exist a finite number of tropical functions $h_1 , \ldots , h_m$ on $U$ such that $f= \min_{i=1 , \ldots , m} h_i$.

Let $(Q , \ell)$ be a pair of a symmetric $\RR$-bilinear form $Q\colon  N_{\RR} \times N_{\RR} \to \RR$ that is integral over $M' \times N$ and an $\RR$-linear form $\ell\colon  N_{\RR} \to \RR$. We say that a $\TT$-valued function $f(x)$ on $x \in N_{\RR}$ is \emph{quasi-periodic} with respect to $(Q,\ell)$ if for any $u' \in M'$, we have
\begin{equation}
\label{eqn:quasi:periodicity:1}
f (x +  u^\prime)
= f(x) -  Q(x,  u^\prime) - \frac{1}{2} Q(u^\prime, u^\prime)
  + \ell(u^\prime).
\end{equation}
A \emph{tropical theta function} with respect to $(Q,\ell)$ means a regular function $\vartheta\colon  N_{\RR} \to \TT$ that is quasi-periodic with respect to $(Q,\ell)$.
One sees that the set of tropical theta functions with respect to $(Q,\ell)$ has a natural structure of $\TT$-semimodule.

Assume that $Q$ is a polarization. Then the homomorphism $\lambda\colon M^\prime \to M$ corresponding to $Q$ is injective. Since $N_{\RR} = M'_{\RR}$, $\lambda$ induces an $\RR$-linear extension $\lambda_{\RR}\colon  N_{\RR} \to M_{\RR}$, which is an isomorphism.  Since $Q$ is nondegenerate,  there exists a unique $r \in N_\RR$ such that $\ell(\cdot) = Q(\cdot, r)$. 
For any $b \in M$, we define a function
$\vartheta_b\colon N_\RR \to \RR$ by 
\begin{equation}
\label{eqn:def:theta:b}
\vartheta_b(x) 
= \min_{u^\prime \in M^\prime} 
\left(
Q\!\left(x,\, \lambda_\RR^{-1}(b) + u^{\prime} \right) + \frac{1}{2} Q\!\left(u^\prime + \lambda_\RR^{-1}(b) - r,\, u^\prime + \lambda_\RR^{-1}(b) - r\right)
\right), 
\end{equation}
where we remark that since $Q$ is positive-definite, the minimum on the right-hand side of \eqref{eqn:def:theta:b} exists.
It is straightforward to see that $\vartheta_b$ is quasi-periodic and thus a tropical theta function with respect to $(Q,\ell)$.

\begin{Remark} \label{remark:independent:representatives}
Let $\mathfrak{B} \subseteq M$ be a complete system of 
representatives of  $M/\lambda(M^\prime)$, which is a finite set. One sees that
$\vartheta_{b+\lambda(u^\prime)}(x) = 
\vartheta_{b}(x)$ 
for any $u^\prime \in M^\prime$, and thus $\{\vartheta_b\}_{b \in \mathfrak{B}}$ is independent of the choice of a representative $\mathfrak{B} \subseteq M$ of $M' / \lambda (M)$. 
\end{Remark}

\begin{Lemma} \label{lemma:even:theta}
Let $Q$ be a principal polarization on $X = N_{\RR}/M'$ and let $\lambda\colon  M' \to M$ be the corresponding homomorphism. Let $\mathfrak{B} \subseteq M$ be a complete system of representatives of $M/2\lambda (M')$. We take any $b \in \mathfrak{B}$, and let $\vartheta_b$ be the tropical theta function with respect to $(2Q,0)$ given by (\ref{eqn:def:theta:b}) for $(2Q,0)$ in place of $(Q,\ell)$. Then $\vartheta_b$ is an even function, i.e., $\vartheta_b (-x) = \vartheta_b(x)$ for any $x \in N_{\RR}$.
\end{Lemma} 

\Proof
We have
\begin{align*}
\vartheta_{b}(x) 
& = \min_{u^\prime \in M^\prime} 
\left(
(2Q)\left(x,\, u^\prime + (2\lambda)_{\RR}^{-1}(b)\right)
+ Q\left(u^\prime + (2\lambda)_{\RR}^{-1}(b),\, u^\prime + (2\lambda)_{\RR}^{-1}(b)\right)
\right) 
\\
& = 
\min_{u^\prime \in M^\prime} 
Q\left(x + u^\prime + (2\lambda)_{\RR}^{-1}(b),\, x + u^\prime + (2\lambda)_{\RR}^{-1}(b)
\right) - 
Q(x, x).
\end{align*}
Since $Q$ is a principal polarization, $\lambda\colon M^\prime \to M$ is an isomorphism, and hence 
there exists a $u_0^\prime \in M^\prime$ such that 
$b = \lambda(u_0^\prime)$. Then noting $(2 \lambda)_\RR^{-1}(b) = \frac{1}{2} u_0^\prime$, we have
{\allowdisplaybreaks 
\begin{align}
\label{eqn:theta:even}
\vartheta_{b}(-x)
&= \min_{u^\prime \in M^\prime} 
Q\left(-x + u^\prime + (2\lambda)_{\RR}^{-1}(b),\, -x + u^\prime + (2\lambda)_{\RR}^{-1}(b)
\right) - 
Q(-x, -x)
\\
\notag
& 
= \min_{u^\prime \in M^\prime} 
Q\left(x - u^\prime - (2\lambda)_{\RR}^{-1}(b),\, x - u^\prime - (2\lambda)_{\RR}^{-1}(b)
\right) - 
Q(x, x)
\\
\notag
&  
= \min_{u^\prime \in M^\prime} 
Q\left(x - (u^\prime + u_0^\prime) + (2\lambda)_{\RR}^{-1}(b),\, x - (u^\prime + u_0^\prime) + (2\lambda)_{\RR}^{-1}(b)\right)
 - 
Q(x, x)
\\
\notag
&  
= \vartheta_{b}(x). 
\end{align}
}
Thus $\vartheta_{b}$ is an even function. 
\QED

\subsection{Tropical descent data} \label{subsection:TDD}

There is another equivalent way to describe tropical theta functions, which will be useful when we consider the relation between the tropical theta functions and the nonarchimdean theta functions
in Sections~\ref{sec:uniformization}--\ref{can:skeleton:essential:skeleton}. 
We explain this description following \cite{FRSS}. 
Let $M, N, N_\RR$, and $M^\prime$ be those as in the beginning of Subsection~\ref{subsec:trop:abelian:var}.

\begin{Definition}[tropical descent datum]
\label{def:trop:descent:datum}
We call a pair $(\lambda, \gamma)$ of a $\ZZ$-linear map $\lambda\colon M^\prime \to M$ and a map $\gamma\colon M^\prime \to \RR$ a {\em tropical descent datum} 
if it satisfies 
\begin{equation}
\label{eqn:trop:descent:datum}
\gamma(u_1^\prime + u_2^\prime) - \gamma(u_1^\prime) -\gamma(u_2^\prime) = 
\langle \lambda(u_2^\prime), u_1^\prime 
\rangle
\qquad (\text{for any $u_1^\prime, u_2^\prime \in M^\prime$}), 
\end{equation}
where $\langle \ndot, \ndot\rangle\colon M \times N_\RR \to \RR$ is the canonical pairing. 
\end{Definition}

Note that by \eqref{eqn:trop:descent:datum}, we have 
$\langle \lambda(u_2^\prime), u_1^\prime 
\rangle = \langle \lambda(u_1^\prime), u_2^\prime 
\rangle$ for any $u_1^\prime, u_2^\prime \in M^\prime$. 

The following lemma can be checked straightforwardly (cf. \cite[Lemma~3.1]{KY}). 

\begin{Lemma}
\label{lem:Q:ell:lambda:gamma}
There exists a natural bijection between the set of pairs $(Q, \ell)$ and that of tropical descent data $(\lambda , \gamma)$, where $Q$ is a symmetric $\RR$-bilinear form that satisfies tropical Riemann's relation and $\ell\colon  N_{\RR} \to \RR$ is an $\RR$-linear map. Here, given $Q$, the map $\lambda$ corresponds to $Q$ by \eqref{eqn:lambda}, and conversely, given $\lambda$, the corresponding $\RR$-bilinear form $Q$ is characterized by $Q(u_1^\prime, u_2^\prime) = \langle \lambda(u_1^\prime), u_2^\prime 
\rangle$. Further, $\ell$ and $\gamma$ are related by 
$\gamma(u^\prime) = \frac{1}{2} Q(u^\prime, u^\prime) - \ell(u^\prime)
$
for $u^\prime \in M^\prime$.
\end{Lemma}

We set $X = N_{\RR}/M'$. Let $(Q, \ell)$ be a pair as in the above lemma and let $(\lambda , \gamma)$ be the corresponding tropical descent datum.

\begin{Remark} \label{remark:notation:H0}
For the sake of compatibility with our interplay with nonarchimedean geometry, we write $H^0 (X , L(\lambda , \gamma))$ for the $\TT$-semimodule of tropical theta functions with respect to $(Q,\ell)$.
\end{Remark}

We remark on the relation among the data $(Q,\ell)$ as above, the tropical descent data,  and the ``tropical line bundles'' on the real torus with an integral structure. We refer to \cite{MZ, Sumi, GS} for details. 
Let $(Q,\ell)$ be a pair as above and let $(\lambda , \gamma)$ be the corresponding tropical descent datum. Then $(\lambda , \gamma)$ induces an $M'$-action on $N_{\RR} \times \TT$ given by, for any $u' \in M'$ and for any $(x,t) \in N_{\RR} \times \TT$,
\[
(x ,t ) \mapsto (x + u' , t - \langle \lambda (u') , x \rangle - \gamma (u'))
.
\]
The natural projection $N_{\RR} \times \TT \to N_{\RR}$ is $M'$-equivariant, and passing through the quotient by this $M'$-action, we obtain a tropical line bundle $L(\lambda,\gamma)$. We can define the notion of regular global sections of $L(\lambda ,\gamma)$, and the regular global sections naturally form a $\TT$-semimodule. Noting the quasi-periodicity (\ref{eqn:quasi:periodicity:1}) and Lemma~\ref{lem:Q:ell:lambda:gamma}, we see that this $\TT$-semimodule is canonically identified with the $\TT$-semimodule of tropical theta functions with respect to $(Q,\ell)$. Further, it is known as a tropical Appell--Humbert theorem that any tropical line bundle on $N_{\RR}/M'$ can be obtained as $L(\lambda,\gamma)$ for some tropical descent datum. 

\section{Rational polyhedral orbifolds and faithful embeddings}
\label{sec:orbifold:rps:ft}

In this section, we define the notion of rational polyhedral orbifolds and that of faithful embeddings. 

\subsection{Rational polyhedral spaces} 
\label{subsec:rps}
In this subsection, we recall the notion of rational polyhedral spaces. For details, we refer to \cite{MZ14, JRS, GS}. 

A {\em rational polyhedron} in $\RR^n$ is a finite intersection of subsets of the form $\{x \in \RR^n \mid [a, x]\leq b \}$, where $a \in \ZZ^n, b \in \RR$, and $[\ndot, \ndot]$ is the standard inner product on $\RR^n$. A {\em rational polyhedral set} in $\RR^n$ is a finite union of rational polyhedra in $\RR^n$. We put on a rational polyhedral set in $\RR^n$ the restriction topology of the Euclid topology of $\RR^n$. 
An {\em open rational polyhedral set}  in $\RR^n$ means an open subset of a rational polyhedral set in $\RR^n$. 

A map $F\colon \RR^m \to \RR^n$ is {\em $\ZZ$-affine} if there exist 
an $m \times n$ matrix $A$ with entries in $\ZZ$ and a $b \in \RR^n$ such that $F(x) = A x + b$ for all $x \in \RR^m$, where $x$ is regarded as a column vector. A map $\varphi\colon U \to V$, where $U \subseteq \RR^m$ and $V \subseteq \RR^n$ are open rational polyhedral sets, is called a {\em tropical morphism} if for any $p \in U$, there are an open neighborhood $U'$ of $p$ in $U$ and a $\ZZ$-affine map $F\colon \RR^m \to \RR^n$ such that $\rest{\varphi}{U'} = \rest{F}{U'}$. When $V = \RR$, a tropical morphism is nothing but a tropical function. We note that the composite of two tropical morphisms is a tropical morphism. 
A tropical morphism between open rational polyhedral spaces that has an inverse tropical morphism is called a \emph{tropical isomorphism}. An \emph{open immersion} of open rational polyhedral spaces means a tropical morphism that is an open map and gives a tropical isomorphism onto its image.

A \emph{rational polyhedral chart} on a topological space $X$ means an injective open continuous map $\alpha\colon  V \to X$ with $V$ being an open rational polyhedral set. We sometimes write $(V , \alpha)$ for a rational polyhedral chart $\alpha\colon  V \to X$.

\begin{Definition}[rational polyhedral atlases]
\label{def:rat:poly:space}
A {\em rational polyhedral atlas} for a topological space $X$ is an indexed family  $\{(V_i, \alpha_i)\}_{i \in \mathcal{I}}$ of rational polyhedral charts on $X$ such that 
setting $U_i := \alpha_i (V_i)$ for any $i \in \mathcal{I}$, we have the following:
\begin{enumerate}
\item[(i)]
the family $\{ U_i \}_{i \in \mathcal{I}}$ is an open covering of $X$;
\item[(ii)]
for any $i, j \in \mathcal{I}$, 
the transition map 
\[
\rest{\alpha_i^{-1} \circ \alpha_j}{\alpha_j^{-1}(U_i \cap U_j)}: \alpha_j^{-1}(U_i \cap U_j) \to \alpha_i^{-1}(U_i \cap U_j)
\] 
is a tropical isomorphism. 
\end{enumerate}
\end{Definition}

Two rational polyhedral atlases for a topological space $X$ are {\em equivalent} if their union forms a rational polyhedral atlas for $X$. We call an equivalence class of rational polyhedral atlases for $X$ a \emph{rational polyhedral structure} on $X$.
A {\em rational polyhedral space} means a topological space equipped with
 a rational polyhedral structure on the topological space.

Let $X$ be a rational polyhedral space. An \emph{atlas} of $X$ means a rational polyhedral atlas that belongs to the rational polyhedral structure of $X$.
A  \emph{chart of $X$} means a rational polyhedral chart $\alpha\colon  V \to X$ on the topological space $X$ such that there exists an atlas of the rational polyhedral space $X$ to which $\alpha$ belongs. 

\begin{Example} \label{example:RPS}
Let $N$ be a free $\ZZ$-module of finite rank. Then $N_{\RR}$ is naturally a rational polyhedral space. Further, for a lattice $M' \subseteq N_{\RR}$, the quotient $N_{\RR}/M'$ is naturally a rational polyhedral space. Thus a real torus with an integral structure is a rational polyhedral space.
\end{Example}

\begin{Remark}
A rational polyhedral space in Definition~\ref{def:rat:poly:space} is called 
a boundaryless rational polyhedral space in \cite{GS}. Since we only consider  boundaryless rational polyhedral spaces in this paper, for simplicity, we omit the word ``boundaryless.'' 
\end{Remark}

Let $f\colon  X \to Y$ be a map between rational polyhedral spaces. We say that $f$ is a \emph{tropical morphism at $x \in X$} if for any chart $\beta\colon  W \to Y$ of $Y$ with $f(x) \in \beta (W)$, there exists a chart $\alpha\colon  V \to X$ such that $x \in \alpha (V)$, $f(\alpha (V)) \subseteq \beta (W)$, and such that the composite
\begin{align} \label{align:tropical:morphism:def:by:charts1}
\begin{CD}
V @>{\alpha}>> \alpha (V) @>{f}>> \beta (W) @>{\beta^{-1}}>> W
\end{CD}
\end{align}
is a tropical morphism of open rational polyhedral sets. One sees that 
$f$ is a tropical morphism at $x$ if and only if for some chart $\beta\colon  W \to Y$ of $Y$ with $f(x) \in \beta (W)$, there exists a chart $\alpha\colon  V \to X$ that has the same property as above.
We call $f\colon X \to Y$ a \emph{tropical morphism} if it is a tropical morphism at any $x \in X$.
If $Y = \RR$, we call $f$ a \emph{tropical function} on $X$. The notion of \emph{tropical isomorphisms} of rational polyhedral spaces is defined in an obvious way.

We remark that a homomorphism of real tori with integral structures is a tropical morphism, and it is an isomorphism if and only if it is a tropical isomorphism.

We rephrase Condition~(ii) in Definition~\ref{def:rat:poly:space}. Let $X$ be a topological space. For 
two rational polyhedral charts $(V^\prime, \alpha^\prime)$ and $(V, \alpha)$ on $X$, we say that 
$(V^\prime, \alpha^\prime) \to (V, \alpha)$ is a {\em tropical embedding} if there exists 
an open immersion $\gamma\colon V^\prime \to V$ of open rational polyhedral sets such that
$\alpha = \alpha^\prime \circ \gamma$.
Let $\gamma\colon (V^\prime, \alpha^\prime) \hookrightarrow (V, \alpha)$ stand for a tropical embedding.  Then
a {rational polyhedral atlas} for $X$ is a collection $\{(V_i, \alpha_i)\}_{i \in \mathcal{I}}$ of charts on $X$ 
satisfying Conditions~(i) in Definition~\ref{def:rat:poly:space} and 
satisfying the following condition (ii)', where we recall that $U_i \colonequals \alpha_i(V_i)$.
\begin{enumerate}
\item[(ii)']
For any $i, j \in \mathcal{I}$
and for any $x \in U_i \cap U_j$, 
there exist 
a chart $(V^\prime, \alpha^\prime)$ on $X$
and a
tropical embedding $(V^\prime, \alpha^\prime) \hookrightarrow (V_i, \alpha_i)$ such that $x \in \alpha^\prime (V^{\prime})$.
\end{enumerate}
We take this view when we define a rational polyhedral orbifold atlas in the next subsection. 

\subsection{Rational polyhedral orbifolds}

In this subsection, we introduce the notion of rational polyhedral orbifolds. 

Let $V$ be a rational polyhedral space. A \emph{tropical action} of a group $G$ on $V$ is an action of $G$ on $V$ such that for any $g \in G$, the self-map $g\colon V \to V$ induced by $g$ is a tropical isomorphism.

Let $Y$ be a topological space.
A {\em rational polyhedral orbifold chart} 
on $Y$ is a triple $(V, G, \beta)$ given by 
a rational polyhedral space $V$, a finite group $G$ with a tropical action on $V$, and a continuous map $\beta\colon V \to Y$ that is $G$-invariant
 (i.e., 
$\beta \circ g = \beta$ for any $g \in G$) and that 
 induces an injective open continuous map $V/G \to Y$ (i.e., $\beta$ gives a homeomorphism from $V/G $ onto $\beta(V)$, which is open in $Y$).  

For two orbifold charts $(V^\prime, G^\prime, \beta^\prime)$ and $(V, G, \beta)$ on $Y$, 
we say that $(V^\prime, G^\prime, \beta^\prime) \to (V, G, \beta)$ is a {\em tropical embedding} if there exist an open immersion $\gamma\colon V^\prime \to V$ of rational polyhedral spaces
and an injective group homomorphism $\rho\colon G^\prime \to G$ such that $\beta^\prime = \beta \circ \gamma$ and $\gamma \circ g^\prime  = \rho(g^\prime) \circ \gamma$ for any $g^\prime \in G^\prime$. We write $(\gamma, \rho)\colon (V^\prime, G^\prime, \beta^\prime) \hookrightarrow (V, G, \beta)$. 

A {\em rational polyhedral orbifold atlas} for $Y$ is a collection 
$\{(V_i, G_i, \beta_i)\}_{i \in \mathcal{I}}$ of orbifold charts on $Y$ with the following properties: $\{ \beta_i(V_i) \}_{i \in \mathcal{I}}$ is an open covering of $Y$; for any $i,j \in \mathcal{I}$ and 
for any $y \in \beta_i (V_i) \cap \beta_j(V_j)$, there exist 
an orbifold chart $(V^\prime, G^\prime, \beta^\prime)$ on $Y$ and 
tropical embeddings $(V^\prime, G^\prime, \beta^\prime) \hookrightarrow (V_i, G_i, \beta_i)$ and 
$(V^\prime, G^\prime, \beta^\prime) \hookrightarrow (V_j, G_j, \beta_j)$ such that $y \in \beta' (V')$. 
Two rational polyhedral orbifold atlases for $Y$ are {\em equivalent} if their union forms a rational polyhedral orbifold atlas, and an equivalence class of rational polyhedral orbifold atlases for $Y$ is called a \emph{rational polyhedral orbifold structure} on $Y$.

\begin{Definition} [rational polyhedral orbifolds]
A {\em rational polyhedral orbifold} is a topological space $Y$ equipped with a rational polyhedral orbifold structure on $Y$. 
\end{Definition}

Let $Y$ be a rational polyhedral orbifold. An \emph{orbifold atlas} of $Y$ means a rational polyhedral orbifold atlas that belongs to the rational polyhedral orbifold structure of $Y$.
A  \emph{orbifold chart of $Y$} means a rational polyhedral orbifold chart $(V ,G , \beta)$ on the topological space $Y$ such that there exists an orbifold atlas of the rational polyhedral orbifold $Y$ to which $\alpha$ belongs.

\begin{Remark} \label{remark:quotient:orbifold}
A typical example of rational polyhedral orbifolds is given as the quotient of a rational polyhedral space by a finite group tropical action. Indeed, if a finite group $G$ tropically acts on a rational polyhedral space on $X$ and $Y=X/G$, then $Y$ has a canonical structure of rational polyhedral orbifold such that $\{ (X , G , q) \}$ is an orbifold atlas of $Y$, where $q\colon X \to Y$ is the quotient of $X$ by $G$. 
\end{Remark}

\subsection{Faithful embeddings of rational polyhedral orbifolds}
In this subsection, we define the notion of a faithful embedding of a rational polyhedral orbifold into tropical projective space. 

We recall the notion of piecewise $\ZZ$-affine maps. First, we recall the notion of finite rational polyhedral decomposition of $\RR^n$. Let $\sigma$ be a rational polyhedron in $\RR^n$ (see Section~\ref{subsec:rps}), and we write $\sigma = \{y \in \RR^n \mid [a_1, y]\leq b_1, \, \ldots, \,  [a_k, y]\leq b_k \}$ for some $a_i \in \ZZ^n, b_i \in \RR$ for $i = 1, \ldots, k$. A {\em face} of $\sigma$ is $\sigma$ itself or a subset defined by changing some of the defining inequalities into equalities. 
A \emph{facet} of $\sigma$ means a face of $\sigma$ of codimension $1$ in $\sigma$. A {\em finite rational polyhedral complex} on $\RR^n$ is a finite collection $\mathcal{C}$ of rational polyhedra in $\RR^n$ with the following properties: 
\begin{enumerate}
\item[(i)]
If $\sigma$ belongs to $\mathcal{C}$, then any face of $\sigma$ belongs to $\mathcal{C}$; 
\item[(ii)]
If $\sigma$ and $\sigma'$ belong to $\mathcal{C}$ and $\sigma\cap \sigma' \neq \emptyset$ , then 
$\sigma\cap\sigma'$ is a face of both $\sigma$ and $\sigma'$. 
\end{enumerate}
The underlying set of a finite rational polyhedral complex $\mathcal{C}$ is defined to be $|\mathcal{C}|  \colonequals \bigcup_{\sigma \in \mathcal{C}} \sigma \subseteq \RR^n$. When $|\mathcal{C}| = \RR^n$, we call $|\mathcal{C}|$ a \emph{finite rational polyhedral decomposition} of $\RR^n$.

Let $V$ be an open rational polyhedral set in $\RR^n$. 
A map $f\colon V \to \RR^m$ is said to be \emph{piecewise $\ZZ$-affine} if there exists a finite rational polyhedral complex $\mathcal{C}$ such that $V$ is an open subset of $|\mathcal{C}|$ and such that for any $\sigma \in \mathcal{C}$, $\rest{f}{\sigma \cap V}\colon \sigma \cap V \to \RR^m$ is a tropical morphism. A map $f\colon V \to \RR^m$ is said to be \emph{locally piecewise $\ZZ$-affine} if for any $p \in V$, there exists an open neighborhood $V'$ of $p$ in $V$ such that $\rest{f}{V'}\colon V' \to \RR^m$ is piecewise $\ZZ$-affine.

We recall the notion of unimodularity. We say that a $\ZZ$-linear map $\ZZ^n \to \ZZ^m$ is {\em unimodular} if it is injective and its cokernel is torsion free.
We say that an $m \times n$-matrix $A$ with entries in $\ZZ$ is unimodular if the corresponding $\ZZ$-linear map $\ZZ^n \to \ZZ^m$ is unimodular, 
and a $\ZZ$-affine map $\RR^n \to \RR^m$ is \emph{unimodular} if the coefficient matrix of the linear part is unimodular.

Let $V$ be an open rational polyhedral set in $\RR^n$.
We say that a tropical morphism $f :V \to \RR^m$ is \emph{unimodular} if for any $p \in V$, there exist an open neighborhood $V'$ of $p$ and a unimodular $\ZZ$-affine map $F\colon \RR^n \to \RR^m$ such that $\rest{f}{V^\prime} = \rest{F}{V^\prime}$. Finally, 
we say that a map $f :V \to \RR^m$ is \emph{unimodular} if for any $p \in V$, there exist an open neighborhood $V'$ of $p$ in $V$ and a finite rational polyhedral complex $\mathcal{C}$ such that $V'$ is an open subset of $|\mathcal{C}|$ with the following properties: 
For any $\sigma \in \mathcal{C}$, the restriction 
$\rest{f}{\sigma \cap V}\colon \sigma \cap V \to \RR^m$ is a unimodular tropical morphism. 
By definition, a unimodular map in this sense is necessarily locally piecewise $\ZZ$-affine. We remark that this notion does not depend on the choice of a  finite rational polyhedral complex $\mathcal{C}$ such that $V^\prime$ is an open subset of $|\mathcal{C}|$ and such that $\rest{f}{V^\prime \cap \sigma}$ is a tropical morphism for any $\sigma \in \mathcal{C}$.

Let $X$ be a rational polyhedral space. Let $f\colon X \to \RR^m$ be a map.
We say that $f$ is \emph{locally piecewise $\ZZ$-affine} (resp. \emph{unimodular}) if for any chart $\alpha : V \to X$ of $X$, the composite $f \circ \alpha\colon V \to \RR^m$ is a locally piecewise $\ZZ$-affine (resp. unimodular). It is equivalent to saying that there exists an atlas of $X$ such that for any chart $\alpha\colon V \to X$ in the atlas, $f \circ \alpha$ is piecewise $\ZZ$-affine (resp. unimodular). 

Let $Y$ be a rational polyhedral orbifold. Let $f\colon Y \to \RR^m$ be a map. We say that $f$ is \emph{locally piecewise $\ZZ$-affine} (resp. \emph{unimodular}) if there exists a rational polyhedral orbifold atlas $\{(V_i, G_i, \beta_i)\}_{i \in I}$ such that for any $i \in I$, the composite map $f \circ \beta_i\colon V_i \to \RR^m$ is piecewise $\ZZ$-affine (resp. unimodular). We easily see that if $f$ is unimodular, then for any orbifold chart $(V, G, \beta)$ of $Y$, the composite $f \circ \beta\colon V \to \RR^m$ is unimodular.
 
\begin{Remark} \label{remark:unimodular:orbifold}
Let $X$ be a rational polyhedral space on which a finite group $G$ tropically acts. Let $q\colon X \to Y := X / G$ be the quotient. We regard $Y$ as a rational polyhedral orbifold with atlas $\{ q\colon X \to Y\}$.  Let $\psi\colon Y \to \TT\PP^m$ be a map. Then $\psi$ is unimodular if and only if $\psi \circ q\colon X \to \TT\PP^m$ is unimodular.
\end{Remark}

We fix an isomorphism between $\RR^m$ and the principal tropical algebraic torus in $\TT\PP^m$ and regard $\RR^m \subseteq \TT\PP^m$ via this isomorphism.

\begin{Definition}[faithful embedding]
\label{def:faithful:embedding}
Let $Y$ be a rational polyhedral orbifold. Let $\psi\colon Y \to \TT\PP^m$ be a map. Assume that $\psi (Y) \subseteq \RR^m$. We call $\psi$ 
a {\em faithful embedding} if it is injective and the map $Y \to \RR^m$ induced by $\psi$ is unimodular. This definition does not depend on the choice of an isomorphism between $\RR^m$ and the principal tropical algebraic torus in $\TT\PP^m$.
\end{Definition}

\section{Faithful embeddings of tropical Kummer surfaces}
\label{sec:fe:trop:Kum:surf}
In this section, we first define tropical Kummer varieties as rational polyhedral orbifolds. Then we state in Theorem~\ref{thm:faith:embeddings} our result on faithful embeddings of tropical Kummer surfaces by tropical theta functions of second order. Further, we define the tropical Kummer quartic surfaces to be the images of the tropical Kummer varieties by such faithful embeddings, and we give a description of them in Theorem~\ref{thm:faith:embeddings2}. 

\subsection{Faithful embedding of tropical Kummer surfaces}
\label{subsec:trop:Kummer:orbifold:rps}

Let $X$ be a real torus with an integral structure. Let $\langle -1 \rangle$ be the cyclic group of order $2$ with generator $-1$. This acts on $X$ as a tropical action by  $x \mapsto -x$.  Let $Y := X/\langle -1 \rangle$ be the quotient with canonical surjection $p:X \to Y$.
Then, noting Example~\ref{example:RPS} and Remark~\ref{remark:quotient:orbifold}, we see that $Y$ has a natural structure of rational polyhedral orbifold.
We call this rational polyhedral orbifold a {\em tropical Kummer variety}. When $\dim (X) = 2$, we call $Y$  \emph{a tropical Kummer surface}.

Let $(X, Q)$ be a principally polarized tropical abelian surface. We write $X = N_\RR/M^\prime$ as before, where $N \cong \ZZ^2$ and $M^\prime \cong \ZZ^2$.  Let $\lambda\colon M^\prime \to M$ be the $\ZZ$-linear map corresponding to $Q$ in \eqref{eqn:lambda}. Since $Q$ is a principal polarization, $\lambda$ is an isomorphism. 

We consider the tropical theta functions with respect to $(2Q,0)$. The homomorphism $2 \lambda\colon M' \to M$ corresponds to $2Q$. Let $\mathfrak{B}$ be a complete system of representatives of $M/2\lambda (M')$.
Since $M / 2\lambda (M') \cong (\ZZ/2\ZZ)^{\oplus 2}$, we have $|\mathfrak{B}| = 4$ and write 
\begin{equation}
\label{eqn:B:dim:2}
\mathfrak{B} = \{b_0, b_1, b_2, b_3\}.
\end{equation} 
Then have tropical theta functions $\vartheta_{b_0}, \vartheta_{b_1}, \vartheta_{b_2}, \vartheta_{b_3}$ with respect to $(2Q , 0)$, i.e., those given in (\ref{eqn:def:theta:b}) for $(2Q , 0)$ in place of $(Q,\ell)$.

Let 
$
Y \colonequals X/\langle -1 \rangle
$
be the tropical Kummer surface associated to $X$ with quotient map $p\colon X \to Y$. 
By the quasi-periodicity of tropical theta functions, the map
\[
N_\RR \to \TT\PP^3, \quad
x \mapsto \left(
\vartheta_{b_0}(x): \vartheta_{b_1}(x): \vartheta_{b_2}(x): \vartheta_{b_3}(x)
\right). 
\]
descends to a map $\varphi\colon X \to \TT\PP^3$. Since, by Lemma~\ref{lemma:even:theta}, $\vartheta_{b_i}$ is even for any $i = 0 , \ldots , 3$,  $\varphi$ descends to a map
\begin{equation}
\label{eqn:def:psi}
\psi  \colon  Y \to \TT\PP^3.
\end{equation}  
These maps give rise to the following commutative diagram.
\[
\begin{tikzpicture}[auto,->]
\node (X) at (0,2) {$X$}; 
\node (P) at (4,2) {$\TT\PP^3$};
\node (Y) at (2,0) {$Y$}; 
\draw (X) --  node{$\scriptstyle \varphi$} (P);
\draw (X) -- node[swap]{$\scriptstyle p$} (Y);
\draw (Y) -- node[swap]{$\scriptstyle \psi$} (P);
\end{tikzpicture}
\]

\begin{Remark} \label{remark:well-defined:varphi:psi}
We see from Remark~\ref{remark:independent:representatives} that the $\varphi$ and hence $\psi$ are independent of the choice of the complete system $\mathfrak{B}$ of representatives of $M/2 \lambda (M')$ modulo the numbering of the coordinates on $\TT\PP^3$. 
\end{Remark}

To state the the main results of this section, let us introduce the notion of 
a principally polarized tropical abelian surface of {\em product type} (resp. being {\em irreducible}). 

\begin{Definition}
\label{def:irr}
We say that a polarized tropical abelian surface $(X, Q)$ is {\em of product type} if there exist polarized tropical elliptic curves $(E_1, Q_1)$ and $(E_2, Q_2)$ such that $(X, Q)$ is isomorphic to $(E_1 \times E_2, Q_1 \times Q_2)$ as polarized abelian surfaces. We say that $(X, Q)$ is {\em irreducible} if it is not of product type. 
\end{Definition}

\begin{Theorem}
\label{thm:faith:embeddings}
Let $(X, Q)$ be a principally polarized tropical abelian surface. Assume that $(X, Q)$ is irreducible. Then the map $\psi\colon Y \to \TT\PP^3$ 
in \eqref{eqn:def:psi} is a faithful embedding (cf. Definition~\ref{def:faithful:embedding}). 
\end{Theorem}

We will prove
Theorem~\ref{thm:faith:embeddings} in Section~\ref{section:proof:(1)}.

\begin{Remark} \label{remark:product-type:case}
In fact, the converse of Theorem~\ref{thm:faith:embeddings} also holds.
To see that, suppose that 
$X$ is product type. 
Then noting Remark~\ref{remark:product-PP}, we may assume that $X = E_1 \times E_2$ for some principally polarized tropical elliptic curves $E_1$ and $ E_2$.
Further, for $i = 1, 2$,  
we may assume that 
$E_i 
= \RR/\varpi_i \ZZ$ for some $\varpi_i > 0$ and that the principal polarization $Q_i\colon \RR \times \RR \to \RR$ on $E_i$ is given by $Q_i(x, y) = \frac{1}{\varpi_i} xy$. For $i=1,2$ and $j_i  \in \{0, 1\}$, we define $\vartheta_i\!\left[j_i\right]\colon \RR\to\RR$ by 
\[
\vartheta_i\!\left[j_i\right](x)  
= \min_{k \in \ZZ} \left(
2  x \left(k + \frac{j_i}{2}\right) + 
\varpi_i x \left(k  + \frac{j_i}{2}\right)^2 
\right).
\]
Further, we define $\vartheta\!\left[
\begin{smallmatrix} j_1 \\ j_2 \end{smallmatrix}\right]\colon \RR^2 \to \RR$ by
 $\vartheta\!\left[
\begin{smallmatrix} j_1 \\ j_2 \end{smallmatrix}\right](x_1, x_2) =
\vartheta_1\!\left[j_1\right](x_1) +  \vartheta_2\!\left[j_2\right](x_2)$. 
Then we see that $\psi\colon (E_1 \times E_2)/\langle -1\rangle \to \TT\PP^3$ in (\ref{eqn:def:psi}) is given by 
\[
\psi ([x_1 , x_2]) = \left( \vartheta\!\left[
\begin{smallmatrix} 0 \\ 0 \end{smallmatrix}\right] (x_1 , x_2): \vartheta\!\left[
\begin{smallmatrix} 0 \\ 1 \end{smallmatrix}\right] (x_1 , x_2)\colon \vartheta\!\left[
\begin{smallmatrix} 1 \\ 0 \end{smallmatrix}\right] (x_1 , x_2)\colon \vartheta\!\left[
\begin{smallmatrix} 1 \\ 1 \end{smallmatrix}\right] (x_1 , x_2) \right).
\] 
Since each $\vartheta_i\!\left[j_i\right]\colon \RR\to\RR$ is an even function, 
it follows that $\psi([x_1, -x_2]) = \psi([x_1, x_2])$ for any $[x_1, x_2] \in (E_1 \times E_2)/\langle -1\rangle$. Thus $\psi$ is not injective. 
\end{Remark}

Here, we give an informal explanation that irreducible principally polarized tropical abelian surfaces are generic in the moduli space $\mathcal{A}_2^{\mathrm{tr}}$
of principally polarized tropical abelian surfaces. To see that, we first recall the tropical Torelli map for tropical curves of genus $2$. We refer to \cite{MZ, BMV, Ch} for details. 
In general, it is known that to a tropical curve $\Gamma$, one assigns a principally polarized tropical abelian variety $(\Jac_{\Gamma} , Q_{\Gamma})$, called the tropical Jacobian of $\Gamma$. This assignment defines  a map from the moduli space of tropical curves of genus $g$ to that of principally polarized tropical abelian varieties of dimension $g$, which is called the tropical Torelli map. 
It is known that when $g=2$, the tropical Torelli map induces a surjective continuous map from the moduli spaces of bridgeless tropical curves of genus $2$ to $\mathcal{A}_2^{\mathrm{tr}}$. 
The bridgeless tropical curves of genus $2$ are divided into two topological types as follows.
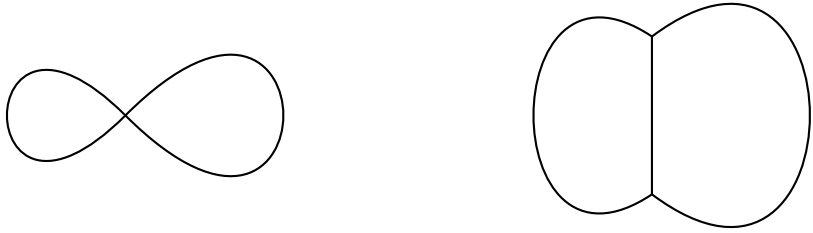
\vspace{-5ex}
\[
\begin{tikzpicture}[scale=0.7]
\draw[shift={(10,0)}] [thick] (0,0) --(0,3) ;
\draw[shift={(10,0)}][thick] (0,0) .. controls (4,-3) and (4, 6) ..  (0,3);
\draw[shift={(10,0)}][thick] (0,0) .. controls (-3,-2) and (-3, 5) ..  (0,3);
\draw[shift={(0,1.5)}] [thick]  (0,0) .. controls (-3,3) and (-3, -3) ..  (0,0);
\draw[shift={(0,1.5)}] [thick]  (0,0) .. controls (4,4) and (4, -4) ..  (0,0);
\end{tikzpicture}
\vspace{-5ex}
\]
We call such a curve a {\em petal curve}  (resp. a {\em dipole curve}) of genus $2$ if its topological type is in the left (resp. right) of the above figure.
It is known that for a bridgeless tropical curve $\Gamma$ of genus $2$, the Jacobian $(\Jac_{\Gamma} , Q_{\Gamma})$ is irreducible if and only if $\Gamma$ is a dipole curve of genus $2$.

Let $(X,Q)$ be a principally polarized tropical abelian surface of product type. 
Then noting the above facts, we find a petal curve $\Gamma_1$ of genus $2$ such that $(X, Q) \cong \left( \Jac_{\Gamma_1} , Q_{\Gamma_1} \right)$. One notes that $\Gamma_1$ is regarded as the limit of a dipole curve $\Gamma'$ of genus $2$ as $e$ contracts to a point, where $e$ is one of the three arcs that connect the two points of valence $3$ in the dipole curve $\Gamma'$. Then one sees that in $\mathcal{A}_2^{\mathrm{tr}}$, $\left( \Jac_{\Gamma'} , Q_{\Gamma'}\right)$ specializes to
$\left( \Jac_{\Gamma_1} , Q_{\Gamma_1} \right) \cong (X,Q)$ as $e$ contracts to a point. Since $\left( \Jac_{\Gamma'} , Q_{\Gamma'}\right)$ is irreducible, this observation indicates that irreducible ones are general in $\mathcal{A}_2^{\mathrm{tr}}$.

\subsection{Tropical Kummer quartic surfaces}
In this subsection, we define the tropical Kummer quartic surfaces and
give an explicit description of them. 
Under the setting of Theorem~\ref{thm:faith:embeddings}, $\psi (Y)$ is the image of an irreducible  principally polarized tropical abelian surface 
by the map given by a basis of the $\TT$-semimodule of the tropical theta functions of second order. In the classical algebraic geometry, the image of an irreducible principally polarized abelian surface by the morphism given by a basis of the vector space of theta functions of second order is called the Kummer quartic surface associated to the abelian surface. Following this classical terminology, we make the following definition.

\begin{Definition}[Tropical Kummer quartic surface]
\label{definition:TKQS}
We call $\psi(Y)$ in Theorem~\ref{thm:faith:embeddings}
a {\em tropical Kummer quartic surface}. 
\end{Definition}

We will see in Theorem~\ref{thm:FTmain2} that the tropical Kummer quartic surface is a subset of the tropicalization of a Kummer quartic surface over some nonarchimedean field. 

We can give a more explicit description of the Kummer quartic surfaces.
Let $S_4$ be the symmetric group of four letters $\{ 0 , 1 , 2, 3\}$, and let 
\[
V_4 = \{{\rm Id}, \; (0\, 1) (2\, 3), \; (0\, 2) (1\, 3), \; (0\, 3) (1\, 2)\} \subseteq S_4
\]
be the Klein $4$-group. Then $V_4$ acts on $\TT\PP^3$ by the change of coordinates. 
Let $\iota\colon \TT\PP^3 
\to \TT\PP^3$ be the tropical Cremona involution, i.e., $\iota (x_0 : \cdots : x_3) = (-x_0:\cdots:-x_3) $. 
Then the group
$V_4 \times \langle \iota \rangle \cong (\ZZ/2\ZZ) \times (\ZZ/2\ZZ) \times(\ZZ/2\ZZ)$ naturally acts
on $\TT\PP^3$. 
We set 
\begin{align*}
\tau_0 & \colonequals (\vartheta_{b_0}(0): \vartheta_{b_1}(0): \vartheta_{b_2}(0): \vartheta_{b_3}(0)), 
\\
\quad
\widetilde{\tau}_0
& \colonequals 
\iota(\tau_0) = 
(-\vartheta_{b_0}(0): -\vartheta_{b_1}(0): -\vartheta_{b_2}(0): -\vartheta_{b_3}(0)). 
\end{align*}
Further, we set 
$\tau_i \colonequals (0\, i) (j\, k) \cdot \tau_0$ and $\widetilde{\tau}_i \colonequals (0\, i) (j\, k) \cdot \widetilde{\tau}_0$ for $i = 1,2,3$, where $\{i, j, k\} = \{1, 2, 3\}$; explicitly,  
\begin{align*}
& \tau_1 = (\vartheta_{b_1}(0): \vartheta_{b_0}(0): \vartheta_{b_3}(0): \vartheta_{b_2}(0)), \quad 
\widetilde{\tau}_1  
= (-\vartheta_{b_1}(0): -\vartheta_{b_0}(0): -\vartheta_{b_3}(0): -\vartheta_{b_2}(0)),
\\
& \tau_2 = (\vartheta_{b_2}(0): \vartheta_{b_3}(0): \vartheta_{b_0}(0): \vartheta_{b_1}(0)), \quad 
\widetilde{\tau}_2
= (-\vartheta_{b_2}(0): -\vartheta_{b_3}(0): -\vartheta_{b_0}(0): -\vartheta_{b_1}(0)),
\\
& \tau_3 = (\vartheta_{b_3}(0): \vartheta_{b_2}(0): \vartheta_{b_1}(0): \vartheta_{b_0}(0)), \quad 
\widetilde{\tau}_3
= (-\vartheta_{b_3}(0): -\vartheta_{b_2}(0): -\vartheta_{b_1}(0): -\vartheta_{b_0}(0)). 
\end{align*}
We have the following description of the tropical Kummer quartic surface.

\begin{Theorem}
\label{thm:faith:embeddings2}
Let $(X, Q)$ be an irreducible principally polarized tropical abelian surface. Let 
$\psi\colon Y \to \TT\PP^3$ be the map 
in \eqref{eqn:def:psi}.
Then the tropical Kummer quartic surface $\psi(Y) \subseteq \TT\PP^3$ is a parallelepiped whose eight vertices are 
$\tau_0, \tau_1, \tau_2, \tau_3, \widetilde{\tau}_0, \widetilde{\tau}_1, \widetilde{\tau}_2, \widetilde{\tau}_3$.
\end{Theorem}

We will prove
Theorem~\ref{thm:faith:embeddings2} in Section~\ref{section:proof:(2)}.

\section{Proof of Theorem~\ref{thm:faith:embeddings}} \label{section:proof:(1)}

In this section, we prove Theorem~\ref{thm:faith:embeddings}. In the proof, the Voronoi cell of a lattice will be play crucial roles.

\subsection{Description using the Voronoi cells}
\label{subsec:reduction}

In this subsection, we give a description of a principally polarized tropical abelian surface in such a way that the polarization is the standard inner product $[\ndot,\ndot]$ on $\RR^2$. 
We consider the Voronoi cells of lattices in $\RR^2$, and we prove a criterion of irreducibility for principally polarized tropical abelian surfaces. 

Let $(X, Q)$ be a principally polarized tropical abelian surface. We let $M, N \subseteq N_{\RR}  , M'$ be as before for $X$ and write $X = N_\RR/M^\prime$.
We fix a $\ZZ$-basis $\mathsf{f}_1^\prime, \mathsf{f}_2^\prime$ of $M^\prime$. Since $Q$ is a positive-definite symmetric $\RR$-bilinear form, the $2 \times 2$ matrix
$
\begin{pmatrix}
Q(\mathsf{f}_i^\prime, \mathsf{f}_j^\prime)
\end{pmatrix}_{ i, j } 
$
is symmetric and positive-definite. It follows that there exists a unique 
positive-definite symmetric matrix $P \in \GL_2(\RR)$ such that 
$
\begin{pmatrix}
Q(\mathsf{f}_i^\prime, \mathsf{f}_j^\prime)
\end{pmatrix}_{ i, j } 
= P^2$. 

We write $P = \begin{pmatrix}\mathsf{p}_1 & \mathsf{p}_2 \end{pmatrix}$. 
Regarding $\mathsf{p}_1 ,\mathsf{p}_2 \in \RR^2$, we define a lattice $\mathsf{P} \subseteq \RR^2$ by $\mathsf{P} = \ZZ \mathsf{p}_1 + \ZZ \mathsf{p}_2$. 
We define a $\ZZ$-linear isomorphism $\alpha\colon M^\prime \to \mathsf{P}$ by $\alpha(\mathsf{f}_j^\prime) = \mathsf{p}_j$ for $j = 1, 2$. Then we have a lattice $\alpha_\RR(N) \subseteq \RR^2$, where $\alpha_{\RR}\colon N_{\RR} \to \RR^2$ is an isomorphism given as the $\RR$-linear extension of $\alpha$. Thus we have a real torus  $\RR^2/\mathsf{P}$ with an integral structure $\mathsf{N} \colonequals \alpha_\RR(N)$.  
The proposition below shows that when we consider principally polarized tropical abelian surfaces, it suffices to consider this type. 

\begin{Proposition}
\label{prop:reduction}
With the above notation, the standard inner product $[\ndot, \ndot]$ is a principal polarization on $\RR^n/\mathsf{P}$. Further, $\alpha_{\RR}\colon N_{\RR} \to \RR^2$ induces an isomorphism $(X,Q) \to \left(\RR^2/\mathsf{P}, [\ndot, \ndot]\right)$ of principally polarized tropical abelian variety.
\end{Proposition}

\Proof
Since $\alpha_{\RR} (N)$ is the integral structure on our $\RR^n/\mathsf{P}$ and since $\alpha_{\RR} (M') = \mathsf{P}$, it suffices to show that for any $x,y \in N_{\RR}$, we have $Q(x,y) = [\alpha_{\RR} (x), \alpha_{\RR} (y)]$. We only have to check when $x=\mathsf{f}_i$ and $y = \mathsf{f}_j$ for $i,j \in \{ 1,2 \}$, and indeed, it follows from the definitions of $\begin{pmatrix} \mathsf{p}_1 & \mathsf{p}_2 \end{pmatrix}$ and $\alpha_{\RR}$ that 
$Q(\mathsf{f}_i^\prime, \mathsf{f}_j^\prime) = [\mathsf{p}_i, \mathsf{p}_j] = [\alpha_{\RR} (\mathsf{f}_i^\prime), \alpha_{\RR} (\mathsf{f}_j^\prime)]$. Thus the proposition holds.
\QED

From here on to the end of this section, let $(\RR^2/\mathsf{P} , [\ndot , \ndot])$ be a principally polarized tropical abelian surface 
with an integral structure $\mathsf{N} \subseteq \RR^2$. 
To look into the structure of this principally polarized tropical abelian surface, 
we use the Voronoi cell.
Let 
\[
\mathsf{V} \colonequals \{y \in \RR^2 \mid \text{$[y, y] \leq  [y-\mathsf{p}, y-\mathsf{p}]$ for any $\mathsf{p} \in \mathsf{P}$}\}
\] 
be the {\em Voronoi cell} of $\mathsf{P}$ around the origin. One sees that $\mathsf{V}$ is a hexagon or a rectangle. A lattice vector $\mathsf{p} \in \mathsf{P}$ is a {\em Voronoi relevant vector} if the line defined by $[y, \mathsf{p}] = \frac{1}{2} [\mathsf{p}, \mathsf{p}]$ contains a facet of $\mathsf{V}$. Such a facet $\sigma$ is uniquely determined from a given Voronoi relevant vector $\mathsf{p}$, and we say that $\sigma$ is \emph{defined} by $\mathsf{p}$.
For distinct Voronoi relevant vectors $\mathsf{q}_1 , \mathsf{q}_2 \in \mathsf{P}$ with $[\mathsf{q}_1 , \mathsf{q}_2] \geq 0$, the facet defined by $\mathsf{q}_1$ and that defined by $\mathsf{q}_2$ intersects at a vertex of $\mathsf{V}$. Further, any vertex of $\mathsf{V}$ is obtained in this way. 
If $\mathsf{V}$ is a hexagon, then the Voronoi relevant vectors are given as $\{\pm \mathsf{u}, \pm \mathsf{v}, \pm \mathsf{w}\}$ with $\mathsf{u} + \mathsf{v} + \mathsf{w} = 0$, 
and any two of $\mathsf{u}$, $\mathsf{v}$, and $\mathsf{w}$ form a basis of $\mathsf{P}$. If $\mathsf{V}$ is a rectangle, then the Voronoi relevant vectors are given as $\{\pm \mathsf{u}, \pm \mathsf{v}\}$, and $\mathsf{u} , \mathsf{v}$ is an orthogonal basis of $\mathsf{P}$. 
See \cite[\S2]{CS92} for details. 

\begin{Remark} \label{remark:angle:voronoi}
From the above description, one can take linearly independent Voronoi relevant vectors $\mathsf{u}, \mathsf{v}$ such that $[\mathsf{u}, \mathsf{v}] \leq 0$. Furthermore, if $[\mathsf{u}, \mathsf{v}] =0$, then the Voronoi cell is a rectangle; if $[\mathsf{u}, \mathsf{v}] < 0$, then the Voronoi cell is a hexagon, and  there exists a Voronoi relevant vector $\mathsf{w}$ such that $\mathsf{u} + \mathsf{v} + \mathsf{w} = 0$. 
\end{Remark}

The proposition below shows that we see whether or not a principally polarized tropical abelian surface is irreducible by looking at the Voronoi cell.

\begin{Proposition}
\label{prop:equiv}
The principally polarized tropical abelian surface
$(\RR^2/\mathsf{P} , [\ndot , \ndot])$ is irreducible if and only if the Voronoi cell $\mathsf{V}$ of $\mathsf{P}$ is a hexagon. 
\end{Proposition}

\Proof
Since $\mathsf{V}$ is a rectangle or a hexagon, it suffices to show that $(\RR^2/\mathsf{P} , [\ndot , \ndot])$ is of product type if and only if $\mathsf{V}$ is a rectangle. 
We first show the ``if'' part. Suppose that $(\RR^2/\mathsf{P} , [\ndot , \ndot])$ is of product type. Then $\mathsf{P}$ has a $\ZZ$-basis $\mathsf{p}_1 , \mathsf{p}_2$ that is also an orthogonal basis of $\RR^2$. Then it is obvious that the Voronoi cell of $\mathsf{P}$ is a rectangle.

To show the converse, suppose that $\mathsf{V}$ is a rectangle. Then as we see in Remark~\ref{remark:angle:voronoi}, there exists a $\ZZ$-basis $\mathsf{q}_1 , \mathsf{q}_2$ of $\mathsf{P}$ that is also an orthogonal basis of $\RR^2$. Since $[\ndot , \ndot]$ is a principal polarization, there exists a $\ZZ$-basis $\mathsf{n}_1, \mathsf{n}_2$ of $\mathsf{N}$ such that 
$[\mathsf{n}_i,\mathsf{q}_j] = \delta_{ij}$, where $\delta_{ij}$ is the Kronecker delta. Since $\mathsf{q}_1, \mathsf{q}_2$ are orthogonal, it follows that $\mathsf{n}_i = \frac{1}{\Vert \mathsf{q}_i \Vert^2} \mathsf{q}_i$ for $i=1,2$. We set $N_i = \ZZ \frac{1}{\Vert \mathsf{q}_i\Vert^2} \mathsf{q}_i$, $M_i^\prime = \ZZ \mathsf{q}_i$, and $Q_i = \rest{[\cdot, \cdot]}{N_{i\RR} \times N_{i\RR}}$  for $i = 1, 2$. Then $(\RR^2/\mathsf{P}, [\ndot , \ndot]) \cong 
(N_{1\RR}/M_1^\prime, Q_1) \times (N_{2\RR}/M_2^\prime, Q_2)$, which proves that 
$(\RR^2/\mathsf{P} , [\ndot , \ndot])$ is of product type. 
\QED

We describe the tropical theta functions on the tropical abelian surface $\RR^2/\mathsf{P}$ with respect to $(2[\ndot , \ndot] , 0)$.
We fix Voronoi relevant vectors $\mathsf{u}, \mathsf{v}$ as in Remark~\ref{remark:angle:voronoi}. We set 
\begin{equation}
\label{eqn:q1:q2}
\mathsf{q}_1 \colonequals  \mathsf{u} 
\quad\text{and}\quad 
\mathsf{q}_2 \colonequals  \mathsf{v}. 
\end{equation}
Then $\mathsf{P} = \ZZ\, \mathsf{q}_1 + \ZZ\, \mathsf{q}_2$, as we see in the previous subsection. 
Noting that  the integral structure of $\RR^2/\mathsf{P}$ is given by $\mathsf{N}$, we set $\mathsf{N}^* \colonequals \Hom_\ZZ(\mathsf{N}, \ZZ)$. Let $\lambda\colon \mathsf{P} \to \mathsf{N}^*$ be the $\ZZ$-linear map corresponding to the polarization $[\ndot, \ndot]$. We note that $2 \lambda$ corresponds to $2[\ndot,\ndot]$. Since the polarization $[\ndot,\ndot]$ is principal, $\lambda$ is an isomorphism.
It follows that $(2 \lambda)_\RR$ restricts to a bijection from $\left\{
\frac{j_1}{2} \mathsf{q}_1 + \frac{j_2}{2} \mathsf{q}_2 \mid j_1, j_2 \in \{0, 1\}
\right\}$ to a complete system of representatives  
of $\mathsf{N}^*/2 \lambda (\mathsf{P})$. 
For $j_1, j_2 \in \{0, 1\}$, 
we define
$\vartheta\!\left[
\begin{smallmatrix}
j_1 \\ j_2 
\end{smallmatrix}
\right]\colon \RR^2 \to \RR$ by
\begin{align}
\label{eqn:appendix:xi}
\vartheta\!\left[
\begin{smallmatrix}
j_1 \\ j_2 
\end{smallmatrix}
\right](x) 
& = 
\min_{\mathsf{p} \in P} 
\left(
2 \left[x, \; \mathsf{p} + \frac{j_1}{2} \mathsf{q}_1 + \frac{j_1}{2} \mathsf{q}_2\right]
+ \left[
 \mathsf{p} + \frac{j_1}{2} \mathsf{q}_1 + \frac{j_2}{2} \mathsf{q}_2, 
 \mathsf{p} + \frac{j_1}{2} \mathsf{q}_1 + \frac{j_2}{2} \mathsf{q}_2
\right]\right)
\\
\notag
& = \min_{ (a_1, a_2) \in \ZZ^2} 
\left[
\sum_{i =1}^2
\left( a_i + \frac{j_{i}}{2} + x_i 
\right)  \mathsf{q}_i, \; 
\sum_{i =1}^2
\left( a_i + \frac{j_{i}}{2} + x_i 
\right)  \mathsf{q}_i 
\right]
- [x, x],  
\end{align}
where $x = x_1 \mathsf{q}_1 + x_2 \mathsf{q}_2 \in \RR^2$. 
With the notation in (\ref{eqn:def:theta:b}), this equals the tropical theta function $\vartheta_b (x)$ with respect to $(2 [\ndot , \ndot] ,0)$ for
$b = 2 \lambda_{\RR} \left( \frac{j_1}{2} \mathsf{q}_1 + \frac{j_2}{2} \mathsf{q}_2 \right) = \lambda (j_1 \mathsf{q}_1 + j_2 \mathsf{q}_2) \in \mathsf{N}^*$.

\begin{Remark} \label{remark:voronoi:theta}
We remark how the Voronoi cell is used to study the tropical theta functions.
Let $j_1,j_2 \in \{ 0,1\}$, $x = \sum_{i=1}^2 x_i \mathsf{q}_i \in \RR^2$, and $(c_1 , c_2) \in \ZZ^2$. Suppose that $\sum_{i =1}^2
\left( c_i + \frac{j_{i}}{2} + x_i 
\right)  \mathsf{q}_i \in \mathsf{V}$. 
Then the minimum in 
$
\min_{ (a_1, a_2) \in \ZZ^2} 
\left[
\sum_{i =1}^2
\left( a_i + \frac{j_{i}}{2} + x_i 
\right)  \mathsf{q}_i, \; 
\sum_{i =1}^2
\left( a_i + \frac{j_{i}}{2} + x_i 
\right)  \mathsf{q}_i 
\right]
$ is attained at $(a_1, a_2) = (c_1, c_2)$, 
and we have 
\[
\vartheta\!\left[
\begin{smallmatrix}
j_1 \\ j_2 
\end{smallmatrix}
\right](x) = 
\left[
\sum_{i =1}^2
\left( c_i + \frac{j_{i}}{2} + x_i 
\right)  \mathsf{q}_i, \; 
\sum_{i =1}^2
\left( c_i + \frac{j_{i}}{2} + x_i 
\right)  \mathsf{q}_i 
\right]
- [x,x].
\] 
In particular, if $x \in \mathsf{V}$, then $\vartheta\!\left[
\begin{smallmatrix}
0 \\ 0 
\end{smallmatrix}
\right](x) = 0$.
\end{Remark}

We embed $\RR^3$ into $\TT\PP^3$ by $(x_1, x_2, x_3) \mapsto (0: x_1: x_2: x_3)$. 
Then the image $\psi([x])$ is identified with 
$\left(
\vartheta\!\left[\begin{smallmatrix}1 \\ 0 \end{smallmatrix}\right](x)
- \vartheta\!\left[\begin{smallmatrix} 0 \\ 0 \end{smallmatrix}\right](x), \, 
\vartheta\!\left[\begin{smallmatrix}0 \\ 1 \end{smallmatrix}\right](x)
- \vartheta\!\left[\begin{smallmatrix} 0 \\ 0 \end{smallmatrix}\right](x), \, 
\vartheta\!\left[\begin{smallmatrix}1 \\ 1 \end{smallmatrix}\right](x)
- \vartheta\!\left[\begin{smallmatrix} 0 \\ 0 \end{smallmatrix}\right](x)
\right) \in \RR^3$. Further, by Proposition~\ref{prop:equiv}, $(\RR^2/\mathsf{P}, [\ndot, \ndot])$ being irreducible 
is the same as the Voronoi cell $\mathsf{V}$ being a hexagon. 
Thus Theorem~\ref{thm:faith:embeddings}  amounts to the proposition below.

\begin{Proposition}
\label{prop:thm:faith:embeddings:1:a}
We keep the above setting and notation.
Assume that $\mathsf{V}$ is a hexagon. We set $X \colonequals \RR^2/\mathsf{P}$ and 
let $Y \colonequals X/\langle -1\rangle$  be the tropical Kummer surface (cf. Proposition~\ref{prop:equiv}). 
Then the map 
\[
\psi \colon Y \to \RR^3, 
\quad
[x] \mapsto \left( 
\vartheta\!\left[\begin{smallmatrix}1 \\ 0 \end{smallmatrix}\right](x)
- \vartheta\!\left[\begin{smallmatrix} 0 \\ 0 \end{smallmatrix}\right](x), 
\vartheta\!\left[\begin{smallmatrix}0 \\ 1 \end{smallmatrix}\right](x)
- \vartheta\!\left[\begin{smallmatrix} 0 \\ 0 \end{smallmatrix}\right](x), 
\vartheta\!\left[\begin{smallmatrix}1 \\ 1 \end{smallmatrix}\right](x)
- \vartheta\!\left[\begin{smallmatrix} 0 \\ 0 \end{smallmatrix}\right](x)
\right)
\]
is injective and unimodular. 
\end{Proposition}

\subsection{Preliminary to the proof of Proposition~\ref{prop:thm:faith:embeddings:1:a}} \label{subsection:for:proof}
We keep the setting and the notation in the previous section. Assume that $\mathsf{V}$ is a hexagon.
We fix the notation that will be used in the proof.
For a $\mathsf{p} \in  \mathsf{P}$, let $H_{\mathsf{p}}$ denote the line in $\RR^2$ defined by $[x, \mathsf{p}] = \frac{1}{2} [\mathsf{p}, \mathsf{p}]$. 
Recall that $\mathsf{u} , \mathsf{v}$ are Voronoi relevant vectors as in Remark~\ref{remark:angle:voronoi} and that they form a basis of $\mathsf{P}$. Further, recall that 
$\{\pm\mathsf{u}, \pm\mathsf{v}, \pm\mathsf{w}\}$ is the set of Voronoi relevant vectors 
of $\mathsf{P}$, where $\mathsf{u} + \mathsf{v} + \mathsf{w} = 0$. 
Let $\mathsf{m}_1, \ldots, \mathsf{m}_6$ be the vertices of 
$\mathsf{V}$ as in the following figures; for example, $\mathsf{m}_1 = H_{\mathsf{w}} \cap H_{\mathsf{-v}}$. By symmetry, we have
$\mathsf{m}_1 + \mathsf{m}_4 = \mathsf{m}_2 + \mathsf{m}_5 = \mathsf{m}_3 + \mathsf{m}_6 = 0$. 
The three lines through the origin respectively with directions $\mathsf{u}$, $ \mathsf{v}$, and $ \mathsf{w}$ divide $\mathsf{V}$ into six quadrilaterals. 
Let $\mathsf{Q}_i$ denote the (closed) quadrilateral that contains $\mathsf{m}_i$  for $i = 1, \ldots, 6$. 

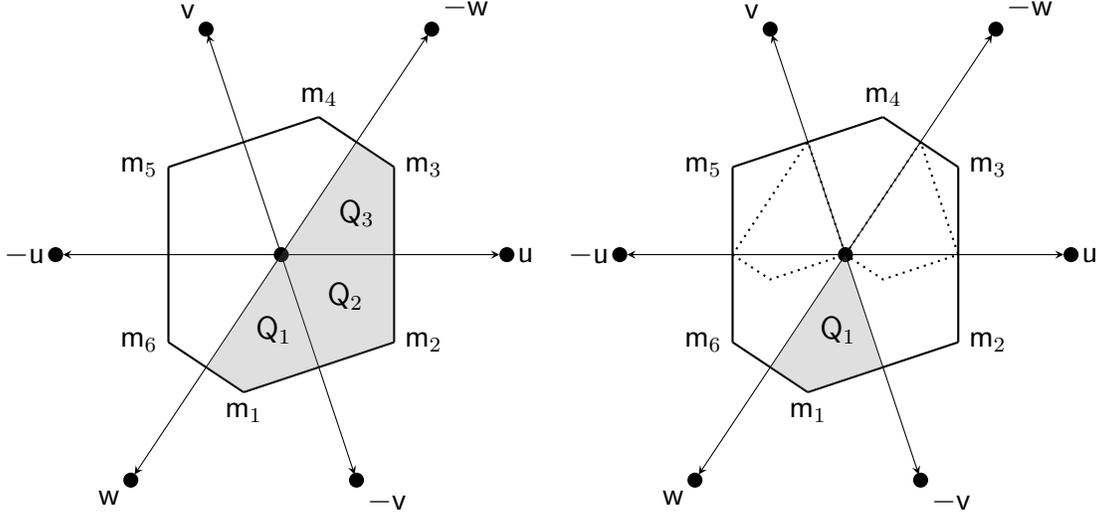
\[
\begin{tikzpicture}[scale = 0.5]
\begin{scope}
  \draw[->, >=stealth] (0,0) -- (5.8,0);
  \draw[->, >=stealth] (0,0) -- (-5.8, 0);
  \draw[->, >=stealth] (0,0) -- (3.9,5.85);
  \draw[->, >=stealth] (0,0) -- (-3.9,-5.85);
  \draw[->, >=stealth] (0,0) -- (-1.95, 5.85);
  \draw[->, >=stealth] (0,0) -- (1.95,-5.85);
  \fill (0,0) circle [radius=0.2];
  \fill (6,0) circle [radius=0.2]
  node[right]{$\mathsf{u}$}; 
  \fill (-6, 0) circle [radius=0.2]
  node[left]{$-\mathsf{u}$};
  \fill (4, 6) circle [radius=0.2]
  node[above right]{$-\mathsf{w}$};
  \fill (-4,-6) circle [radius=0.2]
  node[below left]{$\mathsf{w}$};
  \fill (-2, 6) circle [radius=0.2]
  node[above left]{$\mathsf{v}$};
  \fill (2, -6) circle [radius=0.2]
  node[below right]{$-\mathsf{v}$};
  \draw[thick] (3, -2.33)--(3, 2.33); 
  \draw[thick] (-3, -2.33)--(-3, 2.33); 
  \draw[thick] (3, 2.33)--(1, 3.66); 
  \draw[thick] (-3, 2.33)--(1, 3.66); 
  \draw[thick] (3, -2.33)--(-1, -3.66); 
  \draw[thick] (-3, -2.33)--(-1, -3.66); 
  \draw (-1, -3.66) node[below]{$\mathsf{m}_1$}; 
  \draw (3, -2.33) node[right]{$\mathsf{m}_2$}; 
  \draw (3, 2.33) node[right]{$\mathsf{m}_3$}; 
  \draw (1, 3.66) node[above]{$\mathsf{m}_4$}; 
  \draw (-3, 2.33) node[left]{$\mathsf{m}_5$}; 
  \draw (-3, -2.33) node[left]{$\mathsf{m}_6$}; 
  \fill[gray!100,nearly transparent] 
  (-2,-3)--(-1, -3.66)--(1,-3)--(0, 0)--(-2,-3); 
   \fill[gray!100,nearly transparent] 
  (3, 0)--(3, -2.33)--(1,-3)--(0, 0)--(3, 0); 
   \fill[gray!100,nearly transparent] 
  (3, 0)--(3, 2.33)--(2,3)--(0, 0)--(3, 0); 
    \draw (-0.2, -2) node{$\mathsf{Q}_1$};
  \draw (1.7, -1.1) node{$\mathsf{Q}_2$};
  \draw (2, 1.1) node{$\mathsf{Q}_3$};
    \end{scope}
  \begin{scope}[shift={($(15, 0)$)}]
  \draw[->, >=stealth] (0,0) -- (5.8,0);
  \draw[->, >=stealth] (0,0) -- (-5.8, 0);
  \draw[->, >=stealth] (0,0) -- (3.9,5.85);
  \draw[->, >=stealth] (0,0) -- (-3.9,-5.85);
  \draw[->, >=stealth] (0,0) -- (-1.95, 5.85);
  \draw[->, >=stealth] (0,0) -- (1.95,-5.85);
  \fill (0,0) circle [radius=0.2];
  \fill (6,0) circle [radius=0.2]
  node[right]{$\mathsf{u}$}; 
  \fill (-6, 0) circle [radius=0.2]
  node[left]{$-\mathsf{u}$};
  \fill (4, 6) circle [radius=0.2]
  node[above right]{$-\mathsf{w}$};
  \fill (-4,-6) circle [radius=0.2]
  node[below left]{$\mathsf{w}$};
  \fill (-2, 6) circle [radius=0.2]
  node[above left]{$\mathsf{v}$};
  \fill (2, -6) circle [radius=0.2]
  node[below right]{$-\mathsf{v}$};
  \draw[thick] (3, -2.33)--(3, 2.33); 
  \draw[thick] (-3, -2.33)--(-3, 2.33); 
  \draw[thick] (3, 2.33)--(1, 3.66); 
  \draw[thick] (-3, 2.33)--(1, 3.66); 
  \draw[thick] (3, -2.33)--(-1, -3.66); 
  \draw[thick] (-3, -2.33)--(-1, -3.66); 
  \draw (-1, -3.66) node[below]{$\mathsf{m}_1$}; 
  \draw (3, -2.33) node[right]{$\mathsf{m}_2$}; 
  \draw (3, 2.33) node[right]{$\mathsf{m}_3$}; 
  \draw (1, 3.66) node[above]{$\mathsf{m}_4$}; 
  \draw (-3, 2.33) node[left]{$\mathsf{m}_5$}; 
  \draw (-3, -2.33) node[left]{$\mathsf{m}_6$}; 
  \fill[gray!100,nearly transparent] 
  (-2,-3)--(-1, -3.66)--(1,-3)--(0, 0)--(-2,-3); 
  \draw (-0.2, -2) node{$\mathsf{Q}_1$};
  \draw[thick, dotted] (-3,0)--(-2, -0.66)--(0,0)--(-1, 3)--(-3,0); 
  \draw[thick, dotted] (0,0)--(2, 3)--(3,0)--(1, -0.66)--(0,0);
  \end{scope}
\end{tikzpicture}
\]

We prove an easy lemma on the Voronoi cell $\mathsf{V}$.

\begin{Lemma} \label{lemma:perpendicular:bisector}
Let $\mathsf{q}$ be a Voronoi relevant vector. Then the line $\{ t \mathsf{q} \mid t \in \RR \}$ is the perpendicular bisector of the facet defined by $\mathsf{q}$. Further, if $\mathsf{m}_i$ and $\mathsf{m}_j$ are the vertices of this facet, then $\mathsf{q} = \frac{1}{2} (\mathsf{m}_i + \mathsf{m}_j)$.
\end{Lemma}

\Proof
We keep the above notation.
By symmetry on $\mathsf{u} , \mathsf{v} , \mathsf{w}$, we may assume $\mathsf{q} = \mathsf{u} + \mathsf{v} = - \mathsf{w}$. Since we know that the line $\{ t \mathsf{w} \mid t \in \RR \}$ is perpendicular to the the facet defined by $\mathsf{w}$, it suffices to show that $\mathsf{m}_4 = - \mathsf{w} - \mathsf{m}_3$. Since $\mathsf{m}_3$ is an intersection point of the facet defined by $\mathsf{u}$ and that defined by $-\mathsf{w}$,
we have
$[\mathsf{m}_3 , \mathsf{u} ] = \frac{1}{2} [\mathsf{u} , \mathsf{u}]$ and
$[\mathsf{m}_3 , -\mathsf{w} ] = \frac{1}{2} [-\mathsf{w} , -\mathsf{w} ]$. 
Noting $-\mathsf{w} = \mathsf{u} + \mathsf{v}$, we see from those equalities that
$[\mathsf{m}_3 , \mathsf{v}] = \frac{1}{2} [\mathsf{v}, \mathsf{v}] + [\mathsf{u} , \mathsf{v}]$, 
and hence we compute
\begin{align*}
&[- \mathsf{w} - \mathsf{m}_3 , \mathsf{v}] = [\mathsf{u} + \mathsf{v} , \mathsf{v}] - [\mathsf{m}_3 , \mathsf{v}] = \frac{1}{2} [\mathsf{v} , \mathsf{v}]
\\
&[-\mathsf{w} - \mathsf{m}_3 , - \mathsf{w}] =  [-\mathsf{w}  , -\mathsf{w} ] - [\mathsf{m}_3 , -\mathsf{w} ]
=
\frac{1}{2} [- \mathsf{w}, - \mathsf{w}].
\end{align*}
These equality shows that $- \mathsf{w} - \mathsf{m_3}$ is on the intersection  of the lines containing the facets defined by $\mathsf{v}$ and $\mathsf{w}$. Thus $- \mathsf{w} - \mathsf{m_3} = \mathsf{m}_4$, which proves the lemma.
\QED

The following lemma will be frequently used in the proof.

\begin{Lemma} \label{claim:translate:in:voronoi}
With the above notation, we have
 $\frac{1}{2} \mathsf{v} + \mathsf{Q}_1 \subseteq \mathsf{V}$ and $\frac{1}{2}(\mathsf{-w}) + \mathsf{Q}_1 \subseteq  \mathsf{V}$ (See the figure on the right-hand side). 
\end{Lemma}

\Proof
We see that 
$0 + \frac{1}{2}\mathsf{v} = \frac{1}{2}\mathsf{v} \in \mathsf{V}$, 
$\frac{1}{2} (\mathsf{-v}) + \frac{1}{2}\mathsf{v} = 0 \in \mathsf{V}$, 
and $\frac{1}{2} \mathsf{w} + \frac{1}{2}\mathsf{v} = \frac{1}{2} (\mathsf{-u}) \in \mathsf{V}$. Further, by Lemma~\ref{lemma:perpendicular:bisector}, we have $-\mathsf{v} = \mathsf{m}_1 + \mathsf{m}_2$, and hence $\mathsf{m}_1 + \frac{1}{2} \mathsf{v} = \frac{1}{2}\mathsf{m}_1 - \frac{1}{2} \mathsf{m}_2 = \frac{1}{2} (\mathsf{m}_1 + \mathsf{m}_5) \in \mathsf{V}$.
Since $\mathsf{V}$ is convex, it follows that $\frac{1}{2} \mathsf{v} + \mathsf{Q}_1 \subseteq \mathsf{V}$. 
Similarly, we have  $\frac{1}{2}(\mathsf{-w}) + \mathsf{Q}_1 \subseteq  \mathsf{V}$. Thus the lemma holds. 
\QED

We remark on the relation between $\vartheta\!\left[\begin{smallmatrix}
0 \\ 0 
\end{smallmatrix}\right]$ and other $\vartheta\!\left[\begin{smallmatrix}
j_1 \\ j_2 
\end{smallmatrix}\right]$, which will be used in the proof of Proposition~\ref{prop:thm:faith:embeddings:1:a}.
Recall that in the definition of $\vartheta\!\left[\begin{smallmatrix}
j_1 \\ j_2 
\end{smallmatrix}\right]$, we set $\mathsf{q}_1 := \mathsf{u}$ and $\mathsf{q}_2 := \mathsf{v}$; see (\ref{eqn:q1:q2}). Note also that $\mathsf{w} = - \mathsf{u} - \mathsf{v} = - \mathsf{q}_1 - \mathsf{q}_2$.
Then we see from \eqref{eqn:appendix:xi} that for any $x \in \RR^2$ and $\mathsf{q} \in \mathsf{P}$, 
we have
\begin{equation}
\label{eqn:xi:translation:u}
\vartheta\!\left[\begin{smallmatrix}
1 \\ 0
\end{smallmatrix}\right]\left(x\right)
= 
\vartheta\!\left[\begin{smallmatrix}
0 \\ 0 
\end{smallmatrix}\right](x +  \frac{1}{2} \mathsf{u} + \mathsf{q})
+ [x, \mathsf{u} + 2 \mathsf{q}] + \frac{1}{4} [\mathsf{u} + 2 \mathsf{q}, \mathsf{u} + 2 \mathsf{q}]. 
\end{equation}
Indeed, we compute 
\begin{align*}
\text{(RHS)} 
& = \min_{\mathsf{p}\in \mathsf{P}} 
\left(
2 \left[
x +  \frac{1}{2} \mathsf{u} + \mathsf{q},  \mathsf{p}
\right] + \left[\mathsf{p}, \mathsf{p}\right]
+ [x, \mathsf{u} + 2 \mathsf{q}] + \frac{1}{4} [\mathsf{u} + 2 \mathsf{q}, \mathsf{u} + 2 \mathsf{q}]
\right)
\\
& = 
\min_{\mathsf{p}\in \mathsf{P}} 
\left(
2 \left[
x,  \frac{1}{2} \mathsf{u} + \mathsf{q} + \mathsf{p}
\right] + \left[
\frac{1}{2} \mathsf{u} + \mathsf{q} + \mathsf{p}, \frac{1}{2} \mathsf{u} + \mathsf{q} + \mathsf{p}
\right]
\right)
= \vartheta\!\left[\begin{smallmatrix}
1 \\ 0
\end{smallmatrix}\right]\left(x\right). 
\end{align*}
Similarly, for any $x \in \RR^n$ and $\mathsf{q} \in \mathsf{P}$, we get 
\begin{align}
\label{eqn:xi:translation:vw}
\begin{split}
& \vartheta\!\left[\begin{smallmatrix}
0 \\ 1
\end{smallmatrix}\right]\left(x\right)
= 
\vartheta\!\left[\begin{smallmatrix}
0 \\ 0 
\end{smallmatrix}\right](x +  \frac{1}{2} \mathsf{v} + \mathsf{q})
+ [x, \mathsf{v} + 2 \mathsf{q}] + \frac{1}{4} [\mathsf{v} +2 \mathsf{q}, \mathsf{v} + 2\mathsf{q}], 
\\
&
\vartheta\!\left[\begin{smallmatrix}
1 \\ 1
\end{smallmatrix}\right]\left(x\right)
= 
\vartheta\!\left[\begin{smallmatrix}
0 \\ 0 
\end{smallmatrix}\right](x -  \frac{1}{2} \mathsf{w} + \mathsf{q})
+ [x, -\mathsf{w} + 2 \mathsf{q}] + \frac{1}{4} [-\mathsf{w} +2 \mathsf{q}, -\mathsf{w} + 2\mathsf{q}]. 
\end{split}
\end{align}

In particular, taking $\mathsf{q} = \mathsf{0}$ or $- \mathsf{u}$ in \eqref{eqn:xi:translation:u}, 
we obtain the first line of the following equalities, and we similarly obtain the second and the third lines.  
\begin{align}
\label{eqn:xi:translation}
\begin{split}
& \vartheta\!\left[\begin{smallmatrix}
1 \\ 0
\end{smallmatrix}\right]\left(x\right)
= 
\vartheta\!\left[\begin{smallmatrix}
0 \\ 0 
\end{smallmatrix}\right](x +  \frac{1}{2} \mathsf{u})
+  [x, \mathsf{u}] + \frac{1}{4} [\mathsf{u}, \mathsf{u}] 
= 
\vartheta\!\left[\begin{smallmatrix}
0 \\ 0 
\end{smallmatrix}\right](x -  \frac{1}{2} \mathsf{u})
-  [x, \mathsf{u}] + \frac{1}{4} [\mathsf{u}, \mathsf{u}], 
\\
& \vartheta\!\left[\begin{smallmatrix}
0 \\ 1
\end{smallmatrix}\right]\left(x\right)
= 
\vartheta\!\left[\begin{smallmatrix}
0 \\ 0 
\end{smallmatrix}\right](x +  \frac{1}{2} \mathsf{v})
[x, \mathsf{v}] + \frac{1}{4} [\mathsf{v}, \mathsf{v}]
= 
\vartheta\!\left[\begin{smallmatrix}
0 \\ 0 
\end{smallmatrix}\right](x -  \frac{1}{2} \mathsf{v})
- [x, \mathsf{v}] + \frac{1}{4} [\mathsf{v}, \mathsf{v}],
\\
&
\vartheta\!\left[\begin{smallmatrix}
1 \\ 1
\end{smallmatrix}\right]\left(x\right)
= 
\vartheta\!\left[\begin{smallmatrix}
0 \\ 0 
\end{smallmatrix}\right](x -  \frac{1}{2} \mathsf{w}) 
- [x, \mathsf{w}] + \frac{1}{4} [\mathsf{w}, \mathsf{w}]
= 
\vartheta\!\left[\begin{smallmatrix}
0 \\ 0 
\end{smallmatrix}\right](x +  \frac{1}{2} \mathsf{w}) 
+ [x, \mathsf{w}] + \frac{1}{4} [\mathsf{w}, \mathsf{w}].
\end{split}
\end{align}

\subsection{Proof of the injectivity in Proposition~\ref{prop:thm:faith:embeddings:1:a}}
\label{subsec:proof:injectivity}
In this subsection, we prove the injectivity assertion in Proposition~\ref{prop:thm:faith:embeddings:1:a}. We keep the notation so far in this section.

First, we note that Proposition~\ref{prop:thm:faith:embeddings:1:a} is deduced from the implication (i) $\Longrightarrow$ (ii) of the following proposition.

\begin{Proposition}
\label{prop:xi}
Assume that $\mathsf{V}$ is a hexagon. 
Then for any $y, y^\prime \in \RR^2$,  the following are equivalent: 
\begin{enumerate}
\item[(i)]
$
\vartheta\!\left[
\begin{smallmatrix}
j_1 \\ j_2 
\end{smallmatrix}
\right](y) - \vartheta\!\left[
\begin{smallmatrix}
0  \\ 0 
\end{smallmatrix}
\right](y)
= 
\vartheta\!\left[
\begin{smallmatrix}
j_1 \\ j_2 
\end{smallmatrix}
\right](y^\prime) - \vartheta\!\left[
\begin{smallmatrix}
0  \\ 0 
\end{smallmatrix}
\right](y^\prime)$ 
for any $j_1, j_2 \in \{0, 1\}$; 
\item[(ii)]
there exists a $\mathsf{p} \in \mathsf{P}$ such that 
$y^\prime = y + \mathsf{p}$ or $y^\prime = - y + \mathsf{p}$. 
\end{enumerate}
\end{Proposition}

It is easy to show the direction (ii) $\Longrightarrow$ (i). We give a proof, though we will not need this direction. 
For any $j_1, j_2 \in \{0, 1\}$, it follows from the quasi-periodicity  \eqref{eqn:quasi:periodicity:1} that 
the function
$\vartheta\!\left[
\begin{smallmatrix}
j_1 \\ j_2 
\end{smallmatrix}
\right] -
\vartheta\!\left[ 
\begin{smallmatrix}
0 \\ 0 
\end{smallmatrix}
\right]
$
is invariant under the translation by any element in $\mathsf{P}$, and by Lemma~\ref{lemma:even:theta},
$\vartheta\!\left[
\begin{smallmatrix}
j_1 \\ j_2 
\end{smallmatrix}
\right] -
\vartheta\!\left[ 
\begin{smallmatrix}
0 \\ 0 
\end{smallmatrix}
\right]
$
is an even function.
Thus (ii) implies (i). 

\medskip
{\sl Proof of the direction (i) $\Longrightarrow$ (ii) in Proposition~\ref{prop:xi}.\ }
We take any $y,y^\prime \in \RR^2$. Assume that (i) holds. Since $\mathsf{V} + \mathsf{P} = \RR^2$, we may and do assume that $y , y' \in \mathsf{V}$. 
Then 
$
\vartheta\!\left[
\begin{smallmatrix}
0 \\ 0 
\end{smallmatrix}
\right](y)
= 0$ and $
\vartheta\!\left[
\begin{smallmatrix}
0 \\ 0 
\end{smallmatrix}
\right](y^\prime)
= 0$ 
by Remark~\ref{remark:voronoi:theta}.
It follows that 
Condition (i) is nothing but the condition that
\begin{align}
\label{align:(i)}
\vartheta\!\left[
\begin{smallmatrix}
j_1 \\ j_2 
\end{smallmatrix}
\right](y) 
= 
\vartheta\!\left[
\begin{smallmatrix}
j_1 \\ j_2 
\end{smallmatrix}
\right](y^\prime) 
\end{align} 
for any $j_1, j_2 \in \{0, 1\}$.

{\bf Case 1.}\quad Suppose that $y, y^\prime \in \mathsf{Q}_i$ for some $i = 1 , \ldots , 6$. Without loss of generality, we assume that $i=1$, i.e., $y, y^\prime \in \mathsf{Q}_1$.  
By Lemma~\ref{claim:translate:in:voronoi}, we have $\frac{1}{2} \mathsf{v} + y  \in \frac{1}{2} \mathsf{v} + \mathsf{Q}_1 \subseteq \mathsf{V}$, and the same holds for $y^\prime$ in place of $y$. Hence by Remark~\ref{remark:voronoi:theta}, 
$\vartheta\!\left[
\begin{smallmatrix}
0 \\ 0 \end{smallmatrix}
\right](\frac{1}{2} \mathsf{v} + y) = \vartheta\!\left[
\begin{smallmatrix}
0 \\ 0 \end{smallmatrix}
\right](\frac{1}{2} \mathsf{v} + y^\prime) =0$. 
It follows from \eqref{eqn:xi:translation} that 
\begin{align}
\label{eqn:appendix:01}
\vartheta\!\left[
\begin{smallmatrix}
0 \\ 1 
\end{smallmatrix}
\right](y)
=  \left[y, \mathsf{v}\right] + \frac{1}{4}[\mathsf{v}, \mathsf{v}], 
\qquad 
\vartheta\!\left[
\begin{smallmatrix}
0 \\ 1 
\end{smallmatrix}
\right](y^\prime) 
=  \left[y^\prime, \mathsf{v}\right] + \frac{1}{4}[\mathsf{v}, \mathsf{v}].  
\end{align}
The same argument using Lemma~\ref{claim:translate:in:voronoi} and Remark~\ref{remark:voronoi:theta} together with \eqref{eqn:xi:translation} shows that
\begin{equation}
\label{eqn:appendix:11}
\vartheta\!\left[
\begin{smallmatrix}
1 \\ 1 
\end{smallmatrix}
\right](y)
= -  \left[y, \mathsf{w}\right] + \frac{1}{4}[\mathsf{w}, \mathsf{w}], \qquad
\vartheta\!\left[
\begin{smallmatrix}
1 \\ 1 
\end{smallmatrix}
\right](y^\prime)
= -  \left[y^\prime, \mathsf{w}\right] + \frac{1}{4}[\mathsf{w}, \mathsf{w}]. 
\end{equation}
Since we have (\ref{align:(i)}), it follows from \eqref{eqn:appendix:01} and \eqref{eqn:appendix:11} that   
$
[y, \mathsf{v}] = [y^\prime, \mathsf{v}]$ and 
$[y, \mathsf{w}] = [y^\prime, \mathsf{w}]$. 
Since $\mathsf{v}, \mathsf{w}$ is an $\RR$-basis of $\RR^2$, 
we obtain $y = y^\prime$. 

\medskip
{\bf Case 2.}\quad 
Suppose that $y \in \mathsf{Q}_i$ and $y^\prime \in \mathsf{Q}_j$ for some $i,j = 1, \ldots , 6$ such that $j=i+1$ or such that $(i, j) = (6,1)$; namely, we assume that $y$ and $y'$ belong to distinct $\mathsf{Q}_i$ and $\mathsf{Q}_j$, respectively and that these polyhedra are adjacent.
Then without loss of generality, we assume that $(i,j) =(1,2)$, namely, $y \in \mathsf{Q}_1$ and $y^\prime \in \mathsf{Q}_2$.

By Lemma~\ref{claim:translate:in:voronoi}, we have 
$
\frac{1}{2}\mathsf{v} + y \in \frac{1}{2}\mathsf{v} + \mathsf{Q}_1 \subseteq \mathsf{V}
$. 
Applying that lemma to $\mathsf{Q}_2$ in place of $\mathsf{Q}_1$, we see that
 $\frac{1}{2}\mathsf{v} + y' \in \frac{1}{2}\mathsf{v} + \mathsf{Q}_2 \subseteq \mathsf{V}$. Then the same argument as in  Case~1 leads us to
\begin{equation}
\label{eqn:yv}
 [y, \mathsf{v}] = [y^\prime, \mathsf{v}] .
 \end{equation}

We are going to divide $\mathsf{Q}_1$ into three closed polytopes $\sigma$, $\tau$, and $\rho$, and $\mathsf{Q}_2$ into three closed polytopes $\sigma^\prime$, $\tau^\prime$, and $\rho^\prime$. Let $\mathsf{n}$ be the middle point between $\mathsf{m}_6$ and 
$\mathsf{m}_2$, i.e., $\mathsf{n} = \frac{1}{2}(\mathsf{m}_2 + \mathsf{m}_6)$. Then $\sigma$ is the closed triangle with vertices $0, \mathsf{n}, \frac{1}{2}\mathsf{w}$, and $\tau$ and $\rho$ are as in the figure below. Let $\mathsf{n}^\prime$ be the middle point between $\mathsf{m}_1$ and $\mathsf{m}_3$. Then $\sigma^\prime$ is the closed triangle with vertices $0, \mathsf{n}^\prime, \frac{1}{2}(-\mathsf{v})$, and $\tau^\prime$ and $\rho^\prime$ are as in the figure below. 

\[
\begin{tikzpicture}[scale = 0.5]
\begin{scope}
  \draw[->, >=stealth] (0,0) -- (5.8,0);
  \draw[->, >=stealth] (0,0) -- (-5.8, 0);
  \draw[->, >=stealth] (0,0) -- (3.9,5.85);
  \draw[->, >=stealth] (0,0) -- (-3.9,-5.85);
  \draw[->, >=stealth] (0,0) -- (-1.95, 5.85);
  \draw[->, >=stealth] (0,0) -- (1.95,-5.85);
  \fill (0,0) circle [radius=0.2];
  \fill (6,0) circle [radius=0.2]
  node[right]{$\mathsf{u}$}; 
  \fill (-6, 0) circle [radius=0.2]
  node[left]{$-\mathsf{u}$};
  \fill (4, 6) circle [radius=0.2]
  node[above right]{$-\mathsf{w}$};
  \fill (-4,-6) circle [radius=0.2]
  node[below left]{$\mathsf{w}$};
  \fill (-2, 6) circle [radius=0.2]
  node[above left]{$\mathsf{v}$};
  \fill (2, -6) circle [radius=0.2]
  node[below right]{$-\mathsf{v}$};
  \draw[thick] (3, -2.33)--(3, 2.33); 
  \draw[thick] (-3, -2.33)--(-3, 2.33); 
  \draw[thick] (3, 2.33)--(1, 3.66); 
  \draw[thick] (-3, 2.33)--(1, 3.66); 
  \draw[thick] (3, -2.33)--(-1, -3.66); 
  \draw[thick] (-3, -2.33)--(-1, -3.66); 
  \draw (-1, -3.66) node[below]{$\mathsf{m}_1$}; 
  \draw (3, -2.33) node[right]{$\mathsf{m}_2$}; 
  \draw (3, 2.33) node[right]{$\mathsf{m}_3$}; 
  \draw (1, 3.66) node[above]{$\mathsf{m}_4$}; 
  \draw (-3, 2.33) node[left]{$\mathsf{m}_5$}; 
  \draw (-3, -2.33) node[left]{$\mathsf{m}_6$}; 
\draw[very thick, dotted] (0, 0)--(1, -0.66)--(3, 0); 
\draw[very thick, dotted] (0, 0)--(1, -0.66)--(1, -3); 
\draw[very thick, dotted] (0, 0)--(0, -2.33)--(-2, -3); 
\draw[very thick, dotted] (0, 0)--(0, -2.33)--(1, -3); 
  \fill[gray!100,nearly transparent] 
  (-2,-3)--(-1, -3.66)--(1,-3)--(0, 0)--(-2,-3); 
   \fill[gray!100,nearly transparent] 
  (3, 0)--(3, -2.33)--(1,-3)--(0, 0)--(3, 0); 
\draw (-1.1, -1.7) node[right]{\tiny{$\sigma$}}; 
\draw (-0.5, -3) node{\tiny{$\tau$}}; 
\draw (0.3, -1.9) node{\tiny{$\rho$}}; 
\draw (0.15, -1.0) node[right]{\tiny{$\sigma^\prime$}}; 
\draw (2, -1.5) node{\tiny{$\tau^\prime$}}; 
\draw (1.2, -0.3) node{\tiny{$\rho^\prime$}}; 
 \end{scope}
\end{tikzpicture}
\]

\smallskip
{\bf Subcase 2-1.}\quad 
Suppose that $y \in \rho$ or that $y^\prime \in \sigma^\prime$. 
Noting the arrangement of the polytopes that we are focusing on is symmetric with respect to the line $\{ t \mathsf{v} \mid t \in \RR\}$, we assume that 
$y^\prime \in \sigma^\prime$ without loss of generality.

\begin{Claim} \label{claim:2-1}
We have $\sigma^\prime + \frac{1}{2}(-\mathsf{w}) \subseteq  \mathsf{V}$. 
\end{Claim}

To show the claim, since it is obvious that the translates of the vertices $0$ and $   \frac{1}{2} (-\mathsf{v})$ by $\frac{1}{2}(-\mathsf{w})$ sit in $\mathsf{V}$, we only have to show that $\mathsf{n}^\prime + \frac{1}{2}(-\mathsf{w}) \in \mathsf{V}$. 
By Lemma~\ref{lemma:perpendicular:bisector}. we have $\mathsf{w} = \mathsf{m}_1 + \mathsf{m}_6 = \mathsf{m}_1 - \mathsf{m}_3$, and hence $\mathsf{n}^\prime + \frac{1}{2}(-\mathsf{w}) = \frac{1}{2} (\mathsf{m}_1 + \mathsf{m}_3 + \mathsf{m}_3 - \mathsf{m_1}) = \mathsf{m}_3 \in \mathsf{V}$, indeed. Thus the claim holds.

By Claim~\ref{claim:2-1}, we have $\frac{1}{2} (-\mathsf{w}) +y' \in \mathsf{V}$. Then the same argument as before gives us
$
 [y, \mathsf{w}] = [y^\prime, \mathsf{w}]
$. Since $\mathsf{v}, \mathsf{w}$ is an $\RR$-basis of $\RR^2$, this equality 
together with \eqref{eqn:yv} shows that  $y = y^\prime$. 

Now, we remark that we only have to consider the case where $y \in \sigma \cup \tau$ and $y' \in \rho' \cup \tau'$ in the remaining of Case~2, because we have completed Subcase~2-1.

\smallskip
{\bf Subcase 2-2.}\quad 
Suppose that $y  \in \sigma$ and $y^\prime \in \tau^\prime$ or that $y  \in \tau$ and $y^\prime \in \rho^\prime$. Noting the symmetry of the arrangement of the polytopes, 
we only have to consider the case where $y  \in \sigma$ and $y^\prime \in \tau^\prime$; we assume this.

\begin{Claim} \label{claim:3-3}
We have $\frac{1}{2}(-\mathsf{w}) + \tau^\prime \subseteq \mathsf{u} + \mathsf{V}$. 
\end{Claim}

For this claim, it suffices to show that the translates of the vertices $\mathsf{m}_2$, $\mathsf{n}'$, $\frac{1}{2} \mathsf{u}$,  and $- \frac{1}{2} \mathsf{v}$ 
by $-\mathsf{u} - \frac{1}{2} \mathsf{w}$ sit in $\mathsf{V}$. By Lemma~\ref{lemma:perpendicular:bisector}, we have 
$\mathsf{m}_1 + \mathsf{m}_2 = -\mathsf{v}$ and 
\begin{align}
\label{align:claim:3-3}
\mathsf{m}_1 + \mathsf{m}_6 = \mathsf{w},
\end{align}
and hence $\mathsf{m}_2 - \mathsf{u} = \mathsf{m}_6$. Since $\mathsf{m}_4 = - \mathsf{m}_1$, it follows from (\ref{align:claim:3-3}) that $\mathsf{m}_4 - \mathsf{m}_6 = - \mathsf{w}$, and thus $\mathsf{m}_6 - \frac{1}{2} \mathsf{w} = \frac{1}{2} (\mathsf{m}_4 + \mathsf{m}_6)$.
It follows that 
\begin{align} \label{align:claim:3-3:2}
\mathsf{m}_2 - \mathsf{u} - \frac{1}{2} \mathsf{w} = \frac{1}{2} (\mathsf{m}_4 + \mathsf{m}_6) ,
\end{align} 
which sits in $\mathsf{V}$.
Noting (\ref{align:claim:3-3:2}), we compute
\[
\mathsf{n}' - \mathsf{u} - \frac{1}{2} \mathsf{w}  = \frac{1}{2} (\mathsf{m}_1 + \mathsf{m}_3) - \mathsf{m}_2 + \frac{1}{2} (\mathsf{m}_4 + \mathsf{m}_6)= \frac{1}{2}
((\mathsf{m}_1 + \mathsf{m}_4) + (\mathsf{m}_3 + \mathsf{m}_6)) + \mathsf{m}_5
=\mathsf{m}_5 ,
\]
which lies in $\mathsf{V}$.
We have
$\frac{1}{2} \mathsf{u} - \mathsf{u} - \frac{1}{2} \mathsf{w} =  \frac{1}{2} \mathsf{v} \in \mathsf{V}$. Finally, $-\frac{1}{2} \mathsf{v} -  \mathsf{u} - \frac{1}{2} \mathsf{w} = - \frac{1}{2} \mathsf{u} \in \mathsf{V}$.
This completes the proof of the claim.

By Claim~\ref{claim:3-3}, we have $y^\prime + \frac{1}{2}(-\mathsf{w}) - \mathsf{u}  \in \mathsf{V}$. 
Then by Remark~\ref{remark:voronoi:theta} and \eqref{eqn:xi:translation:vw}, we have  
\[
\vartheta\!\left[
\begin{smallmatrix}
1 \\ 1 \end{smallmatrix}
\right](y^\prime)
= 
-  [y^\prime, 2\mathsf{u}+ \mathsf{w}]
+ \frac{1}{4} [2\mathsf{u}+ \mathsf{w}, 2\mathsf{u}+\mathsf{w}]. 
\]
By Lemma~\ref{claim:translate:in:voronoi}, we have
$y+\frac{1}{2} (-\mathsf{w})  \in \mathsf{Q}_1 + \frac{1}{2} (-\mathsf{w}) \subseteq \mathsf{V}$, and hence by Remark~\ref{remark:voronoi:theta}, we have $\vartheta\!\left[
\begin{smallmatrix}
0 \\ 0 \end{smallmatrix}
\right](y) = 0$. It follows from (\ref{eqn:xi:translation}) that
\begin{align*}
\vartheta\!\left[
\begin{smallmatrix}
1 \\ 1 \end{smallmatrix}
\right](y)
= -  \left[
y, \mathsf{w}\right] + \frac{1}{4} \left[\mathsf{w}, \mathsf{w}
\right]. 
\end{align*}
Thus (\ref{align:(i)}) for $j_1=j_2=1$ implies that 
\begin{equation}
\label{eqn:33:yw}
[y, \mathsf{w}] = 
[y^\prime, \mathsf{u} - \mathsf{v}] + [\mathsf{u}, \mathsf{v}]. 
\end{equation}

The same argument as in the proof of Claim~\ref{claim:2-1} shows that $\sigma + \frac{1}{2} \mathsf{u} \subseteq \mathsf{V}$, and hence
$y + \frac{1}{2}\mathsf{u} \in \mathsf{V}$. 
Applying Lemma~\ref{claim:translate:in:voronoi} to $\mathsf{Q}_2$ in place of $\mathsf{Q}_1$, we have
$y^\prime + \frac{1}{2}(-\mathsf{u}) \in \frac{1}{2} (- \mathsf{u}) + \mathsf{Q}_2 \subseteq 
\mathsf{V}$. 
Then by Remark~\ref{remark:voronoi:theta} and (\ref{eqn:xi:translation}), we have
\begin{align*}
\vartheta\!\left[
\begin{smallmatrix}
1 \\ 0 \end{smallmatrix}
\right](y)
=  \left[y^\prime, \mathsf{u}\right] + \frac{1}{4}[\mathsf{u}, \mathsf{u}], 
\qquad
\vartheta\!\left[
\begin{smallmatrix}
1 \\ 0 \end{smallmatrix}
\right](y^\prime)
= 
-  \left[y^\prime, \mathsf{u}\right] + \frac{1}{4}[\mathsf{u}, \mathsf{u}]. 
\end{align*}
Thus (\ref{align:(i)}) for $(j_1,j_2)=(1,0)$ implies that 
\begin{equation}
\label{eqn:33:yu}
 [y, \mathsf{u}] 
= [-y^\prime, \mathsf{u}]
\end{equation}
By \eqref{eqn:yv} and \eqref{eqn:33:yu}, we have
$[y, \mathsf{w}] = 
[y, -\mathsf{u} - \mathsf{v}]
=  [y^\prime, \mathsf{u}] - [y^\prime, \mathsf{v}]
= [y^\prime, \mathsf{u} - \mathsf{v}] $,  and by \eqref{eqn:33:yw}, it follows that $[\mathsf{u}, \mathsf{v}] =0$. However,
this contradicts that the Voronoi cell
$\mathsf{V}$ is a hexagon (cf. Remark~\ref{remark:angle:voronoi}). This  concludes that Subcase~2-2 does not occur. 

\smallskip
{\bf Subcase 2-3.}\quad 
Suppose that $y \in \sigma$ and $y^\prime \in \rho^\prime$. 
We see that $y + \frac{1}{2} \mathsf{u} \in \mathsf{V}$ and 
$y^\prime + \frac{1}{2}(-\mathsf{u}) \in \mathsf{V}$. 
Then  by Remark 5.4 and (5.4), we have 
\[
\vartheta\!\left[
\begin{smallmatrix}
1 \\ 0 \end{smallmatrix}
\right](y)
= 
[y, \mathsf{u}]
+ \frac{1}{4} [\mathsf{u}, \mathsf{u}].  
\qquad 
\vartheta\!\left[
\begin{smallmatrix}
1 \\ 0 \end{smallmatrix}
\right](y^\prime)
= 
-  [y^\prime, \mathsf{u}]
+ \frac{1}{4} [\mathsf{u}, \mathsf{u}].  
\]
Thus $[y, \mathsf{u}] = [-y^\prime, \mathsf{u}]$. 
We also see that 
$y + \frac{1}{2} (-\mathsf{w}) \in \mathsf{V}$ and 
$y^\prime + \frac{1}{2} \mathsf{w} \in \mathsf{V}$, and we obtain $[y, \mathsf{w}] = [-y^\prime, \mathsf{w}]$. Thus we get $y = -y^\prime$. We have $y = y^\prime = 0$ in Subcase 2-3.

\smallskip
{\bf Subcase 2-4.}\quad
Suppose that $y  \in \tau$ and $y^\prime \in \tau^\prime$. The same argument as the proof of Claim~\ref{claim:3-3} shows that $\frac{1}{2}\mathsf{u} + \tau \subseteq -\mathsf{v} +\mathsf{V}$, and hence $y + \frac{1}{2}\mathsf{u} \in -\mathsf{v} +\mathsf{V}$. Further, the argument applying Lemma~\ref{claim:translate:in:voronoi} to $\mathsf{Q}_2$ in place of $\mathsf{Q}_1$ shows that 
$y^\prime + \frac{1}{2}(-\mathsf{u}) \in \mathsf{V}$. 
Noting  (\ref{align:(i)}) for $(j_1 , j_2) = (1,0)$, Remark~\ref{remark:voronoi:theta}, \eqref{eqn:xi:translation:vw}, and \eqref{eqn:xi:translation}, we 
obtain 
$
[y, \mathsf{u}+ 2\mathsf{v}] + \frac{1}{4} [\mathsf{u}+ 2\mathsf{v}, \mathsf{u}+ 2\mathsf{v}]
= - [y^\prime, \mathsf{u}] + \frac{1}{4} [\mathsf{u}, \mathsf{u}], 
$
and thus
\begin{equation}
\label{eqn:35:yu}
[y, \mathsf{u}+ 2\mathsf{v}] + [\mathsf{v}, \mathsf{u}+\mathsf{v}]
= - [y^\prime, \mathsf{u}]. 
\end{equation}
Similarly, we see that $y + \frac{1}{2}(-\mathsf{w}) \in \mathsf{V}$ and $y^\prime + \frac{1}{2}(-\mathsf{w}) \in \mathsf{u} + \mathsf{V}$. Noting  (\ref{align:(i)}) for $(j_1 , j_2) = (1,1)$, 
Remark~\ref{remark:voronoi:theta}, and \eqref{eqn:xi:translation:vw}, 
we obtain 
$
[y, -\mathsf{w}] + \frac{1}{4} [\mathsf{w}, \mathsf{w}]
= - [y^\prime, \mathsf{u} + 2 \mathsf{w}] + \frac{1}{4} [\mathsf{u} + 2 \mathsf{w}, \mathsf{u} + 2 \mathsf{w}], 
$
and thus
\begin{equation}
\label{eqn:35:yw}
- [y, \mathsf{u}+ \mathsf{v}] 
= [y^\prime, \mathsf{u}-\mathsf{v}] + [\mathsf{u}, \mathsf{v}]. 
\end{equation}
Summing up \eqref{eqn:yv}, \eqref{eqn:35:yu}, and \eqref{eqn:35:yw}, we get $[y, \mathsf{v}] = [y^\prime, \mathsf{v}] = -\frac{1}{2} [\mathsf{v}, \mathsf{v}]$. Then by \eqref{eqn:35:yu}, we have 
$\left[y + \frac{1}{2} \mathsf{v}, \mathsf{u}\right] = - \left[y^\prime  + \frac{1}{2}\mathsf{v}, \mathsf{u}\right]$. Since 
$\left[y + \frac{1}{2} \mathsf{v}, \mathsf{v}\right] = 0 = - \left[y^\prime  + \frac{1}{2}\mathsf{v}, \mathsf{v}\right]$, 
we conclude that $y = y^\prime = \frac{1}{2}(-\mathsf{v})$ 
in this subcase. 

\medskip
{\bf Case 3.}\quad 
Suppose that $y \in \mathsf{Q}_i$ and $y^\prime \in \mathsf{Q}_j$, where $\mathsf{Q}_i$ and $\mathsf{Q}_j$ are distinct and not adjacent. Without loss of generality, we assume that $(i,j) = (1,3)$.  
In this case, we have $- y^\prime \in \mathsf{Q}_6$, and $\mathsf{Q}_6$ is adjacent to $\mathsf{Q}_1$. Note that the tropical theta functions that we consider is even, and the assumption (i) holds for $-y'$ in place of $y'$. Then by Case~2, we have $y=-y'$, which is (ii) in this case.

By the argument from Cases 1 to 3, we see that if $y, y^\prime \in \mathsf{Q}$, then 
$y = y^\prime$ or $y = - y^\prime$. 
This completes the proof of Proposition~\ref{prop:xi}. 
\QED

\subsection{Proof of the unimodularity in Proposition~\ref{prop:thm:faith:embeddings:1:a}}
\label{subsec:proof:unimodularity}

We keep the notation in the previous subsection. Recall that the integral structure on $\RR^2$ and $X := \RR^2/\mathsf{P}$ is given by a fixed lattice $\mathsf{N} \subseteq \RR^2$. Assume that the Voronoi cell $\mathsf{V}$ of $\mathsf{P}$ is a hexagon. As is noted in Remark~\ref{remark:unimodular:orbifold}, it suffices to show that the map 
\[
\varphi\colon \RR^2 \to \RR^3, \; x \mapsto \left(
\vartheta\!\left[
\begin{smallmatrix} 1 \\ 0 \end{smallmatrix}\right](x) 
- \vartheta\!\left[
\begin{smallmatrix} 0 \\ 0 \end{smallmatrix}\right](x), \, 
\vartheta\!\left[
\begin{smallmatrix} 0 \\ 1 \end{smallmatrix}\right](x) 
- \vartheta\!\left[
\begin{smallmatrix} 0 \\ 0 \end{smallmatrix}\right](x), \, 
\vartheta\!\left[
\begin{smallmatrix} 1 \\ 1 \end{smallmatrix}\right](x)
- \vartheta\!\left[
\begin{smallmatrix} 0 \\ 0 \end{smallmatrix}\right](x) 
\right)
\] 
is unimodular.
Note that $\bigcup_{\mathsf{p} \in \mathsf{P}} (\mathsf{p} + \relin (\mathsf{V}))$ is a dense open subset of $\RR^2$. By \cite[Lemma~6.2]{KY}, it suffices to that $\varphi$ is unimodular on $\bigcup_{\mathsf{p} \in \mathsf{P}} (\mathsf{p} + \relin (\mathsf{V}))$. By the quasi-periodicity of the tropical theta functions, it suffices to show that $\rest{\varphi}{\relin (\mathsf{V})}\colon \relin (\mathsf{V}) \to \RR^3$ is unimodular. 
 Further, since $\bigcup_{i=1}^6 \relin (\mathsf{Q}_i)$ is a dense open subset of $\relin (\mathsf{V})$, we only have to show the following proposition.

\begin{Proposition}
\label{prop:thm:faith:embeddings:1:b}
Assume that $\mathsf{V}$ is a hexagon. 
For $1 \leq i \leq 6$, let $\mathsf{Q}_i \subseteq \mathsf{V}$ be 
as in the proof of Proposition~\ref{prop:thm:faith:embeddings:1:a}. Then 
\[
\varphi\colon \relin (\mathsf{Q}_i) \to \RR^3, \; x \mapsto \left(
\vartheta\!\left[
\begin{smallmatrix} 1 \\ 0 \end{smallmatrix}\right](x) 
- \vartheta\!\left[
\begin{smallmatrix} 0 \\ 0 \end{smallmatrix}\right](x), \, 
\vartheta\!\left[
\begin{smallmatrix} 0 \\ 1 \end{smallmatrix}\right](x) 
- \vartheta\!\left[
\begin{smallmatrix} 0 \\ 0 \end{smallmatrix}\right](x), \, 
\vartheta\!\left[
\begin{smallmatrix} 1 \\ 1 \end{smallmatrix}\right](x)
- \vartheta\!\left[
\begin{smallmatrix} 0 \\ 0 \end{smallmatrix}\right](x) 
\right)
\] 
is unimodular. 
\end{Proposition}

\Proof
We keep the notation in the proof of Proposition~\ref{prop:thm:faith:embeddings:1:a}. 
We will prove the assertion only for $\mathsf{Q}_1$ because
the same arguments work for other $\mathsf{Q}_i$'s. 

Let $x \in \mathsf{Q}_1$. By Remark~\ref{remark:voronoi:theta},
\eqref{eqn:appendix:01}, and 
\eqref{eqn:appendix:11}, 
we have 
\[
\vartheta\!\left[
\begin{smallmatrix}
0 \\ 0 
\end{smallmatrix}
\right](x)
= 0,  
\quad
\vartheta\!\left[
\begin{smallmatrix}
0 \\ 1 
\end{smallmatrix}
\right](x)
=  \left[x, \mathsf{v}\right] + \frac{1}{4}[\mathsf{v}, \mathsf{v}], 
\quad
\vartheta\!\left[
\begin{smallmatrix}
1 \\ 1 
\end{smallmatrix}
\right](x)
= -  \left[x, \mathsf{w}\right] + \frac{1}{4}[\mathsf{w}, \mathsf{w}]. 
\]
Thus it suffices to show that the map $\RR^2 \to \RR^3$ given by $x \mapsto \left(\vartheta\!\left[
\begin{smallmatrix} 1 \\ 0 \end{smallmatrix}\right](x)  , [x , \mathsf{v}] , [x, -\mathsf{w}] \right)$ is unimodular. To see that, we only have to show that the $\ZZ$-linear map $h\colon \mathsf{N} \to \ZZ^2$ given by $x \mapsto  ([x , \mathsf{v}] , [x, -\mathsf{w}])$ is an isomorphism. Let $\lambda\colon \mathsf{P} \to \mathsf{N}^*$ be the $\ZZ$-linear map corresponding to $[\ndot , \ndot]$. Then $h = (\lambda (\mathsf{v}) ,\lambda (- \mathsf{w}))$. Since $[\ndot,\ndot]$ is a principal polarization, $\lambda$ is an isomorphism, and since $\mathsf{v}, - \mathsf{w}$ is a basis of the lattice $\mathsf{P}$, $\lambda (\mathsf{v}) ,\lambda (- \mathsf{w})$ is a $\ZZ$-basis of $\mathsf{N}^*$. This shows that the map $h$ is an isomorphism, which completes the proof.
\QED

\section{Proof of Theorem~\ref{thm:faith:embeddings2}} \label{section:proof:(2)}
\label{sec:tropical:Kummer:quartic}
In this section, we prove Theorem~\ref{thm:faith:embeddings2}. 
We may assume that $(X, Q) = (\RR^2/\mathsf{P}, [\ndot, \ndot])$ as in Section~\ref{subsec:reduction}. 
Let $Y \colonequals X/\langle -1\rangle$ be the tropical Kummer surface. 
Then, by Propositions~\ref{prop:thm:faith:embeddings:1:a} and \ref{prop:thm:faith:embeddings:1:b}, 
the map 
\[
\psi\colon Y \to \RR^3, 
\quad 
[x] \mapsto 
\left(
\vartheta\!\left[
\begin{smallmatrix} 1 \\ 0 \end{smallmatrix}\right](x) 
- \vartheta\!\left[
\begin{smallmatrix} 0 \\ 0 \end{smallmatrix}\right](x), \, 
\vartheta\!\left[
\begin{smallmatrix} 0 \\ 1 \end{smallmatrix}\right](x) 
- \vartheta\!\left[
\begin{smallmatrix} 0 \\ 0 \end{smallmatrix}\right](x), \, 
\vartheta\!\left[
\begin{smallmatrix} 1 \\ 1 \end{smallmatrix}\right](x)
- \vartheta\!\left[
\begin{smallmatrix} 0 \\ 0 \end{smallmatrix}\right](x) 
\right)
\]
faithfully embeds $Y$ into $\RR^3$, where $\vartheta\!\left[\begin{smallmatrix} j_1 \\ j_2 \end{smallmatrix}\right](x)$ is defined by \eqref{eqn:appendix:xi} for $j_1 , j_2 \in \{ 0,1 \}$.
To ease notation, we set, for $j_1, j_2 \in \{0, 1\}$ with $(j_1, j_2) \neq (0, 0)$,
\[
\theta_{j_1 j_2} \colonequals \vartheta\!\left[
\begin{smallmatrix} j_1 \\  j_2 \end{smallmatrix}\right](0)
. 
\]
We note that $\vartheta\!\left[\begin{smallmatrix} 0 \\  0 \end{smallmatrix}\right](0) = 0$. Then the eight points in $\RR^3$ in 
Theorem~\ref{thm:faith:embeddings2} are given by 
\begin{align*}
& \tau_{0}  \colonequals ( \theta_{10} , \theta_{01} , \theta_{11}), && 
\tau_{1}  \colonequals (-\theta_{10} , \theta_{11} -\theta_{10}, \theta_{01}-\theta_{10}), 
\\
& \tau_{2}  \colonequals (\theta_{11}-\theta_{01} , -\theta_{01}, \theta_{10}-\theta_{01}), && 
\tau_{3}  \colonequals ( \theta_{01} -\theta_{11}, \theta_{10}-\theta_{11} , \theta_{11}), 
\\
& \widetilde{\tau}_{0}  \colonequals ( -\theta_{10} , -\theta_{01} , -\theta_{11}), 
&& 
\widetilde{\tau}_{1}  \colonequals ( \theta_{10} ,  -\theta_{11} + \theta_{10}, -\theta_{01} + \theta_{10}), 
\\
& \widetilde{\tau}_{2}   \colonequals ( -\theta_{11}+\theta_{01} , \theta_{01}, -\theta_{10}+\theta_{01}), && 
\widetilde{\tau}_{3}  \colonequals ( -\theta_{01} +\theta_{11} , -\theta_{10}+\theta_{11} , \theta_{11}). 
\end{align*}
By Remark~\ref{remark:voronoi:theta}, we also remark that  
\[
\theta_{00} = 0, 
\quad 
\theta_{10} = \frac{1}{4} [\mathsf{u}, \mathsf{u}], 
\quad 
\theta_{01} = \frac{1}{4} [\mathsf{v}, \mathsf{v}],  
\quad 
\theta_{11} = \frac{1}{4} [\mathsf{w}, \mathsf{w}]. 
\]

Then Theorem~\ref{thm:faith:embeddings2} amounts to the following. 

\begin{Proposition} \label{prop:TKQS}
Assume that $\mathsf{V}$ is a hexagon. Then $\psi(Y) \subseteq \RR^3$ is a parallelepiped whose eight vertices are given by 
$
\tau_{0}, \tau_{1}, \tau_{2}, \tau_{3}, 
\widetilde{\tau}_{0}, \widetilde{\tau}_{1}, 
\widetilde{\tau}_{2}, \widetilde{\tau}_{3}.
$ 
\end{Proposition}

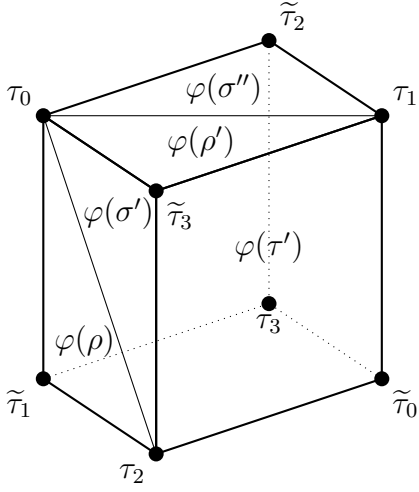
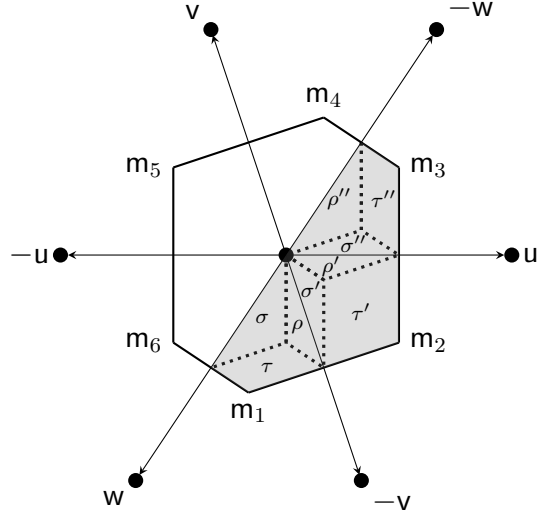
\[
\begin{tikzpicture}[scale = 0.5]
\begin{scope}
\draw[thick] (0, 0)--(0, 7)--(6, 9)--(6, 2)--(0, 0); 
\draw[thick] (0, 0)--(0, 7)--(-3, 9)--(-3, 2)--(0, 0); 
\draw[thick] (0, 7)--(6, 9)--(3, 11)--(-3, 9)--(0, 7); 
\draw (-3, 9)--(6, 9); 
\draw (-3, 9)--(0, 0); 
\draw[dotted] (3, 11)--(3, 4); 
\draw[dotted] (-3, 2)--(3, 4); 
\draw[dotted] (6, 2)--(3, 4); 
\draw (-1.9, 3) node{$\varphi(\rho)$}; 
\draw (-1, 6.4) node{$\varphi(\sigma^\prime)$}; 
\draw (3, 5.5) node{$\varphi(\tau^\prime)$}; 
\draw (1.2, 8.3) node{$\varphi(\rho^\prime)$}; 
\draw (1.8, 9.7) node{$\varphi(\sigma^{\prime\prime})$}; 
\fill (0, 0) circle [radius=0.2] node[below left]{$\tau_{2}$};
\fill (0, 7) circle [radius=0.2] node[below right]{$\widetilde{\tau}_{3}$};
\fill (6, 9) circle [radius=0.2] node[above right]{$\tau_{1}$};
\fill (6, 2) circle [radius=0.2] node[below right]{$\widetilde{\tau}_{0}$};
\fill (3, 11) circle [radius=0.2] node[above right]{$\widetilde{\tau}_{2}$};
\fill (-3, 9) circle [radius=0.2] node[above left]{$\tau_{0}$};
\fill (3,4) circle [radius=0.2] node[below]{$\tau_{3}$};
\fill (-3, 2) circle [radius=0.2] node[below left]{$\widetilde{\tau}_{1}$};
\end{scope}
  \begin{scope}[shift={($(18, 5)$)}]
  \draw[->, >=stealth] (0,0) -- (5.8,0);
  \draw[->, >=stealth] (0,0) -- (-5.8, 0);
  \draw[->, >=stealth] (0,0) -- (3.9,5.85);
  \draw[->, >=stealth] (0,0) -- (-3.9,-5.85);
  \draw[->, >=stealth] (0,0) -- (-1.95, 5.85);
  \draw[->, >=stealth] (0,0) -- (1.95,-5.85);
  \fill (0,0) circle [radius=0.2];
  \fill (6,0) circle [radius=0.2]
  node[right]{$\mathsf{u}$}; 
  \fill (-6, 0) circle [radius=0.2]
  node[left]{$-\mathsf{u}$};
  \fill (4, 6) circle [radius=0.2]
  node[above right]{$-\mathsf{w}$};
  \fill (-4,-6) circle [radius=0.2]
  node[below left]{$\mathsf{w}$};
  \fill (-2, 6) circle [radius=0.2]
  node[above left]{$\mathsf{v}$};
  \fill (2, -6) circle [radius=0.2]
  node[below right]{$-\mathsf{v}$};
  \draw[thick] (3, -2.33)--(3, 2.33); 
  \draw[thick] (-3, -2.33)--(-3, 2.33); 
  \draw[thick] (3, 2.33)--(1, 3.66); 
  \draw[thick] (-3, 2.33)--(1, 3.66); 
  \draw[thick] (3, -2.33)--(-1, -3.66); 
  \draw[thick] (-3, -2.33)--(-1, -3.66); 
  \draw (-1, -3.66) node[below]{$\mathsf{m}_1$}; 
  \draw (3, -2.33) node[right]{$\mathsf{m}_2$}; 
  \draw (3, 2.33) node[right]{$\mathsf{m}_3$}; 
  \draw (1, 3.66) node[above]{$\mathsf{m}_4$}; 
  \draw (-3, 2.33) node[left]{$\mathsf{m}_5$}; 
  \draw (-3, -2.33) node[left]{$\mathsf{m}_6$}; 
\draw[very thick, dotted] (0, 0)--(1, -0.66)--(3, 0); 
\draw[very thick, dotted] (0, 0)--(1, -0.66)--(1, -3); 
\draw[very thick, dotted] (0, 0)--(0, -2.33)--(-2, -3); 
\draw[very thick, dotted] (0, 0)--(0, -2.33)--(1, -3); 
\draw[very thick, dotted] (2, 3)--(2, 0.66)--(0, 0); 
\draw[very thick, dotted] (2, 3)--(2, 0.66)--(3, 0); 
      \fill[gray!100,nearly transparent] 
  (-2,-3)--(-1, -3.66)--(1,-3)--(0, 0)--(-2,-3); 
   \fill[gray!100,nearly transparent] 
  (3, 0)--(3, -2.33)--(1,-3)--(0, 0)--(3, 0); 

   \fill[gray!100,nearly transparent] 
  (3, 0)--(3, 2.33)--(2,3)--(0, 0)--(3, 0); 
\draw (-1.1, -1.7) node[right]{\tiny{$\sigma$}}; 
\draw (-0.5, -3) node{\tiny{$\tau$}}; 
\draw (0.3, -1.9) node{\tiny{$\rho$}}; 
\draw (0.1, -0.9) node[right]{\tiny{$\sigma^\prime$}}; 
\draw (2, -1.5) node{\tiny{$\tau^\prime$}}; 
\draw (1.2, -0.3) node{\tiny{$\rho^\prime$}}; 
\draw (0.8, 1.5) node[right]{\tiny{$\rho^{\prime\prime}$}}; 
\draw (1.8, 0.3) node{\tiny{$\sigma^{\prime\prime}$}}; 
\draw (2.6, 1.5) node{\tiny{$\tau^{\prime\prime}$}}; 
 \end{scope}
\end{tikzpicture}
\]

\Proof
We keep the notation in the proof of Proposition~\ref{prop:thm:faith:embeddings:1:a}. 
We may and do assume that $x \in \mathsf{V}$. Then by Remark~\ref{remark:voronoi:theta}, we note that $\vartheta\!\left[\begin{smallmatrix} 0 \\ 0 \end{smallmatrix}\right](x) = 0$, and hence on $\mathsf{V}$, the map $\varphi\colon \RR^2 \to \RR^3$ is given by $x \mapsto \left(
\vartheta\!\left[
\begin{smallmatrix} 1 \\ 0 \end{smallmatrix}\right](x) , 
\vartheta\!\left[
\begin{smallmatrix} 0 \\ 1 \end{smallmatrix}\right](x) , 
\vartheta\!\left[
\begin{smallmatrix} 1 \\ 1 \end{smallmatrix}\right](x)
\right)$. 
Recall that $\mathsf{Q}_i$ denotes the quadrilateral containing $\mathsf{m}_i$. 
Let $(T_1 , T_2 ,T_3)$ be the standard coordinates on $\RR^3$.

{\bf Case 1.}\quad 
Suppose that $x \in \mathsf{Q}_1$. Then by \eqref{eqn:appendix:01} and \eqref{eqn:appendix:11}, we have 
\begin{align*}
& \vartheta\!\left[
\begin{smallmatrix} 0 \\ 1 \end{smallmatrix}\right](x)
=  \left[
x, \mathsf{v}
\right] + \frac{1}{4} \left[
\mathsf{v}, \mathsf{v}
\right], 
\quad 
\vartheta\!\left[
\begin{smallmatrix} 1 \\ 1 \end{smallmatrix}\right](x)
=  \left[
x, -\mathsf{w}
\right] + \frac{1}{4} \left[
\mathsf{w}, \mathsf{w}
\right]. 
\end{align*}

\smallskip
{\bf Subcase 1-1.}\quad 
Suppose that $x \in \sigma$. Since $x + \frac{1}{2}\mathsf{u} \in \mathsf{V}$, we have 
\[
\vartheta\!\left[\begin{smallmatrix} 1 \\ 0 \end{smallmatrix}\right](x)
=  \left[
x, \mathsf{u}
\right] + \frac{1}{4} \left[
\mathsf{u}, \mathsf{u} 
\right].
\]
Then we have 
\[
\vartheta\!\left[\begin{smallmatrix} 1 \\ 0 \end{smallmatrix}\right](x)
+\vartheta\!\left[\begin{smallmatrix} 0 \\ 1 \end{smallmatrix}\right](x)
- \vartheta\!\left[\begin{smallmatrix} 1 \\ 1 \end{smallmatrix}\right](x)
= \frac{1}{4} \left[\mathsf{u}, \mathsf{u} \right]
+  \frac{1}{4} \left[\mathsf{v}, \mathsf{v} \right]
-  \frac{1}{4} \left[\mathsf{w}, \mathsf{w} \right]
= 
\theta_{10} + \theta_{01} - \theta_{11}
\]
It follows that $\varphi(\sigma)$ is the triangle on the plane 
\[
 T_1 + T_2 - T_3 = \theta_{10} + \theta_{01} - \theta_{11}
\]
with vertices $\varphi(0), \varphi(\frac{1}{2}\mathsf{w})$, and $\varphi(\mathsf{n})$, where we recall that $\mathsf{n} = \frac{1}{2} (\mathsf{m}_2 + \mathsf{m}_6 )$. 

\smallskip
{\bf Subcase 1-2.}\quad 
Suppose that $x \in \tau$. 
Since $x + \frac{1}{2}\mathsf{u} \in -\mathsf{v} +\mathsf{V}$, 
we have 
\[
\theta\!\left[\begin{smallmatrix} 1 \\ 0 \end{smallmatrix}\right](x)
= 
[y, \mathsf{u}+ 2\mathsf{v}]
+ \frac{1}{4} [\mathsf{u}+ 2\mathsf{v}, \mathsf{u}+ 2\mathsf{v}]. \]
Then we have 
\[
\theta\!\left[\begin{smallmatrix} 1 \\ 0 \end{smallmatrix}\right](x)
-\theta\!\left[\begin{smallmatrix} 0 \\ 1 \end{smallmatrix}\right](x)
-\theta\!\left[\begin{smallmatrix} 1 \\ 1 \end{smallmatrix}\right](x)
= \frac{1}{4} [\mathsf{u}+ 2\mathsf{v}, \mathsf{u}+ 2\mathsf{v}]
-  \frac{1}{4} \left[\mathsf{v}, \mathsf{v} \right]
-  \frac{1}{4} \left[\mathsf{w}, \mathsf{w} \right]
= - \theta_{10} + \theta_{01} + \theta_{11}
.
\]
It follows that $\varphi(\sigma)$ is the quadrilateral on the plane 
\[
 T_1 - T_2 - T_3 = 
- \theta_{10} + \theta_{01} + \theta_{11}
\]
with vertices $\varphi(\frac{1}{2}\mathsf{w}), \varphi(\frac{1}{2}(-\mathsf{v})), \varphi(\mathsf{n})$, and $\varphi(\mathsf{m}_1)$.

We can compute similarly, and we get Table~\ref{cap:table} below. 
Note that Case 1-3 corresponds to $x \in \rho$ and that 
Cases 2-1, 2-2, 2-3 correspond to $x \in \sigma^\prime, \tau^\prime, \rho^\prime$, respectively, 
and  Cases 3-1, 3-2, 3-3 correspond to $x \in \sigma^{\prime\prime}, \tau^{\prime\prime}, \rho^{\prime\prime}$, respectively. 
\begin{table}[htb!]
\begin{tabular}{|l|l|l|l|}
\hline
Cases       & Faces             & Equations                       & Descriptions                                 \\ \hline

\begin{tabular}[c]{@{}l@{}}
1-1 
\\
3-3
\end{tabular}
& 
\begin{tabular}[c]{@{}l@{}}
$\varphi(\sigma)$ 
\\
$\varphi(\rho^{\prime\prime})$ 
\end{tabular}
& 
\begin{tabular}[c]{@{}l@{}}
\vspace{-3mm}\\
$T_1 + T_2 - T_3$ 
\\ 
$=  \theta_{10} +  \theta_{01} - \theta_{11} 
$
\vspace{2mm}
\end{tabular}
 & 
\begin{tabular}[c]{@{}l@{}} 
\vspace{-3mm}\\
quadrilateral with vertices $\tau_{0}, \tau_{3}, \widetilde{\tau}_1, \widetilde{\tau}_2$
 \\
 (parallel to $\varphi(\tau^\prime)$)
 \vspace{2mm}
 \end{tabular}
 \\ \hline
1-2    &  $\varphi(\tau)$  &   
\begin{tabular}[c]{@{}l@{}}
\vspace{-3mm}\\
$T_1 - T_2  - T_3$ 
\\ 
$= - \theta_{10} + \theta_{01} + \theta_{11}$ 
\vspace{2mm}
\end{tabular}
& 
\begin{tabular}[c]{@{}l@{}} 
 quadrilateral with vertices $\tau_{2}, \tau_{3}, \widetilde{\tau}_0,  \widetilde{\tau}_1$   
\\
(parallel to $\varphi(\rho^\prime)\cup\varphi(\sigma^{\prime\prime})$)
\end{tabular}
 \\ 
\hline
\begin{tabular}[c]{@{}l@{}}
1-3         \\
2-1
\end{tabular}
&  
\begin{tabular}[c]{@{}l@{}}
$\varphi(\rho)$ \\
$\varphi(\sigma^\prime)$
\end{tabular}
  &              
\begin{tabular}[c]{@{}l@{}}
\vspace{-3mm}\\
$T_1 - T_2  + T_3$ 
\\ 
$=  \theta_{10} - \theta_{01} +  \theta_{11}$ 
\vspace{2mm}
\end{tabular}                  &       
\begin{tabular}[c]{@{}l@{}}                        
quadrilateral  with vertices $\tau_{0}, \tau_{2}, \widetilde{\tau}_1, \widetilde{\tau}_3$   
\\
(parallel to $\varphi(\tau^{\prime\prime})$)     
\end{tabular}           
  \\ \hline
2-2        &   $\varphi(\tau^\prime)$         
       &     \begin{tabular}[c]{@{}l@{}}
\vspace{-3mm}\\
$T_1 + T_2  - T_3$ 
\\ 
$= - \theta_{10}
- \theta_{01} + \theta_{11}$
\vspace{2mm}
\end{tabular}                             &     
\begin{tabular}[c]{@{}l@{}}             
quadrilateral with vertices $\tau_{1}, \tau_{2}, \widetilde{\tau}_0, \widetilde{\tau}_3$         \\
(parallel to $\varphi(\sigma) \cup \varphi(\rho^{\prime\prime})$)                          
\end{tabular} 
    \\ \hline
 \begin{tabular}[c]{@{}l@{}}   
2-3       
\\
3-1
\end{tabular} 
 &  
 \begin{tabular}[c]{@{}l@{}}   
 $\varphi(\rho^\prime)$ 
 \\
 $\varphi(\sigma^{\prime\prime})$ 
\end{tabular}       
          &    \begin{tabular}[c]{@{}l@{}}
\vspace{-3mm}\\
$T_1 - T_2  - T_3$ 
\\ 
$= 
\theta_{10}
- \theta_{01} - \theta_{11}$
\vspace{2mm}
\end{tabular}                              &                  
\begin{tabular}[c]{@{}l@{}}   
quadrilateral with vertices $\tau_{0}, \tau_{1}, \widetilde{\tau}_2, \widetilde{\tau}_3$   
\\
(parallel to $\varphi(\tau)$)
\end{tabular}  
 \\ \hline
3-2         &  $\varphi(\tau^{\prime\prime})$         
        &   \begin{tabular}[c]{@{}l@{}}
\vspace{-3mm}\\
$T_1 - T_2  + T_3$ 
\\ 
$= - \theta_{10}+ \theta_{01}
- \theta_{11} $
\vspace{2mm}
\end{tabular}                               &   
\begin{tabular}[c]{@{}l@{}}
quadrilateral with vertices $\tau_{1}, \tau_{3}, \widetilde{\tau}_0, \widetilde{\tau}_2$        
\\
(parallel to $\varphi(\rho)\cup\varphi(\sigma^\prime)$)       
  \end{tabular}           
\\ \hline
\end{tabular}
\smallskip
\caption{Defining equations of $\psi(Y)$}
\label{cap:table}
\end{table}
This completes the proof.
\QED

\section{Uniformization theory for totally degenerate abelian varieties}
\label{sec:uniformization}
In the following sections, we study the relation between 
Kummer surfaces defined over nonarchimedean valued fields and 
tropical Kummer quartic surfaces. 
In this section, we recall the uniformization theory for totally degenerate abelian varieties with a focus on the descent theory of line bundles and nonarchimedean theta functions.

Let $k$ be an algebraically closed field that is complete with respect to a nontrivial nonarchimedean absolute value $|\ndot|_k$. 
We sometimes call such a valued field an \emph{algebraically closed, nontrivial valued, and nonarchimdean field}.

\subsection{Analytification}
\label{subsec:analytification}
We mean by an algebraic scheme over $k$ a scheme over $k$ that is separated and of finite type. 
For an algebraic scheme $X$ over $k$, let $X^{\an}$ denote 
the analytic space associated to $X$ in the sense of Berkovich. 
We refer to \cite{Berkovich} for a general reference on Berkovich spaces. 
As a set, 
\[
X^{\an} =\left\{ (p ,|\ndot | ) \mid \text{$p \in X$, $|\ndot|$ is an absolute value on $\kappa (p)$ that extends $|\ndot|_k$}\right\}, 
\]
where $\kappa (p)$ denotes the residue field of $X$ at $p$. For $x = (p,|\ndot|) \in X^{\an}$ and $f \in \OO_{X, p}$, we write $|f(x)| \colonequals |f(p)|$, where $f(p)$ is the class of $f$ in $\kappa (p)$. The topology on $X^{\an}$ is 
the weakest topology such that for any 
Zariski open subset $U \subseteq X$ and for any $f \in \OO_{X}(U)$, $U^{\an}$ is an open subset of $X$ and the map $|f|\colon U^{\an} \ni x \mapsto |f(x)| \in \RR$ is continuous.

If $x \in X(k)$, then $\kappa(x) = k$ and 
$(x, |\ndot|_k) \in X^{\an}$.  Via $x \mapsto (x, |\ndot|_k)$, 
we regard $X(k)$ as a subset of $X^{\an}$. A point in $X(k)$ is called a \emph{classical point}.

Let $f\colon X \to Y$ be a morphism of algebraic schemes over $k$. Then $f$ induces a morphism $f^{\an}\colon X^{\an} \to Y^{\an}$ by $f^{\an} (p,|\ndot|) = (f(p) , |f^* (\ndot)|)$. We note that $f^{\an} (X(k)) \subseteq Y(k)$.

\subsection{Nonarchimedean Appell--Humbert theory for totally degenerate abelian varieties}
\label{subsec:na:ah}
We recall the uniformization theory by Raynaud \cite{Raynaud} and Bosch--L\"utkebhomert \cite{BoschLutke-DAV}. We also recall the nonarchimedean Appell--Humbert theory. Our main reference is \cite{BoschLutke-DAV}. We also refer to \cite{FRSS} and \cite{KY}

We begin by recalling the valuation map for an algebraic torus. Let $M$ be a free $\ZZ$-module and set $N_{\RR} \colonequals \mathrm{Hom}(M , \RR)$. Let $T = \Spec (k [\chi^M])$ be the algebraic torus over $k$ with character lattice $M$. Then each $x \in T^{\an}$ defines a homomorphism
$M \ni u \mapsto -\log |\chi^{u} (x)| \in \RR$, and this gives rise to a continuous surjective map 
\begin{equation} 
\label{eqn:valuation:map}
\val\colon T^{\an} \to N_{\RR} = \mathrm{Hom}(M , \RR),
\end{equation}
called the \emph{valuation map}.
This map restricts to a group homomorphism $T(k) \to N_{\RR}$.

Let $A$ be an abelian variety over $k$. We say that $A$ is \emph{totally degenerate} if it is unifomized by an algebraic torus, i.e., there exists a homomorphism $p\colon T^{\an} \to A^{\an}$ of analytic groups that is a universal cover of the pointed topological space $(A^{\an},0)$, where $T$ is an algebraic torus. Such a homomorphism $p$ is unique up to canonical isomorphism, called \emph{uniformization} of $A^{\an}$. Let $M$ be the character lattice of $T$ and $N_{\RR} := \Hom (M,\RR)$. The uniformization $p$ has the following properties: $\Ker (p) \subseteq T(k)$; $M' \colonequals \val (\Ker (p))$ is a lattice (i.e., a full lattice) in $N_{\RR}$; $\rest{\val}{\Ker (p)}\colon \Ker (p) \to M'$ is a group isomorphism.  

Assume that $A$ is totally degenerate. In the sequel, let $p\colon T^{\an} \to A^{\an}$ denote the uniformization, where $T$ is an algebraic torus with character lattice $M$, and let $\val: T^{\an} \to N_{\RR}$ denote the valuation map, where $N_{\RR} := \mathrm{Hom} (M,\RR)$. Let $\widetilde{\Phi}\colon M' \to T^{\an}$ be the group homomorphism given by the inverse of $\rest{\val}{\Ker (p)}\colon \Ker (p) \to M'$. Then $M'$ acts on $T^{\an}$ via $\widetilde{\Phi}$, and we have $A^{\an} = T^{\an}/M'$. The valuation map $\val\colon T^{\an} \to N_{\RR}$ descends to a continuous map 
\begin{equation}
\label{eqn:val;A:an}
\val_{A^{\an}}\colon A^{\an} \to N_{\RR} / M' , 
\end{equation}
and we have a commutative diagram
\[
\begin{CD}
0 @>>>  \Ker(p) @>>> T^{\an} @>{p}>> A^{\an} @>>> 0\\
@. @VV{\cong}V @VV{\val}V @VV{\val_{A^{\an}}}V @. \\
0 @>>>  M^\prime @>>> N_\RR @>>> N_\RR/M^\prime @>>> 0 , \\
\end{CD}
\]
where the horizontal lines are exact.
We call $N_{\RR}/M'$ the \emph{tropicalization} of $A$ or $A^{\an}$.
Since $M'$ is a full lattice in $N_{\RR}$, $N_{\RR}/M'$ is a real torus with an integral structure $N$. We will see in Corollary~\ref{cor:val:tav} that $N_{\RR}/M'$ is in fact a tropical abelian variety, i.e., it has a polarization. 

We set 
$T' \colonequals \Spec (k [\chi^{M'}])$.
Given $u \in M$, we have a $k$-algebra homomorphism
\[
k [\chi^{M'}] \to k, \;\;  
\chi^{u'} \mapsto \chi^{u} (\widetilde{\Phi} (u')) 
\qquad \text{(for any $u^\prime \in M^\prime$)}. 
\]
This gives a $k$-valued point of $T^\prime$, denoted by 
$\widetilde{\Phi}'(u)$. Then we have a homomorphism 
$\widetilde{\Phi}^\prime\colon M \to T^\prime(k)$. 
It is straightforward to see that 
$\chi^{u} (\widetilde{\Phi} (u')) 
= \chi^{u'} (\widetilde{\Phi}' (u))$. 
We define a bilinear form $t\colon M' \times M \to k^{\times}$ by 
\begin{equation}
\label{eqn:def:t}
t\colon M' \times M \to k^{\times}, \quad 
(u' , u) \mapsto \chi^u (\widetilde{\Phi} (u')) = \chi^{u'} (\widetilde{\Phi}'(u)). 
\end{equation}

The following lemma is essentially a special case of \cite[Proposition~3.7]{BoschLutke-DAV}.

\begin{Lemma}
\label{lemma:symmetric:form}
We keep the above notation. Let $\lambda\colon M' \to M$ be a homomorphism. Then the bilinear form $M' \times M' \to k^{\times}$ given by $(u_1' , u_2') \mapsto t(u_1' , \lambda (u_2'))$ is symmetric.
\end{Lemma}

\Proof
Let $\Lambda\colon T \to T'$ be the morphism induced by $\lambda\colon M' \to M$. Then for any $u_1' , u_2' \in M'$, we have
\begin{align*}
t (u_1' , \lambda (u_2')) 
&= \chi^{\lambda (u_2')} (\widetilde{\Phi} (u_1')) =
\chi^{u_1'} (\widetilde{\Phi}' (\lambda (u_2'))) 
& \text{(by \eqref{eqn:def:t})}
\\
& = \chi^{u_1'} (\Lambda (\widetilde{\Phi} (u_2'))) =
\chi^{\lambda (u_1')} (\widetilde{\Phi} (u_2'))) 
= t (u_2' , \lambda (u_1')), & 
\end{align*}
which completes the proof. 
\QED

Next, we recall the theory of line bundles over totally degenerate abelian varieties. We begin by  recalling the notion of rigidified line bundles. 
Let $X$ be an algebraic variety over $k$, i.e., a reduced and irreducible scheme that is separated and of finite type over $k$. 
We fix a reference point $0 \in X(k)$. 
A \emph{rigidification} for a line bundle $L$ on $X$ means an isomorphism $\rest{L}{0} \cong k$. A {\em rigidified line bundle} is 
a line bundle with a rigidification. When $X$ is an algebraic group, we take the identity element of the algebraic group as the reference point. We similarly 
define rigidified line bundles over $X^{\an}$. 
Suppose that $X$ is proper over $k$. Then since $H^{0} (X,k) = k$, a homomorphism of rigidified line bundles is unique, and hence we will often use equality ``$=$'' in place of ``$\cong$.''

Let $A$ be a totally degenerate abelian variety. Let $p\colon T^{\an} \to A^{\an}$, $M$, $N_{\RR}$, $\val\colon T^{\an} \to N_{\RR}$, and $M'$ be as before.
A \emph{nonarchimedean descent datum} is 
a pair $(\lambda, c)$ consisting of a homomorphism $\lambda\colon M' \to M$ and a map $c\colon M' \to k^{\times}$ such that  
\begin{equation}
\label{eqn:c:condition}
c(u_1' + u_2')\, c(u_1')^{-1}\, c(u_2')^{-1} = t (u_1' , \lambda (u_2'))
\qquad (\text{for any $u_1' , u_2' \in M'$}). 
\end{equation}

Let $(\lambda, c)$ be a nonarchimedean descent datum. The datum defines an $M^\prime$-action on the trivial (rigidified) line bundle on $T^{\an}$ characterized by the condition that
\begin{align} \label{align:M'-linearlization}
(x, t) \mapsto \big(\widetilde{\Phi}(u^\prime) \, x, \; c(u')^{-1} \chi^{- \lambda (u')}\, t\big)
\end{align}
for $(x, t) \in T(k) \times k$ and $u' \in M'$.
Taking the quotient, we obtain a rigidified line bundle 
on $A^{\an}$. By GAGA (see \cite[Chap.~3]{Berkovich}), there exists a rigidified line bundle $\OO_A (\lambda, c)$ on $A$ whose analytification coincides with that line bundle
on $A^{\an}$.

We define a bilinear form 
$t_{\trop}\colon M' \times M \to \RR$ by
\[
t_{\trop}\colon
(u',u) \mapsto - \log |t(u',u)|. 
\]
Let $\langle \ndot , \ndot \rangle\colon M \times N_{\RR} \to \RR$ denotes the canonical pairing.
The following lemma is a special case of \cite[Lemma~9.1]{KY}.

\begin{Lemma}
\label{lemma:DD-tropicalization}
We have $t_{\trop} (u' ,u) = \langle u ,u' \rangle$ for any $(u' ,u) \in M' \times M$. 
\end{Lemma}

\Proof
Let $(u' ,u) \in M' \times M$. Note that $\tilde{\Phi}$ is defined by $\val (\tilde{\Phi} (u')) =u'$ and $\val$ is defined by $\langle u ,\val (x) \rangle = - \log |\chi^u (x)|$. Then 
\[
t_{\trop} (u' ,u) 
\colonequals - \log |t(u', u)|
= - \log |\chi^{u} (\widetilde{\Phi} (u'))| 
= \langle u, \val(\tilde{\Phi}(u^\prime))\rangle
= \langle u , u' \rangle, 
\]
as desired. 
\QED

We recall the nonarchimedean Appell--Humbert theory. For any homomorphism $\lambda :M' \to M$, we define an $\RR$-bilinear form $Q_{\lambda}\colon N_{\RR} \times N_{\RR} \to \RR$ by $Q_{\lambda} (u_1' , u_2') = t_{\trop} (u_1' , \lambda (u_2'))$ for any $u_1' , u_2' \in M' \subseteq N_{\RR}$. By Lemma~\ref{lemma:symmetric:form}, $Q_\lambda$ is symmetric. By Lemma~\ref{lemma:DD-tropicalization}, we note that $\lambda$ and $Q_{\lambda}$ correspond to each other in Lemma~\ref{lem:Q:ell:lambda:gamma}.

\begin{Theorem} [Nonarchimedean Appell--Humbert theory~{\cite[Proposition~6.7]{BoschLutke-DAV}}] 
\label{thm:NAH}
Let $A$ be a totally degenerate abelian variety over $k$. 
\begin{enumerate}
\item
Let $(\lambda_1 , c_1)$ and $(\lambda_2, c_2)$ be 
nonarchimedean descent data. Then so is $(\lambda_1 + \lambda_2 , c_1 c_2)$, and we have 
\[
\OO_{A} (\lambda_1 , c_1) \otimes \OO_{A} (\lambda_2 , c_2) = \OO_{A} (\lambda_1 +\lambda_2 , c_1 c_2)
.
\]
Further, $\OO_{A} (\lambda_1 , c_1) = \OO_{A} (\lambda_2 , c_2)$ if and only if $\lambda_1 = \lambda_2$ and $c_2 =  \big( \chi^u \circ \widetilde{\Phi} \big) c_1$ for some $u \in M$. 
\item
For any line bundle $\tilde{L}$ on $A$, there exists a descent datum $(\lambda, c)$ such that $\tilde{L} \cong \OO_A(\lambda, c)$. 
Furthermore, $\tilde{L}$ is ample if and only if the symmetric $\RR$-bilinear form $Q_{\lambda}$
is positive-definite.
\end{enumerate}
\end{Theorem}

\begin{Corollary}
\label{cor:val:tav}
The tropicalization $N_{\RR}/M'$ of $A$
is a tropical abelian variety. 
\end{Corollary}

\Proof
Since $A$ is an abelian variety, Theorem~\ref{thm:NAH}~(2) implies that there exists a descent datum $(\lambda, c)$ such that the corresponding symmetric $\RR$-bilinear form $Q_{\lambda}$ is positive-definite. By Lemmas~\ref{lemma:symmetric:form} and~\ref{lemma:DD-tropicalization}, for any $u' \in M^\prime$ and $v \in N$,  we have 
\[
Q_{\lambda} (u', v) = Q_{\lambda} (v, u') = 
t_{\trop} (\lambda(u'), v) = \langle \lambda(u'), v \rangle \in \ZZ. 
\]
which means that $Q_{\lambda}$ is integral over $M' \times N$. Thus
$Q_{\lambda}$ is a polarization on $N_{\RR}/M'$, which proves that  $N_{\RR}/M'$ is a tropical abelian variety. 
\QED

Finally in this subsection, we see how the involution acts on the line bundles given by nonarchimedean descent data. Let $[m]\colon A \to A$ denote the $m$-times endomorphism on the abelian variety $A$ for an $m \in \ZZ$.

\begin{Lemma} \label{lemma:pullback:involution}
Let $m$ be an integer. For any descent datum $(\lambda , c)$, we have 
\[
[m]^{\ast} (\OO_{A} (\lambda ,c)) = \OO_A \left( m^2 \lambda , 
t (\ndot , \lambda (\ndot))^{\frac{m(m-1)}{2} }
c^{m} \right),
\] where $t (\ndot , \lambda (\ndot))$ denotes the function $M' \to k^{\times}$ given by 
$u' \mapsto t (u' , \lambda (u'))$.
\end{Lemma}

\Proof
The line bundle $\OO_A (\lambda ,c)^{\an}$ is constructed to be the quotient of the trivial line bundle by the $M'$-action given in (\ref{align:M'-linearlization}).
Pulling back this $M'$-action by the homomorphism $t \mapsto t^m$, we obtain an 
$M'$-action on the trivial line bundle characterized by
\[
1  \mapsto c(u')^{-m}  \left( \chi^{\lambda(u')} \left( \widetilde{\Phi} (u') \right)^{\frac{m(m-1)}{2}}  \right)^{-1} \chi^{-m^2\lambda (u')} 
\]
for any $u' \in M'$. Note that $\chi^{\lambda(u')} \left( \widetilde{\Phi} (u')\right) ^{\frac{m(m-1)}{2}} = t (u' , \lambda (u'))^{\frac{m(m-1)}{2} }$ for any $u' \in M'$. Then, passing through the quotient construction, we see that the required equality holds.
\QED

\begin{Proposition} \label{prop:symmetirc:descentdatum}
Let $\tilde{L}_1$ be a line bundle on $A$. Then there exists a descent datum $(\lambda , c)$ that gives $\tilde{L}_1 \otimes [-1]^{\ast} (\tilde{L}_1)$ and such that $c$ is even, i.e., $c(-u') = c(u')$ for any $u' \in M'$.
\end{Proposition}

\Proof
There exists a descent datum $(\lambda_1 , c_1)$ such that $\tilde{L}_1 \cong \OO_{A} (\lambda_1 , c_1)$. We set $\lambda := 2 \lambda_1$ and define $c\colon M' \to k^{\times}$ by $c(u') = t (u' , \lambda_1 (u'))
$. It is straightforward to see that $(\lambda , c)$ is a nonarchimedean descent datum with $c(-u') = c(u')$ and that $(\lambda , c) = (2 \lambda_1 , ( t (\ndot , \lambda_1 (\ndot)) c_1^{-1}  ) c_1)$.
By Lemma~\ref{lemma:pullback:involution} for $m=-1$, we have $[-1]^{\ast} (\OO_A (\lambda_1 , c_1)) = \OO_A (\lambda ,t (u' , \lambda_1 (u')) c_1^{-1})$. Noting Theorem~\ref{thm:NAH}~(1), we obtain the conclusion.
\QED

\subsection{Nonarchimedean theta functions and their tropicalizations}
In this subsection, we recall the tropicalization of nonarchimedean theta functions and show some properties.

Let us begin by recalling nonarchimedean theta functions on an algebraic torus. 
Let $(\lambda, c)$ be a nonarchimedean descent datum. 
A function $f \in H^{0} (T^{\an} , \OO_{T^{\an}})$ is called a 
\emph{nonarchimedean theta function} with respect to $(\lambda, c)$ 
if $f$ satisfies the \emph{quasi-periodicity} condition with respect to $(\lambda, c)$, i.e., $f$ satisfies the condition that
$T_{\widetilde{\Phi}(u')}^{*} (f) = c(u')^{-1} \chi^{- \lambda (u')} f$
for any $u' \in M'$, 
where $T_{\widetilde{\Phi}(u')}\colon T^{\an} \to T^{\an}$ is the translation by $\widetilde{\Phi}(u') \in T (k) \subseteq T^{\an}$. 

We take a rigidified line bundle $\OO_{A} (\lambda, c)$ on $A$. 
By its construction and GAGA (see \cite[Chap.~3]{Berkovich}), we have a canonical inclusion $H^0 (A , \OO_{A} (\lambda, c)) \subseteq H^0 (T^{\an} , \OO_{T^{\an}})$. Via this inclusion, $H^0 (A , \OO_{A} (\lambda, c))$ equals the set of nonarchimedean theta functions with respect to~$(\lambda, c)$.

In \cite{FRSS}, Foster, Rabinoff, Shokrieh, and Soto define the tropicalization of a nonarchimedean theta function. Let $f$ be a nonarchimedean theta function with respect to $(\lambda, c)$. The \emph{tropicalization} of $(\lambda, c)$ is defined to be the pair $(\lambda , c_{\rm trop})$, where $c_{\rm trop}\colon M' \to \RR$ is a map given by $c_{\rm trop} (u') = - \log |c(u')|$. By Lemma~\ref{lemma:DD-tropicalization}, we have
\begin{align} \label{align:tropicalization:NADD}
c_{\trop} (u_1' + u_2') - c_{\trop} (u_1') - c_{\trop} (u_2') = \langle \lambda (u_2') , u_1' \rangle
\end{align}
for any $u_1' , u_2' \in M'$,
which shows that 
$(\lambda , c_{\rm trop})$ is indeed a tropical descent datum (cf. Definition~\ref{def:trop:descent:datum}). 

\begin{Remark} \label{remark:tropicalization:even}
Let $(\lambda , c)$ be a  nonarchimedean descent datum. Suppose that $c(u') = c(-u')$ for any $u' \in M'$. Then the tropical descent datum $(\lambda , c_{\trop})$ satisfies that  $c_{\trop} (u') = c_{\trop} (-u')$ for any $u' \in M'$. It follows from (\ref{align:tropicalization:NADD}) that
$c_{\trop}(u^\prime) = \frac{1}{2} \langle u^\prime, \lambda(u^\prime) \rangle$ for any $u' \in M'$. Further, we see that if $(Q_{\lambda} , \ell)$ corresponds to $(\lambda , c_{\trop})$ in Lemma~\ref{lem:Q:ell:lambda:gamma}, then $\ell = 0$.
\end{Remark}

For a function $f \in H^{0} (T^{\an} , \OO_{T^{\an}})$, we have a unique Fourier expansion
\begin{align*} 
f = \sum_{u \in M} a_u \chi^u 
\qquad (a_u \in k)
,
\end{align*}
and we define the tropicalization of $f$ to be
\begin{equation}
\label{eqn:f:tropicalization}
f_{\trop} (v) \colonequals \min_{u \in M} ( \langle u , v\rangle - \log |a_u| )
\end{equation}
as a function on $v \in N_{\RR}$.
Indeed, the minimum exists because of the convergence of the Fourier expansion.
Suppose that $f$ is a nonarchimedean theta function with respect to $(\lambda , c)$. Then by \cite[Theorem~4.7]{FRSS}, $f_{\trop} \in H^0 (X , L (\lambda , c_{\trop}))$, where $X = N_{\RR} / M'$ and see Remark~\ref{remark:notation:H0}
for the notation.
In this setting, we call $f$ a \emph{nonarchimedean lift} of $f_{\trop}$.

Assume now that $\OO_A (\lambda, c)$ is ample, or equivalently, by Theorem~\ref{thm:NAH}~(2), assume that the corresponding symmetric $\RR$-bilinear form $Q_{\lambda}\colon N_{\RR} \times N_{\RR} \to \RR$ is positive-definite. It follows that for any $b \in M$, 
\[
\min_{u' \in M'} ( \langle b +\lambda (u') , v \rangle + c_{\rm trop} (u') + \langle b, u' \rangle )
\]
exists for any $v$. 
Using this, one shows that for $b \in M$, we see that the infinite sum
\begin{align} \label{def:Suminonarchitheta}
\tilde{\vartheta}_{b} :=
\sum_{u' \in M'} \chi^{b + \lambda (u')}  c(u') \chi^b \left( \widetilde{\Phi} (u') \right)
\end{align}
converges in $H^0 (T^{\an} , \OO_{T^{\an}})$ and that this is in fact a nonarchimedean theta function with respect to $(\lambda , c)$. We refer to \cite[Section~9.2]{KY} for details.

The following lemma is essentially proved in the proof of \cite[Remark~5.4]{BoschLutke-DAV}.

\begin{Lemma} \label{lemma:basis-theta}
Under the above setting, let $\mathfrak{B} \subseteq M$ a complete system of representatives of $M/\lambda (M')$. Then $\{ \tilde{\vartheta}_{b} \}_{b \in \mathfrak{B}}$ is a basis of of $H^0 (A , \OO_A (\lambda, c))$.
\end{Lemma}

\Proof
We write $\mathfrak{B} = \{ b_0 , \ldots , b_{m-1}\}$, where $m = |\mathfrak{B}|$. By \cite[Remark~5.4]{BoschLutke-DAV}, we have $|\mathfrak{B}| = \dim_k (H^0 (A , \OO (\lambda , c)))$, and hence it suffices to show that $\tilde{\vartheta}_{b_0} , \ldots , \tilde{\vartheta}_{b_{m-1}}$ is linearly independent over $k$. 
Suppose that $\alpha_0 , \ldots , \alpha_{m-1} \in k$ and
$\sum_{i=0}^{m-1} \alpha_i \tilde{\vartheta}_{b_i} = 0$.
Fix any $i=0 , \ldots , m-1$. Since $\{ b_0 , \ldots , b_{m-1}\} \subseteq M$ is a complete system of representatives of $M/\lambda (M')$, the coefficients at $b_i$ of the Fourier expansion of the left-hand side comes only from $\tilde{\vartheta}_{b_i}$, and it equals $\alpha_i c(0) \chi^{b_i} (0)$. It follows from the above equality that this coefficient equals $0$. Thus $\alpha_i=0$. This completes the proof.
\QED

A direct computation shows that $\left( \tilde{\vartheta}_{b} \right)_{\trop} \in H^0 (X , L(\lambda ,c_{\trop}))$. More precisely, we have the following.

\begin{Proposition} [cf. Corollary~9.8 in \cite{KY}] \label{prop:lifting:constant}
Let $(\lambda, c)$ be a nonarchimedean descent datum. Suppose that $c(u') = c(-u')$ for any $u' \in M'$. 
Let $b \in M$. Then there exists a $\tau \in T(k)$ such that $\lambda (\val (\tau^2) ) = b$. Further, we have
\[
\left( \chi^b (\tau) \tilde{\vartheta}_b \right)_{\trop} = \vartheta_b,
\]
where $\vartheta_b$ is the tropical theta function with respect to $(Q_{\lambda},0)$ defined in (\ref{eqn:def:theta:b}) (cf. Remark~\ref{remark:tropicalization:even}).
In particular, $\vartheta_b$ has a nonarchimedean lift. Further, if $\mathfrak{B} \subseteq M$ is a complete system of representatives of $M / \lambda (M')$, then $\{ \chi^b (\tau) \tilde{\vartheta}_b\}_{b \in \mathfrak{B}}$ is a basis of $H^0 (A , \OO_A (\lambda , c))$.
\end{Proposition}

\Proof
Since the group $T(k)$ is divisible, there exists a $\tau \in T(k)$ such that $\lambda (\val (\tau^2) ) = b$, which shows the first assertion. The second assertion 
can be checked by a direct computation. Indeed, 
for any $v \in N_\RR$, we have
{\allowdisplaybreaks
\begin{align*}
& \left( \chi^b (\tau) \tilde{\vartheta}_b \right)_{\trop}(v) \\
& \quad = 
\langle b, \val(\tau) \rangle + 
\min_{u' \in M'} ( \langle b +\lambda (u') , v \rangle + c_{\rm trop} (u') + \langle b, u' \rangle )
&&  \text{(by \eqref{eqn:f:tropicalization} and \eqref{def:Suminonarchitheta})} \\
& \quad 
= \frac{1}{2} \langle b, \lambda_\RR^{-1}(b) \rangle + 
\min_{u' \in M'} ( \langle b +\lambda (u') , v \rangle + \frac{1}{2} \langle u^\prime, \lambda(u^\prime) \rangle + \langle b, u' \rangle )
&& \text{(by Remark~\ref{remark:tropicalization:even})} \\
& \quad = 
\min_{u' \in M'} \left( 
Q_{\lambda}(u^\prime + \lambda^{-1}_\RR(b), v) 
+ \frac{1}{2} Q_{\lambda}(u^\prime + \lambda^{-1}_\RR(b), u^\prime + \lambda^{-1}_\RR(b)) \right)
&&  \text{(by $Q_{\lambda}(u_1^\prime, u_2^\prime) =  \langle u_1^\prime, \lambda(u^\prime_2) \rangle$)} \\
& \quad = \vartheta_b(v), && 
\end{align*}
}
as desired. 

The last assertion immediately follows from  
the second assertion and 
Lemma~\ref{lemma:basis-theta}.
\QED

\section{Faithful tropicalization of the canonical skeleton of a Kummer surface}
\label{sec:faithful:nonarch:Kummer}
In this section, we define the canonical skeleton of a Kummer variety. Further, we prove a theorem on faithful tropicalization for Kummer surfaces.

As in Section~\ref{sec:uniformization}, let $k$ be an algebraically closed field that is complete with respect to a nontrivial nonarchimedean absolute value $|\ndot|_k$, namely, an algebraically closed, nontrivially valued, and complete nonarchimedean field, unless otherwise specified. In this section, we further assume that the residue characteristic of $k$ is not equal to $2$. 

\subsection{Canonical skeleton of a Kummer variety}
In this section, we define the canonical skeleton of the Kummer variety associated to an abelian variety. It will be constructed as the image of the canonical skeleton of the abelian variety by the analytification of the quotient morphism, and we will see that it will equal the quotient of the canonical skeleton of the abelian variety by the involution group. Further, we define the notion of faithful tropicalization of the canonical skeletons of Kummer varieties.

We begin by recalling the notion of the Kummer variety associated to an abelian variety.
Let $A$ be an abelian variety over $k$ of dimension $g$. Let $\langle -1 \rangle$ denote the cyclic group of order $2$ with generator $-1$, which we call the \emph{involution group}. The involution group acts on an abelian variety such that $-1$ corresponds to $[-1]\colon A \to A$.
The \emph{Kummer variety} associated to $A$ is defined to be the quotient $q\colon A \to \tilde{Y} \colonequals A / \langle -1 \rangle$ of $A$. Then the 2-torsion points are the fixed points by the involution consisting of  $4^{\dim (A)}$ points. When $g = 2$, the Kummer variety is called the \emph{Kummer surface}.

In the sequel, assume that $A$ is totally degenerate. Let $T$ be an algebraic torus with character lattice $M$ of rank $g$ and with cocharacter $N$. Let $p\colon T^{\an} \to A^{\an}$ be the uniformization and let $\val\colon T^{\an} \to N_{\RR}$ be the valuation map. We take a lattice $M' \subseteq N_{\RR}$ and a homomorphism $\widetilde{\Phi}\colon M' \to T(k)$ such that $M' = \val (\Ker (p))$, as 
in Section~\ref{sec:uniformization}.

We recall the canonical skeleton of $A$. See \cite[Example~5.2.12 and Theorem~6.5.1]{Berkovich}, \cite[Example 7.2]{Gubler-NACM}, \cite[Section~4.1]{FRSS}, and \cite[Section~8.1]{KY} for details.
For each $v \in N_{\RR}$, we define a map $\sigma (v)\colon k [\chi^M] \to\RR_{\geq 0}$ by
$\sum_{u \in M} a_u \chi^u \mapsto \max_{u \in M}  |a_u| e^{- \langle u , v \rangle} $.
Then $\sigma (v)$ is a multiplicative norm on $k[\chi^M]$ extending $|\ndot|_k$ and hence gives a point in $T^{\an}$. This assignment defines a map $\sigma\colon N_{\RR} \to T^{\an}$, which is a continuous section of $\val\colon T^{\an} \to N_{\RR}$ and induces an homeomorphism from $N_{\RR}$ to its image. We call $\sigma$ the \emph{canonical section} of $\val$ and its image $\Sigma \colonequals \sigma (N_{\RR})$ the \emph{canonical skeleton} of $T^{\an}$. One sees that $\sigma$ is $M'$-equivariant and thus descends to a section $\sigma_{A^{\an}}\colon N_{\RR}/M' \to A^{\an}$ of $\val_{A^{\an}}\colon A^{\an} \to N_{\RR}/M'$. We call $\sigma_{A^{\an}}$ the \emph{canonical section} of $\val_{A^{\an}}$. Further, the image of $\sigma_{A^{\an}}$ equals $\Sigma / M'$, where $\Sigma \subseteq T^{\an}$ is the canonical skeleton of $T^{\an}$. We call 
\begin{equation}
\label{eqn:can:skeleton:abel}
\Sigma_{A^{\an}} \colonequals \Sigma / M' = \sigma_{A^{\an}} (N_{\RR}/M^\prime)
\end{equation}
the \emph{canonical skeleton} of~$A^{\an}$.

Now, we consider Kummer varieties.
Let $\tilde{Y} \colonequals A/ \langle -1 \rangle$ be the Kummer variety associated to $A$ and let $q\colon A \to \tilde{Y}$ denote the quotient morphism. Let $\Sigma_{A^{\an}}$ be the canonical skeleton of $A^{\an}$. We set
\[
\Sigma_{\tilde{Y}^{\an}} \colonequals q^{\an} (\Sigma_{A^{\an}}),
\]
where $q^{\an}\colon A^{\an} \to \tilde{Y}^{\an}$ is the analytification of $q$. This is a compact subspace of $\tilde{Y}^{\an}$. 

We want to show that $\Sigma_{\tilde{Y}^{\an}}$ is the quotient of $\Sigma_{A^{\an}}$ by the involution group. To do that, we use the following lemma.

\begin{Lemma} \label{lemme:finite:quotient:Berkovich}
Let $G$ be a finite group acting on a quasi-projective variety $X$ over $k$ and let $X \to Y$ be the quotient of $X$ by $G$. Let $X^{\an} \to Y^{\an}$ be its analytification. Then we have a natural map $|X^{\an}|/G \to |Y^{\an}|$ of topological spaces, where for an analytic space $\mathbf{X}$, we let $|\mathbf{X}|$ denote the underlying topological space. Further, this map is a homeomorphism.
\end{Lemma}

\Proof
Since an automorphism of an algebraic variety naturally induces an automorphism of the analytification, $G$ acts naturally on $X^{\an}$. Further, the map $X^{\an} \to X$ given by $(p, |\ndot|) \mapsto p$ is $G$-equivariant. 

Let $q\colon X \to Y$ denote the quotient. Then $q^{\an}\colon X^{\an} \to Y^{\an}$ is a morphism of Berkovich analytic spaces. 
Let $|q^{\an}|\colon |X^{\an}| \to |Y^{\an}|$ denote the continuous map between the underlying topological spaces. We have a $G$-action on $|X^{\an}|$. Let $G$ acts on $|Y^{\an}|$ trivially. Since $|q^{\an}|$ is $G$-equivariant and the $G$-action on $|Y^{\an}|$ is trivial, the continuous map $|q^{\an}|$ factors through a unique continuous map $|X^{\an}| / G \to Y^{\an}$.

We take any $y \in Y^{\an}$. We prove that $|q^{\an}|^{-1}(y)$ is a $G$-orbit in $|X^{\an}|$. By (*) in the proof of Proposition~3.4.6 in \cite{Berkovich}, the canonical morphism
\[
(q^{\an})^{-1} (y) \to \left( X \times_{Y} \Spec (\mathscr{H} (y)) \right)^{\an}
\]
is an isomorphism, and in particular, the map 
\begin{align} \label{align:bijection:fibers}
\left| (q^{\an})^{-1} (y) \right|
\to 
\left| 
\left( X \times_{Y} \Spec (\mathscr{H} (y)) \right)^{\an}
\right|
\end{align}
is a homeomorphism. 
Note that the action of $G$ on $|X^{\an}|$ restricts to an action on $|(q^{\an})^{-1} (y)|$ and the action of $G$ on $X$ induces an action on $\left( X \times_{Y} \Spec (\mathscr{H} (y)) \right)^{\an}$. Note also that the continuous map in (\ref{align:bijection:fibers}) is $G$-equivariant.
As is mentioned in \cite[Definition/Exercise~4.1.7.1.~(ii)]{Temkin}, we have
\[
\left| (q^{\an})^{-1} (y) \right| = |q^{\an}|^{-1} (y).
\]
Since the morphism $X \to Y$ is finite, the natural map
\[
\left( X \times_{Y} \Spec (\mathscr{H} (y)) \right)^{\an}
\to X \times_{Y} \Spec (\mathscr{H} (y)) 
\]
is a $G$-equivariant bijection of finite sets.
It follows that the natural map
\begin{align} \label{align:G-bijection}
|q^{\an}|^{-1} (y)  \to X \times_{Y} \Spec (\mathscr{H} (y)) 
\end{align}
is a $G$-equivariant bijection of finite $G$-sets. Since $X \times_{Y} \Spec (\mathscr{H} (y))$ is a $G$-orbit of the $G$-action on $X(\mathscr{H} (y))$, the action of $G$ on $X \times_{Y} \Spec (\mathscr{H} (y))$ is transitive. Via the $G$-bijection in (\ref{align:G-bijection}), the $G$-action on $|q^{\an}|^{-1} (y)$ is transitive. Since the $G$-action on $|q^{\an}|^{-1} (y)$ is induced by the $G$-action on $|X^{\an}|$ by restriction, this proves that $|q^{\an}|^{-1} (y)$ is a $G$-orbit in $|X^{\an}|$.
 
It follows from what we have just proved that the induced map $|X^{\an}| / G \to |Y^{\an}|$ is injective. Since this map is surjective, it is bijective. Since $X \to Y$ is finite, $X^{\an} \to Y^{\an}$ is a proper map by \cite[3.4.7~Proposition]{Berkovich}. It follows that the continuous bijective map $|X^{\an}|/G \to |Y^{\an}|$ is proper. Thus it is a homeomorphism. This completes the proof.
\QED

We use the notation and the setting before Lemma~\ref{lemme:finite:quotient:Berkovich}. Using this lemma, we get the following.

\begin{Proposition} \label{prop:canonicalskeleton:kummer}
The action of $\langle -1 \rangle$ on $A^{\an}$ preserves the canonical skeleton $\Sigma_{A^{\an}} \subseteq A^{\an}$. Further, the induced map
\[
\Sigma_{A^{\an}} / \langle -1 \rangle \to \Sigma_{\tilde{Y}^{\an}}
\]
by $q^{\an}\colon A^{\an} \to \tilde{Y}^{\an}$ is a homeomorphism.
\end{Proposition}

\Proof
It is obvious from the construction of $\Sigma_{A^{\an}}$ that the involution $[-1]$ preserves $\Sigma_{A^{\an}}$. By Lemma~\ref{lemme:finite:quotient:Berkovich}, the map $q^{\an}\colon A^{\an} \to \tilde{Y}^{\an}$ factors through the quotient $|A^{\an}| \to |A^{\an}| / \langle -1 \rangle$ as a topological space and the induced map $|A^{\an}| / \langle -1 \rangle \to \tilde{Y}^{\an}$ is a homeomorphism. By restriction, $\Sigma_{A^{\an}} / \langle -1 \rangle \to \Sigma_{\tilde{Y}^{\an}}$ is a homeomorphism.
\QED

We set $Y \colonequals (N_{\RR}/M')/\langle -1 \rangle$, which is a tropical Kummer variety (see Section~\ref{subsec:trop:Kummer:orbifold:rps}). 
Since $\sigma_{A^{\an}}\colon N_{\RR} /M' \to \Sigma_{A^{\an}}$ is $\langle -1 \rangle$-equivariant, 
there exists a unique continuous map $\sigma_{\tilde{Y}^{\an}}\colon Y \to \tilde{Y}^{\an}$ such that the diagram
\[
\begin{CD}
A^{\an} @>{q^{\an}}>> \tilde{Y}^{\an} \\
@A{\sigma_{A^{\an}}}AA @AA{\sigma_{\tilde{Y}^{\an}}}A \\
N_{\RR} / M' @>>> Y
\end{CD}
\]
commutes, where 
the bottom horizontal arrow is the canonical surjection. Then $\sigma_{\tilde{Y}^{\an}}$ restricts to a continuous map $Y \to \Sigma_{\tilde{Y}^{\an}}$, and since $\Sigma_{A^{\an}} = \sigma_{A^{\an}} (N_{\RR}/M')$, the map $Y \to \Sigma_{\tilde{Y}^{\an}}$ is a homeomorphism by Proposition~\ref{prop:canonicalskeleton:kummer}.
Recall that $Y$ is a rational polyhedral orbifold. We put a structure of 
rational polyhedral orbifold on $\Sigma_{\tilde{Y}^{\an}}$ via $\sigma_{\tilde{Y}^{\an}}$. 

\begin{Definition} [Canonical skeleton of a Kummer variety]
\label{def:can:skeleton:Kummer}
We call $\Sigma_{\tilde{Y}^{\an}}$  
the \emph{canonical skeleton} of $\tilde{Y}^{\an}$.
\end{Definition}

\subsection{Faithful tropicalization of the canonical skeleton of a Kummer surface}

In this subsection, we discuss the faithful tropicalization for Kummer varieties and prove the faithful tropicalization theorem for Kummer surfaces..

First, we define the notion of the faithful tropicalizations of the canonical skeletons of Kummer varieties. Let $A$ be an abelian variety over $k$. Assume that it is totally degenerate. Let $q\colon A \to \tilde{Y}$ be the associated Kummer variety.
Let $\Sigma_{\tilde{Y}^{\an}}$ be the canonical skeleton of $\tilde{Y}^{\an}$. 
Let 
\begin{equation}
\label{eqn:trop:morphism}
\trop\colon 
(\PP^n_k)^{\an} \to \TT\PP^n, \qquad x \mapsto (-\log|X_0(x)|: \cdots: -\log|X_n(x)|)
\end{equation}
be the tropicalization map, where $X_0, \ldots, X_n$ are the homogeneous coordinate functions on~$\PP^n_k$. 

\begin{Definition}[faithful tropicalization]
\label{def:faithful:tropicalization}
We say that a morphism $\tilde{\varphi}\colon \tilde{Y} \to \PP^n_k$ 
{\em faithfully tropicalizes} the canonical skeleton $\Sigma_{\tilde{Y}^{\an}}$ 
if 
\[
\rest{\trop\circ\tilde{\varphi}^{\an}}{\Sigma_{\tilde{Y}^{\an}}}\colon 
\Sigma_{\tilde{Y}^{\an}} \to \TT\PP^n
\]
is a faithful embedding of $\Sigma_{\tilde{Y}^{\an}}$, where 
$\tilde{\varphi}^{\an}\colon \tilde{Y}^{\an} \to (\PP^n_k)^{\an}$ is the morphism induced by $\tilde{\varphi}$ and via the canonical section $\Sigma_{\tilde{Y}^{\an}}$ is regarded as a tropical Kummer variety $(\Sigma/M^\prime)/\langle -1\rangle$ (see Equation~\eqref{eqn:can:skeleton:abel} and Proposition~\ref{prop:canonicalskeleton:kummer}). 
\end{Definition}

We recall the descent of symmetric line bundles on an abelian variety to the Kummer variety. Let $A$ be an abelian variety over $k$ and let $q\colon A \to \tilde{Y}$ be the associated Kummer variety.
Let $\tilde{L}$ be a rigidified symmetric line bundle on $A$. Then the involution group acts on $\tilde{L}$ uniquely so that the automorphism $[-1]\colon A \to A$ preserves the rigidification. Since the translation by a $2$-torsion point commutes with $[-1]$, it follows that the action of the involution preserves the fiber of $\tilde{L}$ over any $2$-torsion point. Hence $\tilde{L}$ descends to a unique rigidified line bundle $\tilde{L}_{\tilde{Y}}$ on $\tilde{Y}$, where the rigidification over $\tilde{Y}$ is taken at $q(0)$. 
We call $\tilde{L}_{\tilde{Y}}$ the \emph{descent} of $\tilde{L}$. There exists a canonical homomorphism $q^{\ast}\colon H^0 (\tilde{Y} , \tilde{L}_{\tilde{Y}}) \to H^0 (A , \tilde{L})$, which is injective.

Assume that $A$ is totally degenerate. 
We say that a line bundle $\tilde{L}_{\tilde{Y}}$ {\em admits a faithful tropicalization} for the
canonical skeleton $\Sigma_{\tilde{Y}^{\an}}$ if there exist $s_0, \ldots, s_n \in H^0(\tilde{Y}, \tilde{L}_{\tilde{Y}})$ that do not have common zero points (and hence $\tilde{L}_{\tilde{Y}}$ is basepoint free) and such that
$\tilde{\varphi}_{s_0, \ldots, s_n}\colon \tilde{Y} \to \PP^n_k$ defined by $y \mapsto (s_0(y):\cdots:s_n(y))$ 
faithfully tropicalizes $\Sigma_{\tilde{Y}^{\an}}$.

We define the polarization type of an ample line bundle.
Let $\tilde{L}$ be an ample line bundle on a totally degenerate abelian variety $A$. 
By Theorem~\ref{thm:NAH}, there exists a descent datum $(\lambda, c)$ such that $\tilde{L} \cong \OO_A (\lambda, c)$. 
Since $\lambda$ is uniquely determined from $\tilde{L}$, so is the corresponding symmetric  $\RR$-bilinear form $Q_{\lambda}$, and since $\tilde{L}$ is ample, $Q_{\lambda}$ is a polarization; see Theorem~\ref{thm:NAH}~(2). We call $Q_{\lambda}$ the \emph{polarization given by $\tilde{L}$}. Let $(d_1, d_2, \ldots, d_g) \in (\ZZ_{\geq 1})^g$ be the type of $Q_{\lambda}$,
where $g := \dim A$. 
We call $(d_1 , \ldots , d_g)$ the \emph{polarization type} of $\tilde{L}$.

\begin{Remark} \label{remark:give:PP}
Let $A'$ be the dual abelian variety of $A$. For an ample line bundle $\tilde{L}$, we define a homomorphism $\phi_{\tilde{L}} \colon A \to A'$ by $\phi_{\tilde{L}} (x) \colonequals T_x^{\ast} (\tilde{L}) \otimes \tilde{L}^{-1}$.
By \cite[Theorem~6.15]{BoschLutke-DAV}, $\tilde{L}$ gives a principal polarization in the above sense if and only if $\phi_L$ is an isomorphism.
\end{Remark}

Assume that $\tilde{L}_1$ is a symmetric ample line bundle on $A$ and gives a principal polarization. We set $\tilde{L} \colonequals \tilde{L}_1^{\otimes 2}$.
It is well known that $\tilde{L}$ is basepoint free and that $\dim_k (H^0 (A , \tilde{L} ))= 2^g$, where $g:= \dim (A)$. Further, every section in $H^0 (A , \tilde{L} )$ is $\langle -1 \rangle$-invariant (see \cite[\S4.8]{BL} when $k = \CC$, and \cite[\S2]{Mu66} for general, together with Remark~\ref{remark:give:PP}). 
Thus the canonical homomorphism $H^0 (\tilde{Y} , \tilde{L}_{\tilde{Y}} ) \to H^0 (A , \tilde{L} )$ is an isomorphism. In particular, $\tilde{L}_{\tilde{Y}}$ is also basepoint free, and a basis $s_0 , \ldots , s_{2^g-1}$ of $H^0 (\tilde{Y} , \tilde{L}_{\tilde{Y}} )$ induces a morphism $\tilde{\varphi}_{s_0 , \ldots , s_{2^g-1}}\colon \tilde{Y} \to \PP^{2^g-1}_k$.

To state our faithful tropicalization result for Kummer surfaces, 
we define the notion of tropically irreducibility for a totally degenerate abelian surface. We keep the assumption that $A$ is totally degenerate and
assume in addition that $\dim (A) = 2$.
Assume 
that a line bundle $\tilde{L}$ on $A$ gives a principal polarization $Q$. Let $X = N_\RR/M^\prime$ be the tropicalization of $A^\an$. Note that $(X,Q)$ is a principally polarized tropical abelian surface.
Then we say that the pair
$(A, \tilde{L})$ is {\em tropically of product type} (resp. {\em tropically irreducible}) if $(X, Q)$ is of product type (resp. irreducible); see Definition~\ref{def:irr}.

The following is one of our main results.

\begin{Theorem} \label{thm:main1:FT}
Let $A$ be a totally degenerate abelian surface over $k$. Let $\tilde{L}_1$ be a symmetric ample line bundle on $A$ and assume that it gives a principal polarization.
Assume that 
the pair $(A, \tilde{L}_1)$  
is tropically irreducible. Let $\tilde{Y} \colonequals A/\langle -1 \rangle$ be the Kummer surface associated to $A$ and let $\tilde{L}_{\tilde{Y}}$ on $\tilde{Y}$ be the descent of $\tilde{L} \colonequals \tilde{L}_1^{\otimes 2}$ on $A$. Then there exists a basis $s_0, \ldots, s_3$ of $H^0(\tilde{Y}, \tilde{L}_{\tilde{Y}})$ such that the induced morphism ${\varphi}_{s_0, \ldots, s_3}\colon \tilde{Y} \to \PP^3_k$ 
faithfully tropicalizes 
$\Sigma_{\tilde{Y}^{\an}}$. In particular, 
$\tilde{L}_{\tilde{Y}}$ admits a faithful tropicalization for the canonical skeleton of $\tilde{Y}^{\an}$
\end{Theorem}

\Proof
Let $p\colon T^{\an} \to A^{\an}$ be the uniformization and let $M$ be the character lattice of $T$.
Let $X=N_{\RR}/M'$ be the tropicalization of $A$. 
Since $\tilde{L}_1$ is symmetric, it follows from Proposition~\ref{prop:symmetirc:descentdatum} that there exists a nonarchimedean descent datum $(c,\lambda)$ such that $\tilde{L} \cong \OO_{A} (\lambda , c)$ and $c (u) = c (-u)$. 
Since the polarization $Q_{\lambda}$ corresponding to $\lambda\colon M' \to M$ is the twice of a principal polarization, $\lambda$ is the twice of an isomorphism, and hence $|M / \lambda (M')| = 4$. We fix a complete system $\{ b_0 , \ldots , b_3 \} \subseteq M$ of representatives of $M / \lambda (M')$. 
For each $i=0,\ldots , 3$, Proposition~\ref{prop:lifting:constant} gives us a basis $s_0 , \ldots , s_3$ of $H^0 (A , L)$ such that $(s_i)_{\trop} = \vartheta_{b_i}$, where $\vartheta_{b_i}$ is the tropical theta function with respect to $(Q_{\lambda} , 0)$ defined in (\ref{eqn:def:theta:b}).  

Let  $\sigma_{A^{\an}}\colon X \to \Sigma_{A^{\an}}$ be the canonical section of the valuation map $A^{\an} \to X$. 
Then $\sigma_{A^{\an}}$ induces a homeomorphism $Y \to \Sigma_{\tilde{Y}^{\an}}$, where $Y \colonequals X/\langle -1 \rangle$.
Then since $\tilde{\varphi}_{s_0 , \ldots , s_3}^{\trop} \circ \sigma_{\tilde{Y}^{\an}} = \varphi_{\vartheta_{b_0} , \ldots , \vartheta_{b_3}}$,
Theorem~\ref{thm:faith:embeddings} shows that $\tilde{\varphi}^{\trop}$, which is the composite of $\tilde{\varphi}$ with $\trop\colon (\PP^3_k)^{\an} \to \TT\PP^3$, restricts to a faithful embedding of $\Sigma_{\tilde{Y}^{\an}}$ into $\TT\PP^3$. This completes the proof.
\QED

Next, we consider a lifting of tropical Kummer quartic surfaces. 
We use the following proposition. 

\begin{Proposition}
\label{prop:lifting:abel}
Let $(X,Q)$ be a principally polarized tropical abelian surface. 
Then there exists a pair $(A,\tilde{L})$ of a totally degenerated abelian surface $A$ over some algebraically closed, nontrivially valued, and complete nonarchimdean field 
and an ample line bundle on $A$ such that the following hold: $\Sigma_{A^{\an}} \cong X$; there exists a nonarchimedean descent datum $(\lambda , c)$ such that $\tilde{L} \cong \OO_{A} (\lambda,c)$, $\lambda$ corresponds to $Q$, and such that $c_{\trop}(u') = \frac{1}{2} \langle \lambda (u') , u'\rangle$ for any $u' \in M'$.
\end{Proposition}

\Proof
We take an algebraically closed, nontrivially valued, and complete nonarchimedean field $k'$ that is a field extension of $k$ such that
the value group $\{-\log |\alpha| \mid \alpha \in k^\times\}$ of $k$ is equal to $\RR$. There exists such a valued field; see \cite{Ma}, for example. Replacing $k$ by $k'$, we may and do assume that $k' = k$. 

We write $X = N_\RR/M^\prime$ as before. Recall from Section~\ref{subsec:na:ah} that  $T = \Spec(k[\chi^M])$ and $N_{\RR} = \Hom (M,\RR)$. We have a $\ZZ$-linear map $\lambda\colon M' \to M$ corresponding to $Q$, i.e., $\lambda$ is characterized by the condition that  $Q(u_1^\prime, u_2^\prime) = \langle \lambda(u_1^\prime), u_2^\prime\rangle$ 
for $u_1^\prime, u_2^\prime \in M^\prime$. Since 
$Q$ is a principal polarization, $\lambda$ is an isomorphism.

It suffices to construct a group homomorphism $\tilde{\Phi}\colon M^\prime \to T(k) \subseteq T^{\an}$ and a map $c\colon M^\prime \to k^\times$ such that 
$\val(\tilde{\Phi}(u^\prime))= u^\prime$ for any $u^\prime \in M^\prime$, the bilinear form 
$M^\prime \times M^\prime \ni (u^\prime_1, u^\prime_2) \mapsto \chi^{\lambda(u^\prime_1)}(
\tilde{\Phi}(u^\prime_2))$ is symmetric, $c$ satisfies \eqref{eqn:c:condition}, and $c_{\trop} (u^\prime) = \frac{1}{2}  \langle \lambda(u^\prime), u^\prime\rangle$ for any $u^\prime \in M^\prime$. Indeed, we then have a proper commutative analytic group 
$T^{\an}/\tilde{\Phi}(M^\prime)$, and since $(\lambda, c)$ is a nonarchimedean descent datum, we have a line bundle $\OO_A(\lambda, c)$. Since the symmetric bilinear form $Q$ is positive-definite by the assumption, it follows from \cite[Theorem~2.4]{BoschLutke-DAV} that $\OO_A(\lambda, c)$ is ample. Thus $T^{\an}/\tilde{\Phi}(M^\prime)$ is an analytification of an abelian variety, and this is necessarily totally degenerate. 

We fix a $\ZZ$-basis $\mathsf{f}_1^\prime, \mathsf{f}_2^\prime$ of $M^\prime$ and set $\mathsf{e}_1 = \lambda(\mathsf{f}_1^\prime)$ and $\mathsf{e}_2 = \lambda(\mathsf{f}_2^\prime) \in M$. Since $\lambda$ 
is a $\ZZ$-linear isomorphism, $\mathsf{e}_1, \mathsf{e}_2$ is a $\ZZ$-basis 
of $M$. 
Since $\langle \mathsf{e}_i, \mathsf{f}_j^\prime\rangle \in \RR = \{ - \log |\alpha| \mid \alpha \in k^{\times} \} $, there exists a $2 \times 2$ symmetric matrix $\begin{pmatrix} \alpha_{ij} \end{pmatrix}_{i,j}$ with entries in $k^\times$ 
such that $-\log|\alpha_{ij}| = \langle \mathsf{e}_i, \mathsf{f}_j^\prime\rangle$. 
Since $(\chi^{\mathsf{e}_1}, \chi^{\mathsf{e}_2})\colon T \to \GG_{m, k}^2$ is an isomorphism 
of algebraic tori, there exist $t_1, t_2 \in T(k)$ such that 
$\chi^{\mathsf{e}_i}(t_j) = \alpha_{ij}$ for any $i, j \in \{1, 2\}$.

We define a group homomorphism $\tilde{\Phi}\colon M^\prime \to T(k) \subseteq T^{\an}$ by 
$\tilde{\Phi}(a \mathsf{f}_1^\prime + b \mathsf{f}_2^\prime) = t_1^a t_2^{b}$ 
for any $a, b \in \ZZ$. 
Then we have 
\begin{align*}
\val(\tilde{\Phi}(a \mathsf{f}_1^\prime + b \mathsf{f}_1^\prime))(\mathsf{e}_i)
& = \val(t_1^a t_2^b) (\mathsf{e}_i) = - \log\left|\chi^{\mathsf{e}_i}(t_1^a  t_2^b)\right|\\
&  = a (-\log |\alpha_{i1}|)  + b (-\log |\alpha_{i2}|)
= \langle \mathsf{e}_i, a \mathsf{f}_1^\prime + b \mathsf{f}_2^\prime \rangle, 
\end{align*}
so $\val(\tilde{\Phi}(u^\prime)) = u^\prime$ for any $u^\prime \in M^\prime$. Further, since  
$\chi^{\lambda(\mathsf{f}_i^\prime)}(\tilde{\Phi}(\mathsf{f}_j^\prime))
= \chi^{\mathsf{e}_i}(t_j) = \alpha_{ij}$, the bilinear form 
$M^\prime \times M^\prime \ni (u_1^\prime, u_2^\prime)  \mapsto \chi^{\lambda(u_1^\prime)}(\tilde{\Phi}(u_2^\prime)) \in  k^\times$ is symmetric.  

For any $i, j \in \{ 1,2 \}$, we fix $\sqrt{\alpha_{ij}} \in k^\times$ such that $\sqrt{\alpha_{ij}} = \sqrt{\alpha_{ji}}$. Then  
the $2 \times 2$ matrix $\begin{pmatrix} \sqrt{\alpha_{ij}} \end{pmatrix}_{i,j}$ with entries in $k^\times$ 
is symmetric. We define a symmetric $\ZZ$-bilinear form $b\colon M^\prime \times M^\prime \to k^\times$ by $b(\mathsf{f}_i^\prime, \mathsf{f}_j^\prime) = \sqrt{\alpha_{ij}}$ for 
any $i, j \in \{1, 2\}$. We note that $b(u')^2 = t(u', \lambda (u'))$ for any $u' \in M'$, where $t$ is defined in (\ref{eqn:def:t}). We define $c\colon M^\prime \to k^\times$ by $c(u^\prime) = b(u^\prime, u^\prime)$ for $u^\prime \in M^\prime$. Then we have 
$c(u_1^\prime + u_2^\prime)c(u_1^\prime)^{-1}c(u_2^\prime)^{-1} 
= b(u_1^\prime, u_2^\prime)^2 
=  \chi^{\lambda(u_2^\prime)}(\tilde{\Phi}(u_1^\prime))$. 
Further, $c_{\trop}(u^\prime) = -\log \left|b(u^\prime, u^\prime)\right| = \frac{1}{2} \langle \lambda(u^\prime), u^\prime\rangle$. 
This completes the proof. 
\QED

Let $(X, Q)$ be an irreducible principally polarized tropical abelian surface. As before, we use the notation $M$ and $X = N_\RR/M^\prime$, where $N_{\RR} = \Hom (M,\RR)$. Let $\lambda\colon M^\prime \to M$ be the $\ZZ$-linear map corresponding to $Q$. Fix a numbering of the elements of $M/2\lambda (M')$. Then we have a map $\psi\colon X \to \TT\PP^3$ given by (\ref{eqn:def:psi}) (cf. Remark~\ref{remark:well-defined:varphi:psi}), and the image $\psi (X)$ is called a tropical quartic Kummer surface (cf. Definition~\ref{definition:TKQS}).

\begin{Theorem} \label{thm:FTmain2}
Under the above setting, there exists a Kummer quartic surface 
$\tilde{Z} \subseteq \PP^3$ over some algebraically closed, nontrivially valued,  and complete nonarchimdean  field such that $\psi(Y) \subseteq \trop(\tilde{Z}^{\an})$, where $\trop$ is the map in \eqref{eqn:trop:morphism} for $n=3$. 
\end{Theorem}

\Proof
By Proposition~\ref{prop:lifting:abel}, there exist: a 
totally degenerate abelian surface $A$ over some algebraically closed, nontrivially valued, and complete nonarchimedean field; 
and 
a nonarchimedean descent datum $(\lambda , c)$ on $A^{\an}$ such that $c_{\trop}(u^\prime) = \frac{1}{2} \langle 
\lambda(u^\prime), u^\prime\rangle$ for any $u^\prime \in M^\prime$.
We take a basis $s_0, \ldots, s_3 \in H^0(A, \tilde{L}_1^{\otimes 2})$ 
as in the proof of Theorem~\ref{thm:main1:FT} and set
$\tilde{Z} \colonequals \tilde{\varphi}_{s_0,\ldots ,s_3}(\tilde{A})$. 
Then $\tilde{Z}$ is a Kummer quartic surface. 
Since the diagram 
\[
\begin{CD}
A^{\an} @>\tilde{\varphi}_{s_0,\ldots , s_3}>> (\PP^3)^{\an} \\
@AA\sigma_{A^{\an}}A @VV{\trop}V \\
X  @>
{\psi}
>> \TT\PP^3
\end{CD}
\]
commutes, we obtain $\psi(Y) \subseteq \trop(\tilde{Z}^{\an})$. 
\QED

\section{Canonical skeletons and Kontsevich--Soibelman skeletons}
\label{can:skeleton:essential:skeleton}

In this section, we remark on the relation between the canonical skeleton and the Kontsevich--Soibelman skeleton of a Kummer surface and prove a theorem on the faithful tropicalization of the Kontsevich--Soibelman skeletons.

\subsection{Uniformization and the canonical skeletons over algebraically non-closed base fields} \label{subsection:non-closed}

In this subsection, we remark on analytic geometry when the base field is not necessarily algebraically closed.

Let $k_0$ be a valuation field that is complete with respect the valuation and assume that the valuation is nontrivial. 
For any algebraic variety $X_0$ over $k_0$, we have the same description of $X_0^{\an}$ as the case where $k_0$ is algebraically closed, and thus we have an analytic theory in the sense of Berkovich.

Let $A_0$ be an abelian variety.
We say that $A_0$ is \emph{uniformized by a split torus} if there exists a morphism $p_0\colon  T_0^{\an} \to A_0^{\an}$ with $T_0$ being a split torus over $k_0$ such that $p_0$ is the universal cover. We call $p_0$ the uniformization of $A_0^{\an}$. Let $k$ be an algebraically closed complete valuation field that extends $k_0$. Note that there exists such a field and it is an algebraically closed, nontrivially valued, and complete nonarchimedean field. 
Suppose that the abelian variety $A_0 \otimes_{k_0} k$ over $k$ is totally degenerate. Then it is known that replacing $k_0$ by a finite separable extension if necessary, $A_0$ is uniformized by a split torus.

Assume that $A_0^{\an}$ is uniformized by the a split torus and let $p_0\colon T_0^{\an} \to A_0^{\an}$ denote the uniformization, where $T_0 = \Spec (k_0 [\chi^M])$ is a split torus over $k_0$ with character lattice $M$. In the same way as the case where the base field is algebraically closed, 
the valuation map $\val\colon T_0^{\an} \to N_{\RR}$ exists, $\Ker (p_0) \subseteq T(k)$, and the map $\val$ restricts to a group isomorphism $\Ker (p_0) \to M'$ for some full lattice $M'$ in $N_{\RR}$.
The canonical section $\sigma_0\colon N_{\RR} \to T^{\an}$ is defined, and the canonical skeleton $\Sigma_{T_0^{\an}}$ of $T^{\an}$ is defined to be the image of $\sigma_0$.
The map $\val\colon T^{\an} \to N_{\RR}$ descends to a valuation map $A_0^{\an} \to N_{\RR} / M'$, and the canonical section of $\val\colon T^{\an} \to N_{\RR}$ also descends to a section $N_{\RR} / M' \to A_0^{\an}$. The canonical skeleton $\Sigma_{A_0^{\an}}$ of $A_0^{\an}$ is defined to be the image of this section, which is a real torus with an integral structure. 

\begin{Remark} \label{remark:AH:non-closed}
The Appell--Humbert theory holds for line bundles on $A_0^{\an}$. Further, Lemma~\ref{lemma:basis-theta} holds even over $k_0$.
\end{Remark}

\begin{Remark} \label{remark:canonicalskeleton-basechage}
Let $k$ be an algebraically closed, nontrivially valued and complete nonarchimedean field that extends the valued field $k_0$. Set $A := A_0 \otimes_{k_0} k$. It follows from the construction of the canonical skeletons as above that the canonical morphism $A^{\an} \to A_0^{\an}$ induces an isomorphism between the canonical skeletons.
\end{Remark}

\subsection{Faithful tropicalization of the Kontsevich--Soibelman skeletons of Kummer surfaces}

In this section, we set $k_0 : = \CC (\!(t)\!)$ and fix a $t$-adic absolute value $|\ndot|_0$ on it. This is a discrete valuation. 
We set 
$k_0^\circ := \{\alpha \in k_0 \mid |\alpha| \leq 1\}$, the valuation ring of $k_0$.

For a smooth and proper algebraic variety $Z_0$ over $k_0$, Musta\c{t}\u{a} and Nicaise define in \cite{MN} the essential skeleton 
$\Sk(Z_0) \subseteq Z_0^{\an}$. 
When $Z_0$ is a $K3$ surface, the essential skeleton $\Sk(Z_0)$ coincides with the Kontsevich--Soibelman skeleton $\Sk(Z_0 , \omega)$, where $\omega$ is a nonzero regular $2$-form on $Z_0$; see \cite{KS} and 4.7.1, 4.7.3, Theorem~4.5.7, and Definition~4.10 in \cite{MN} for details. In this paper, when $Z_0$ is a $K3$-surface, we call the essential $\Sk(Z_0)$ the \emph{Kontsevich--Soibelman skeleton}.

Let $A_0$ be an abelian surface over $k_0$ throughout this subsection. 
Let $q_0\colon A_0 \to \tilde{Y}_0$ be the quotient by the involution group $\langle -1 \rangle$. We call $\tilde{Y}_0$ the Kummer surface associated to $A_0$. Let $\tilde{Y}_0' \to \tilde{Y}_0$ be the minimal resolution. Then $\tilde{Y}_0'$ is a smooth $K3$ surface. 
For an algebraic variety $X_0$ over $k_0$, we set $X_0^{\rm bir} := \{ (p , |\ndot|) \in X_0^{\an}  \mid \text{$p$ is the generic point of $X$} \}$. 
By definition (\cite[Definition~4.10]{MN}) that $\Sk (\tilde{Y}_0') \subseteq (\tilde{Y}_0')^{\rm bir}$. Since the continuous map $(\tilde{Y}_0')^{\an} \to (\tilde{Y}_0)^{\an}$ restricts to a homeomorphism $(\tilde{Y}_0')^{\rm bir} \to (\tilde{Y}_0)^{\rm bir}$, we have a homeomorphism $\Sk(\tilde{Y}_0') \to \Sk(\tilde{Y}_0)$. We call $\Sk(\tilde{Y}_0)$ the \emph{Kontsevich--Soibelman skeleton} of~$\tilde{Y}_0^{\an}$. 

A \emph{strict Kulikov model} of $\tilde{Y}_0'$ is a regular strictly semistable projective model $\tilde{\mathcal{Y}}_0' \to \Spec(k_0^\circ)$ of $\tilde{Y}_0'$ with trivial relative dualizing sheaf; see \cite{Overkamp} for more details of  strict Kulikov models. 
 With this notation, we have $\Sk({\tilde{Y}_0}) = S_{\tilde{Y}_0^{\an}} (\tilde{\mathcal{Y}}_0')$, where $\tilde{\mathcal{Y}}_0'$ is a strict Kulikov model of $\tilde{Y}_0'$.

Note that there exists an algebraically closed field extension $k$ of $k_0$ equipped with an absolute value $|\ndot|_k$ extending $|\ndot|_0$ such that $k$ is complete with respect to $|\ndot|_k$. We fix such a $k$ in the sequel. We set $A := A_0 \otimes_{k_0} k$ and $\tilde{Y} := \tilde{Y}_0 \otimes_{k_0} k$.
The following proposition shows that our canonical skeleton $\Sigma_{\tilde{Y}^{\an}}$ 
of $\tilde{Y}^{\an}$ in Definition~\ref{def:can:skeleton:Kummer} coincides with the Kontsevich--Soibelman skeleton of $\tilde{Y}_0^{\an}$  under a suitable assumption if we ignore difference of the base fields.

\begin{Proposition} 
\label{thm:canonical-essential}
With the above notation, assume that $A$ is totally degenerate. Let $\nu^{\an}\colon \tilde{Y}^{\an} \to \tilde{Y}_0^{\an}$ be the natural map. Then the following holds.
\begin{enumerate}
\item
The map $\nu^{\an}\colon \tilde{Y}^{\an} \to \tilde{Y}_0^{\an}$ restricts to a surjective continuous map $\Sigma_{\tilde{Y}^{\an}} \to \Sk({\tilde{Y}_0})$. 
\item
Suppose that $A_0^{\an}$ is uniformized by a split torus over $k_0$. Then  
the map $\Sigma_{\tilde{Y}^{\an}} \to \Sk({\tilde{Y}_0})$ is a homeomorphism.
\end{enumerate}
\end{Proposition}

\Proof
There exists a finite extension $k_1$ of $k_0$ in $k$ such that $(A_0 \otimes_{k_0} k_1)^{\an}$ is uniformized by a split torus over $k_1$. 
By \cite[Proposition~3.1.3]{HN},  the canonical morphism $\tilde{A}_0 \otimes_{k_0} k_1 \to \tilde{A}_0$ induces a surjective continuous map $\Sk (\tilde{A}_0 \otimes_{k_0} k_1) \to \Sk (\tilde{A}_0)$. It follows that it suffices to prove (2).

We assume that $A_0^{\an}$ is uniformized by a split torus over $k_0$. There exists a finite extension $k_1$ of $k_0$ such that the $2$-torsion points of $A$ lie in $A_0 (k_1)$. Further, replacing $k_1$ by a finite extension if necessary, we may and do assume,  by K\"unnemann's construction in \cite{Ku}, that there exists a regular strictly semistable projective model $\mathcal{A}_1 \to \Spec(k_1^{\circ})$ of $A_1 := A_0 \otimes_{k_0} k_1$ with the following properties:
\begin{enumerate}
\item[(i)]
the skeleton $S(\mathcal{A}_1)$ equals the canonical skeleton $\Sigma_{A_1^{\an}}$ of $A_1^{\an}$ (see Section~\ref{subsection:non-closed} for the canonical skeleton over $k_1$, which is not algebraically closed); 
\item[(ii)]
the $\langle -1 \rangle$-action on $A_1$ extends to 
$\mathcal{A}_1$; 
\item[(iii)] 
if $\mathcal{A}_1'$ denote the blow-up of $\mathcal{A}_1$ along the closure of the $2$-torsion subgroup $A_1 [2]$, then the $\langle -1 \rangle$-action on $\mathcal{A}_1$ lifts to that on $\mathcal{A}_1'$. 
\end{enumerate}
See \cite[Proposition~3.8, its proof, and Example~3.5]{NXY} for property (i) in the above.
Then, as in \cite{Overkamp}, 
the quotient $\mathcal{A}_1' / \langle -1 \rangle$ is a strict Kulikov model of the $K3$ surface $\tilde{Y}_1'$, and the map $q_1^{\an}\colon A_1^{\an} \to \tilde{Y}_1^{\an}$ restricts to a map $S(\mathcal{A}_1) \to S_{\tilde{Y}_1^{\an}}(\mathcal{A}_1' / \langle -1 \rangle)$. By \cite[Theorem~4.8]{Overkamp}, this map induces a homeomorphism $S(\mathcal{A}_1) / \langle -1 \rangle \to \Sk(\tilde{Y}_1)$. Since $S(\mathcal{A}_1) = \Sigma_{A_1^{\an}}$, we thus have a homeomorphism
$\beta_1\colon \Sigma_{A_1^{\an}} / \langle -1 \rangle \to \Sk(\tilde{Y}_1)$
induced by the quotient $q_1\colon A_1 \to \tilde{Y}_1$.

Let $\Sigma_{A_0^{\an}} \subseteq A_0^{\an}$ be the canonical skeleton of $A_0^{\an}$. 
We have a commutative diagram
\[
\begin{CD}
\Sigma_{A^{\an}} @>>> \Sigma_{A^{\an}} / \langle -1 \rangle \\
@V{\cong}VV @V{\cong}VV \\
\Sigma_{A_1^{\an}} @>>> \Sigma_{A_1^{\an}} / \langle -1 \rangle @>{\beta_1}>{\cong}> \Sk(\tilde{Y}_1) \\
@V{\cong}VV @V{\cong}VV @VVV \\
\Sigma_{A_0^{\an}} @>{\alpha_0}>> \Sigma_{A_0^{\an}} / \langle -1 \rangle @>{\beta_0}>> \tilde{Y}_0^{\an} ,
\end{CD}
\]
where the vertical arrows are maps induced by the canonical morphisms by basechanges, the horizontal arrows are induced by the quotient morphisms of the abelian surfaces by the $\langle -1 \rangle$-actions. Indeed, since $A_0^{\an}$, $A_1$, and $A$ are uniformized by split tori over $k_0$, $k_1$, and $k$, respectively, the maps $\Sigma_{A_1^{\an}} \to \Sigma_{A_0^{\an}}$, $\Sigma_{A^{\an}} \to \Sigma_{A_1^{\an}}$, and $\Sigma_{A^{\an}} \to \Sigma_{A_0^{\an}}$ induced by the canonical morphisms $A_1 \to A_0$, $A \to A_1$, and $A \to A_0$, respectively, are all isomorphisms of rational polyhedral spaces (cf. Remark~\ref{remark:canonicalskeleton-basechage}). Since the maps $\Sigma_{A_1^{\an}} \to \Sigma_{A_0^{\an}}$ and $\Sigma_{A^{\an}} \to \Sigma_{A_1^{\an}}$ are $\langle -1 \rangle$-equivariant, it follows that the induced maps $\Sigma_{A_1^{\an}} / \langle -1 \rangle \to \Sigma_{A_0^{\an}} / \langle -1 \rangle$ and $\Sigma_{A^{\an}} / \langle -1 \rangle \to \Sigma_{A_0^{\an}} / \langle -1 \rangle$ are isomorphisms of rational polyhedral orbifolds. 
As we remarked above, $\beta_1$ is a homeomorphism. Further, we see that the same statement as Lemma~\ref{lemme:finite:quotient:Berkovich} holds also over $k_0$. Since $\alpha_0$ is the quotient by the $\langle -1 \rangle$-action, it follows  that $\beta_0$ is just an inclusion of a topological space. By \cite[Proposition~3.1.3]{HN}, the image of the right vertical arrow equals the Kontsevich--Soibelman skeleton $\Sk (\tilde{Y}_0)$ of $\tilde{Y}_0$. By the commutativity of the diagram, it follows that $\beta_0$ induces an homeomorphism 
$\Sigma_{A_0^{\an}} / \langle -1 \rangle \to \Sk (\tilde{Y}_0)$.
Since $\Sigma_{A^{\an}} / \langle -1 \rangle = \Sigma_{A_0^{\an}} / \langle -1 \rangle$, this completes the proof. 
\QED

Finally, we state our faithful tropicalization result for the Kontsevich--Soibelman skeletons. We defined the notion of faithful tropicalization when the base valued field is algebraically closed, but the definition makes sense even when the base field is not algebraically closed.

\begin{Theorem}
Assume that the abelian surface $A_0^{\an}$ over $k_0$ is uniformized by a split torus over $k_0$. Let $\tilde{L}_1$ be a line bundle on $A_0$ that gives a principal polarization. Assume that $(A \otimes_{k_0} k , \tilde{L}_1 \otimes_{k_0} k )$ is tropically irreducible. We set $\tilde{L} \colonequals \tilde{L}_1^{\otimes 2}$ and let $(\tilde{L})_{\tilde{Y}_0}$ be the descent to $\tilde{Y}_0$. Then there exist a basis $s_0 , \ldots , s_3 \in H^0 (\tilde{Y}_0 , (\tilde{L})_{\tilde{Y}_0})$ such that the induced morphism $\varphi_{s_0 , \ldots , s_3}$ faithfully tropicalizes the Kontsevich--Soibelman skeleton of $\tilde{Y}_0^{\an}$.
\end{Theorem}

\Proof
We take a basis $s_0 , \ldots , s_3 \in H^0 (\tilde{Y}_0 , (\tilde{L})_{\tilde{Y}_0})$ as in Lemma~\ref{lemma:basis-theta}; see also Remark~\ref{remark:AH:non-closed} as $k_0$ is not algebraically closed.
Let $(\tilde{L})_{\tilde{Y}}$ be the pullback of $(\tilde{L})_{\tilde{Y}_0}$ by the canonical morphism $\tilde{Y} \to \tilde{Y}_0$. Then we have $s_0' , \ldots , s_3' \in H^0 (\tilde{Y} , (\tilde{L})_{\tilde{Y}})$, where $s_i'$, for $i=0, \ldots , 3$, denotes the image of $s_i'$ under the natural map $H^0 (\tilde{Y}_0 , (\tilde{L})_{\tilde{Y}_0}) \hookrightarrow H^0 (\tilde{Y} , (\tilde{L})_{\tilde{Y}})$. Note that those global sections differ from the global sections that we use in the proof of Theorem~\ref{thm:main1:FT} only by nonzero constant-multiples, it follows from this theorem that the morphism $(\varphi_{s_0 , \ldots , s_3})_k\colon \tilde{Y} \to \PP_k^3$ induced by the global sections $s_0 , \ldots , s_3$ over $\tilde{Y}$ faithfully tropicalizes the canonical skeleton $\Sigma_{Y^{\an}}$.

Note that we have a natural commutative diagram
\[
\begin{CD}
\tilde{Y}^{\an} @>{(\varphi_{s_0 , \ldots , s_3})_k^{\trop}}>> \TT\PP^m \\
@V{\nu}VV @VV{\id}V \\
\tilde{Y}_0^{\an} @>{\varphi_{s_0 , \ldots , s_3}^{\trop}}>> \TT\PP^m.
\end{CD}
\]
By Theorem~\ref{thm:canonical-essential}, $\nu$ restricts to an isomorphism $\Sigma_{\tilde{Y}^{\an}} \to \Sk (\tilde{Y}_0)$ by Theorem~\ref{thm:canonical-essential}. 
Thus the theorem holds.
\QED

\end{document}